\documentclass[a4paper,11pt,reqno]{amsart}

\usepackage{lmodern}

\usepackage[utf8]{inputenc} 
\usepackage[T1]{fontenc}

\usepackage[english]{babel}

\usepackage{amsfonts}
\usepackage{amssymb}
\usepackage{amscd}
\usepackage{amsthm}
\usepackage{a4wide}
\usepackage{mathrsfs}
\usepackage{amsmath}
\usepackage{amssymb}
\usepackage{amsthm}
\usepackage{textcomp}
\usepackage{graphicx}
\usepackage{enumerate}
\usepackage{mathrsfs}
\usepackage{frcursive}
\usepackage{tikz}
\usepackage[cyr]{aeguill}
\usepackage{xspace}
\usepackage{hyperref}
\usepackage{appendix}
\usepackage{esint}
\usepackage{cancel}
\usepackage{calc}
\usepackage[colorinlistoftodos,prependcaption]{todonotes}

\newtheorem{defn}{Definition}[section]
\newtheorem{lemma}[defn]{Lemma}
\newtheorem{prop}[defn]{Proposition}
\newtheorem{theo}[defn]{Theorem}
\newtheorem{coro}[defn]{Corollary}
\newtheorem{conj}[defn]{Conjecture}
\newtheorem{claim}{Claim}
\newtheorem{rk}[defn]{Remark}
\newtheorem{ques}[defn]{Question}

\newcounter{mtheorem}
\newtheorem{mtheorem}[mtheorem]{Theorem}

\def\RR{\mathbb{R}}
\def\Ric{\mathop{\rm Ric}\nolimits}

\def\Rm{\mathop{\rm Rm}\nolimits}
\def\tr{\mathop{\rm tr}\nolimits}

\def\vol{\mathop{\rm vol}\nolimits}

\def\dim{\mathop{\rm dim}\nolimits}
\def\vol{\mathop{\rm Vol}\nolimits}

\def\inj{\mathop{\rm inj}\nolimits}

\def\div{\mathop{\rm div}\nolimits}

\def\Id{\mathop{\rm Id}\nolimits}

\def\Li{\mathop{\rm \mathscr{L}}\nolimits}

\def\Ric{\mathop{\rm Ric}\nolimits}

\def\Rm{\mathop{\rm Rm}\nolimits}
\def\tr{\mathop{\rm tr}\nolimits}

\def\vol{\mathop{\rm vol}\nolimits}

\def\dim{\mathop{\rm dim}\nolimits}
\def\vol{\mathop{\rm Vol}\nolimits}

\def\inj{\mathop{\rm inj}\nolimits}

\def\div{\mathop{\rm div}\nolimits}

\def\Id{\mathop{\rm Id}\nolimits}

\def\Li{\mathop{\rm \mathscr{L}}\nolimits}

\def\R{\mathop{\rm R}\nolimits}

\def\cprime{$'$}
\def\cprime{$^$}

\newsavebox\CBox
\newcommand\hcancel[2][0.5pt]{%
  \ifmmode\sbox\CBox{$#2$}\else\sbox\CBox{#2}\fi%
  \makebox[0pt][l]{\usebox\CBox}%  
  \rule[0.5\ht\CBox-#1/2]{\wd\CBox}{#1}}

\numberwithin{equation}{section}

\begin{document}
\title{Non self-similar solutions to the Ricci flow coming out of metric cones}
\date{}
\author{Alix Deruelle}
\address{Universit\'e Paris-Saclay, CNRS, Laboratoire de math\'ematiques d'Orsay, 91405, Orsay, France}
\email{alix.deruelle@universite-paris-saclay.fr}
\author{Tristan Ozuch}
\address{MIT, Dept. of Math., 77 Massachusetts Avenue, Cambridge, MA 02139-4307.}
\email{ozuch@mit.edu}
\author{Felix Schulze}
\address{Mathematics Institute, Warwick University, Gibbet hill road, Coventry CV4 7AL, UK}
\email{felix.schulze@warwick.ac.uk}
\maketitle

\begin{abstract}
    Given an asymptotically conical \textit{unstable} expanding gradient Ricci soliton, we provide the existence of a family of Ricci flows starting at the same metric cone up to isometry. Moreover, the family of solutions we construct is shown to be parametrized by positive modes of the weighted Lichnerowicz operator associated to the background expanding gradient Ricci soliton. Each solution of this family is shown not to be induced by an expanding \textit{gradient} Ricci soliton. If the background expander has no non-trivial infinitesimal deformations, we show furthermore that the solutions we provide are not self-similar. We finally prove that solutions are distinct up to isometries preserving the potential function of the background expanding gradient Ricci soliton.
\end{abstract}

\maketitle

\setcounter{tocdepth}{1}
\markboth{Alix Deruelle, Tristan Ozuch and Felix Schulze}{Non self-similar solutions to the Ricci flow}

\tableofcontents
\section{Introduction}

\subsection{Background}
%Stability of a fixed point of a continuous dynamical system is a key property as it governs the flow of nearby small perturbations. In an infinite dimensional setting, this is particularly striking for solutions of geometric partial differential equations associated to a (geometric) flow such as the mean curvature flow or the Ricci flow. In this paper, we focus on forward self-similar solutions to the Ricci flow. 

%Let us make precise the notion of self-similar solutions we use here.
%\todo[inline]{cite Kotschwar's result on solutions terminating on a cone \cite{Kot-Con}}
Recall that an expanding gradient Ricci soliton is a triple $(M^n,g,\nabla^gf)$ where $(M^n,g)$ is a complete Riemannian manifold and $f:M\rightarrow \mathbb{R}$ is a smooth function (called the potential function of the soliton) satisfying the following Bakry-\'Emery condition on $M$:
\begin{equation}\label{egs-eqn}
\Ric(g)-\nabla^{g,\,2}f=-\frac{g}{2}.
\end{equation}
If $(\varphi_t)_{t>0}$ denotes the flow generated by $-t^{-1}\nabla^gf$ such that $\varphi_t\big|_{t=1}=\Id_M$, then $g(t):=t\varphi_t^*g$ is a self-similar solution to the Ricci flow for $t>0$ which evolves only by diffeomorphisms and scaling. Such forward self-similar solutions have an initial condition at $t=0$ that may be singular both at the topological and metric level. We restrict here to singular initial conditions that are metric cones which are smooth outside the tip. Examples of such solutions have been constructed using ODE techniques (see for example \cite{Buz-Dan-Gal-Wan, Ang-Kno, Nie-Win}), via blow-down methods \cite{Sch-Sim}, continuity arguments \cite{Der-Smo-Pos-Cur-Con, con-der} as well as more recently using degree theoretic arguments \cite{Bam-Che}. In sharp contrast to the setting of shrinking gradient Ricci solitons, where Kotschwar \cite{Kot-Con} has shown that any Ricci flow terminating smoothly in a cone (with quadratic curvature decay) must be a gradient shrinking soliton, we exhibit an abundance of non self-similar solutions emerging from the same cone forward in time whenever the background soliton is linearly unstable.

One of the main motivations to study such solutions is that they appear naturally as first order models smoothing out conical singularities. These arise for example when considering Ricci flows from non-smooth initial data but are also natural in the aim to understand and construct Ricci flow through singularities. See \cite{FIK} for an example of an asymptotically conical self-similarly shrinking solution which limits at time zero to a metric cone (smooth outside of the tip) and which is resolved by a smooth self-similarly expanding solution for positive time.  In \cite{Sto-Sgs} it is shown that any asymptotically conical self-similarly shrinking solution arises as a singularity model at a finite time singularity for a Ricci flow on a compact manifold. Even more, as part of the aim to establish the existence of a Ricci flow through singularities in dimension four it is expected that such singularities arise generically, see \cite{Naff-Ozuch1} and \cite{Kim-Ozuch2}. A further important aspect where expanding gradient solitons appear is the large time asymptotics of immortal solutions, see for example \cite{Chen-Zhu, Lott-Ricci-pinched, Lee-Topping1, Deruelle-Schulze-Simon-HamConj, Lee-Topping2}.

A property which is expected to be central in establishing the existence and uniqueness of a Ricci flow from initial data with conical singularities is the linear stability of the self-similarly expanding solitons modeling the resolution of the singularities. In this article we show that \emph{linear instability} of an asymptotically conical gradient soliton leads to the existence of non self-similar solutions coming out of the initial metric cone, which in turn implies non-uniqueness. Section \ref{sec:unstable-examples} provides an example of an unstable expanding soliton. 

Furthermore, as pointed out to us by E.~Chen, their work with Bamler \cite{Bam-Che} yields that in some cases non-uniqueness for the existence problem also implies the existence of a linearly unstable expander as follows: Consider an asymptotic cone $C$ which is a regular value of the asymptotic projection from the moduli space of expanders (which is expected to be a generic property) and two distinct expanding solitons $(X^4,g_i)_{i=1}^2$ asymptotic to $C$, where $X^4$ has expander degree one. Then the parity identity \cite[Eq.~(1.10)]{Bam-Che} implies that the Morse indices of the two expanders differ by an odd number, hence at least one has to be unstable. Similarly, Angenent and Knopf \cite{Ang-Kno} have shown that from dimension at least five on there exist cones $C$, smooth outside the tip, with an arbitrary number of non-isometric expanding gradient Ricci solitons coming out of $C$. It is expected that these expanding solitons are linearly unstable \cite{Ang-pers-com}.

\subsection{Main result}

Given an expanding gradient Ricci soliton $(M^n,g,\nabla^gf)$, the weighted Lichnerowicz operator denoted by $L_f$ acting on symmetric $2$-tensors is defined by:
\begin{equation*}
\begin{split}
L_fh:=\Delta_gh+\nabla^g_{\nabla^gf}h+2\Rm(g)\ast h,\quad (\Rm(g)\ast h)_{ij}:=\Rm(g)_{ipqj}h^{pq},\quad h\in S^2T^*M.
\end{split}
\end{equation*}
Here $\Delta_g$ denotes the rough Laplacian on tensors induced by the metric $g$.
It is straightforward to check that $L_f$ is symmetric with respect to the $L^2$-weighted scalar product induced by the weighted measure $e^fd\mu_g$, which we denote by $L^2_f(S^2T^*M)$. Also noteworthy is the fact that under mild assumptions on the curvature, the $L^2_f(S^2T^*M)$-spectrum of $L_f$ is discrete.

We say that an expanding gradient Ricci soliton is (linearly) \textit{stable} provided $L_f$ is a nonpositive operator on $L^2_f(S^2T^*M)=:L^2_f$.

Our main result is the following existence of an unstable manifold whenever an expanding gradient Ricci soliton is not stable:
\begin{mtheorem}\label{theo-main}
Let $(M^n,g_0,\nabla^{g_0}f_0)$, $n\geq 3$, be an asymptotically conical expanding gradient Ricci soliton with tangent cone at infinity $(C_0,d_c)$. Denote by $E_+$ the finite dimensional vector space of eigentensors associated to positive eigenvalues of the weighted Lichnerowicz operator $L_{f_0}$ with corresponding orthogonal projection $\pi_+$. Assume $\dim E_+>0$. 

Then there exists $\varepsilon=\varepsilon(n,g_0)>0$ such that the following assertions hold:
\begin{enumerate}
\item \label{item-exi} For each $h_0\in E_+$ with $\|h_0\|_{L^2_{f_0}}<\varepsilon$, there exists  a solution $(g(t))_{t\in(0,1]}$ to the Ricci flow such that $g(1)=g_0+h(0)$ with $\pi_+(h(0))=h_0$ and which smoothly converges to $(C_0,d_c)$ outside the tip as $t$ tends to $0$, up to isometry. 

\item  \label{item-grad} If $(g(t))_{t\in(0,1]}$ induces an expanding \textbf{gradient} Ricci soliton $(M^n,g_1,X_1=\nabla^{g_1}f_1)$ with $g_1=g_0+h(0)$, $\pi_+(h(0))=h_0$ and $\|h_0\|_{L^2_{f_0}}<\varepsilon$ then $h(0)\equiv0$, i.e. $g_1=g_0$ and $X_1=\nabla^{g_0}f_0$. 

 \item \label{item-non-grad} If $(M^n,g_0,\nabla^{g_0}f_0)$ satisfies $\ker L_{f_0}\cap L^2_{f_0}=\{0\}$ then if $(g(t))_{t\in(0,1]}$ induces an expanding Ricci soliton $(M^n,g_1,X_1)$ with $g_1=g_0+h(0)$, $\pi_+(h(0))=h_0$, $\|h_0\|_{L^2_{f_0}}<\varepsilon$ then $h(0)\equiv0$, i.e. $g_1=g_0$ and $X_1-\nabla^{g_0}f_0$ is a Killing field orthogonal to $\nabla^{g_0}f_0$. 
 
\item \label{item-mod-spa} For $h_1$ and $h_2$ in $E_+$ with $\|h_i\|_{L^2_{f_0}}<\varepsilon$, $i=1,2$, let $(g_i(t))_{t\in(0,1]}$ be the solutions to the Ricci flow associated to each $h_i$ as in \eqref{item-exi}. Then there exists a smooth diffeomorphism $\varphi$ of $M$ preserving the zero set of $\nabla^{g_0}f_0$ such that $g_2(t)=\varphi^*g_1(t)$ for every $t\in(0,1]$ if and only if $\varphi$ is an isometry of $g_0$ (i.e. $\varphi^*g_0=g_0$) which preserves the potential function $f_0$ (i.e. $\varphi^*f_0=f_0$) and $h_2=\varphi^*h_1$.
\end{enumerate}

Moreover, if $(C_0,d_c)$ has positive scalar curvature then there exists a neighborhood $U_g$ of $g$ in a suitable topology such that each solution $(g(t))_{t\in(0,1]}$ in $U_g$ has positive scalar curvature too.
%\todo[inline]{F: I guess we are using here the recent result of Cao et al that the expander has positive scalar curvature?}

\end{mtheorem}
% Given $\chi_1>0$, there exists $\varepsilon=\varepsilon(n,g_0,\chi_1)>0$ such that if $(g(t))_{t\in(0,1]}$ induces an expanding Ricci soliton $(M^n,g_1,X_1)$ with $g_1=g_0+h(0)$, $\pi_+(h(0))=h_0$, $\|h_0\|_{L^2_{f_0}}<\varepsilon$ and such that $|X_1|_{g_0}\leq \chi_1f_0^{1/2}$ then $h(0)\equiv0$, i.e. $g_1=g_0$ and $X_1-\nabla^{g_0}f_0$ is a Killing field orthogonal to $\nabla^{g_0}f_0$. Moreover, if $X_1$ is gradient, $\varepsilon=\varepsilon(n,g_0)$ only.
Recall that a non self-similar solution is a solution that cannot be induced by an expanding Ricci soliton $(M^n,g',X')$ whether $X'$ is gradient or not.
A more formal statement can be found in Theorem \ref{theo-ancient-sol} regarding the existence part of Theorem \ref{theo-main}. The convergence to the cone of the solutions provided by Theorem \ref{theo-ancient-sol} will be made precise in Section \ref{sec-geo-sol}.

  In Theorem \ref{theo-main}, we present one potential source of nonuniqueness: the linear instability of the expanding soliton asymptotic to cones. Indeed, such expanders may be seen as Ricci flows coming out of their tangent cone at (spatial) infinity, and we show that when the expander in question is unstable, a plethora of solutions to the Ricci flow exist while coming out of the same metric cone (up to isometry). As these Ricci flows are not self-similar, this is a genuine new type of flow out of the cone. The proof of the non self-similarity is subtle in the sense that one needs to get rid of the action of symmetries induced by Killing vector fields of the initial metric cone,  a phenomenon that is specific to Ricci flow and which does not seem to occur for other geometric flows such as the mean curvature flow. 

As recalled in Section \ref{sec-geo-sol}, an expanding gradient Ricci soliton with nonnegative sectional curvature is strictly stable, i.e. the $L^2_f$-spectrum of $L_f$ is contained in $\RR^{-}$. Observe that  nonnegative sectional curvature is implied by the following curvature constraints: nonnegative curvature operator and weakly PIC2. Recently, Cao and Xie \cite{Cao-Xie} have proved that a three-dimensional expanding gradient Ricci soliton with positive scalar curvature must have positive sectional curvature. In particular, a three-dimensional expanding gradient Ricci soliton with positive scalar curvature is strictly stable. The strict stability of all (real) $4$-dimensional Kähler Ricci solitons has recently been established in \cite{naff2025linearstabilityinstabilitykahler}, the case of higher dimensional Kähler expanders have been settled even more recently in \cite{Lavoyer-Thompson}.

We have also recently learnt from Lavoyer \cite{Lav-pers-com} that a weakly PIC1 expanding gradient Ricci soliton is strictly stable: its proof is direct in the sense that he proves that the Rayleigh quotient associated to $-L_f$ is strictly positive. Observe that the combination of Theorem \ref{theo-main} and \cite{Chan-Peachey-Lee} implies by a contradiction argument that an asymptotically conical weakly PIC1 expanding gradient Ricci soliton is stable. As this is a weaker statement than the aforementioned strict stability result of Lavoyer, we refrain ourselves from stating it as a corollary but we believe that proving \textit{linear} stability of a Ricci soliton through dynamics of the flow is worth pointing out.

%\todo[inline]{TO: with Keaton, we'll soon post the strict stability of Kähler expanding solitons in 4D as well. A: ok, you can cite this (upcoming) work too. Is it reasonable to think that this should hold in any dimension (similar to Kähler-Einstein metrics on closed manifolds as in Dai-Ma-Wei?)  TO: Not sure, honestly, I don't think that parallel spin$^c$ spinors play super well with the \textit{weighted} Lichnerowicz Laplacian, but who knows?}

As a consequence of the previous discussion, we believe that the following question is important for understanding the uniqueness of the Ricci flow through isolated singularities: 
\begin{ques}\label{ques-1}
Are weakly PIC expanding gradient Ricci solitons (strictly) stable?  
\end{ques}

To put Question \ref{ques-1} into perspective, let us remark that  we do not know whether there exists an unstable expanding gradient Ricci soliton in dimension at least four with positive scalar curvature.
We end this section by recalling (in a non exhaustive way) the rapidly growing literature on the links between linear (in)stability and dynamical (in)stability around a fixed point of the Ricci flow.
The articles \cite{Ses-Ric-Fla} and \cite{Has-Mul} investigated dynamical stability near Ricci flat metrics on a closed manifold and \cite{Kro-Sta-Ric-Sol} systematically investigated Ricci solitons on a closed manifold. In the non-compact setting, stability of Euclidean space was established in \cite{Sch-Sch-Sim-Euc} and \cite{App}. Dynamical (in)stability of Asymptotically Locally Euclidean Ricci flat metrics are the subject of \cite{Kro-Pet} in the integrable case and \cite{Der-Ozu-Sta} in dimension $n\geq 5$. The articles \cite{Sch-Sch-Sim-Hyp} and \cite{Bam-Sta-Sym-Spa} proved stability of Hyperbolic space (and that of symmetric spaces). Dynamical stability of non-compact strictly stable expanding gradient Ricci solitons was studied in \cite{Der-Sta-Egs} and can be viewed as a companion of Theorem \ref{theo-main}. 

Finally, uniqueness of the solutions produced in Theorem \ref{theo-main} is unknown and the work \cite{Choi-Lai} on uniqueness of ancient Ricci flows on a closed manifold with tangent flow at $-\infty$ an integrable gradient shrinking Ricci soliton should be key here.
%\todo[inline]{A:  Which other ref am I missing here? Write them in the todo and I'll add them}

\subsection{Outline of paper}
We begin in Section \ref{sec-prelim} by recalling the basics about asymptotically conical expanding gradient Ricci solitons. Prerequisites on the spectral properties of the linearized operator associated to such a soliton are then explained in Section \ref{sec-spec-l-f}. Section \ref{sec-set-up-eqn} states the equivalence between the (DeTurck)-Ricci flow and the normalized (DeTurck-)Ricci flow in this setting. The relevant function spaces that are required to prove the existence part of Theorem \ref{theo-main} are defined in Section \ref{sec-defn-fct-spa}.

Section \ref{sec-exi-uns-mfd} is devoted to proving the ``unstable theorem'' aspect of Theorem \ref{theo-main}: this is the main content of Theorem \ref{theo-ancient-sol}. The proof of the latter comprises energy estimates established in Section \ref{sec-ene-est} and pointwise estimates detailed in Section \ref{sec-pt-est}. The aforementioned estimates are then combined to give a proof of Theorem \ref{theo-ancient-sol} in Section \ref{sec-pf-thm}. Let us observe that the proof of Theorem \ref{theo-ancient-sol} is rooted in the original proof of the existence of an unstable manifold for ordinary differential equations. In the setting of geometric flows, we follow the approach of \cite[Theorem $1.6$]{Cho-Man} which considers the mean curvature flow near \textit{closed} unstable minimal submanifolds of a fixed ambient Riemannian manifold. However, because of the non-compactness of the underlying manifold here, the choice of function spaces is delicate and requires some care: the density $e^f$ grows exponentially quadratically fast at infinity which, for instance, prevents the embedding of polynomially decaying tensors in a Hilbert space modeled on $L^2_f(S^2T^*M)$.

A similar issue has been recently treated in \cite{Kno-Ang-Ric-Fla} in the setting of asymptotically flat Ricci flat metrics. However, their approach is different from ours in that the linearized operator at stake has a continuous $L^2$-spectrum equal to the half-line $(-\infty,0]$ with no embedded eigenvalues in it. 

The short Section \ref{sec-dec-spa} is devoted to improving the decay both at spatial infinity  and in the ``past'' of solutions given by Theorem \ref{theo-ancient-sol}. The result is sharp in time and almost sharp in space. The estimate in Proposition \ref{prop-alm-sharp-rate} we provide is also quantitative in terms of an upper bound on the (weighted $L^2$ norm of an) initial unstable perturbation and its proof makes use of a suitable Gagliardo-Nirenberg inequality whose existence is justified in Section \ref{sec-gag-nir-ine}. 

The content of Section \ref{sec-geo-sol} takes on a more geometric flavour starting with curvature estimates for expanding gradient Ricci solitons in Section \ref{sec-prelim-egs-curv-con}. The latter serve as a preliminary step to prove Proposition \ref{prop-curv-con} which shows that solutions produced by Theorem \ref{theo-ancient-sol} inherits the same curvature constraint as those satisfied by the background unstable expanding gradient Ricci soliton. This section ends by carefully analyzing the convergence of the solutions from Theorem \ref{theo-ancient-sol} to the normalized DeTurck-Ricci flow, respectively the DeTurck-Ricci flow and the Ricci flow: see Section \ref{sec-conv-init-time}. 

Section \ref{sec-non-self-sim-sol} proves the non self-similarity of the non trivial solutions of Theorem \ref{theo-ancient-sol} under the assumption that $\ker L_{f_0}\cap L^2_{f_0}=\{0\}$. The proof goes by contradiction and assumes the existence of two expanding Ricci solitons $(M^n,g_i,X_i)_{i=0,1}$ where $X_0$ is gradient and where $(M^n,g_0,X_0)$ is unstable. Recall that the assumption is that  the difference of the metrics $g_1-g_0=:h_0$ ha a non-trivial projection on $E_+$, the finite dimensional vector space of eigentensors of the linearized operator $L_{f_0}$ associated to positive eigenvalues. 

One of the difficulties lies in the fact that the soliton vector field $X_1$ is not assumed to be gradient: even if $h_0\equiv 0$, the soliton vector fields $X_i$, $i=0,1$ could still differ by a non trivial Killing vector field. This explains the necessity of subtracting a canonical extension of a given Killing vector field orthogonal to the radial direction on the cone, a fact justified by Sections \ref{sec-kill-metric-cone} and \ref{sec-kill-ext}. An additional technical but fundamental issue is that the system of equations we are dealing with contains a non local term given by the difference $X:=X_1-X_0$: it is shown that the vector field $X$ satisfies an elliptic equation as demonstrated in Lemma \ref{lem-bochner-lie-der} which can be coupled with the (difference of) the soliton equations \eqref{egs-eqn} together with a further elliptic equation satisfied by the weighted divergence of $X$ (see Lemma \ref{lem-equiv-gauge-eqn}) to result in a system of three (degenerate) elliptic equations. 

Similar but more involved issues enter the proof of Theorem \ref{theo-main}, item \eqref{item-mod-spa} which explains the length of Section \ref{sec-item-mod-spa}: there, one  requires linearizing one solution around another in a parabolic way instead of focusing at the elliptic equation obtained by linearizing the soliton equation as explained above.

The final Section \ref{sec:unstable-examples} legitimizes Theorem \ref{theo-main} by describing an unstable asymptotically conical expanding gradient Ricci soliton in dimension $4$.

\subsection{Acknowledgements}
The authors would like to thank Ronan Conlon for suggesting A. Rajan's article on cohomogeneity one expanders as a source of instability. AD acknowledges support from the Institut Universitaire de France and is partially supported by ``Fondation Charles Defforey'' together with the French Agence Nationale de la Recherche (ANR)  ANR-24-CE40-0702 (Project OrbiScaR). TO was partially supported by NSF Grant DMS-2405328. FS has received funding by the European Union (European Research Council (ERC), GENREG -- 101200301). Views and opinions expressed are however those of the authors only and do not necessarily reflect those of the European Union or the European Research Council. Neither the European Union nor the granting authority can be held responsible for them.

\subsection{AI disclosure} The argument in Appendix \ref{sec-lin-growth}, which removes the linear growth hypothesis on the soliton vector field from Theorem \ref{theo-main}, was suggested by a LLM. The authors checked, simplified and wrote the appendix themselves. No other part of the preparation of the article was assisted by an LLM. The authors retain full responsibility for the correctness of the final work. This article does not contain any AI written text.

\section{Preliminaries}\label{sec-prelim}

\subsection{Asymptotically conical expanding gradient Ricci solitons}\label{sec-prelim-1}

In this section, we explain the results obtained in \cite{CDS} to identify the soliton vector field associated to an asymptotically conical expanding gradient Ricci soliton with the radial vector field $\frac{1}{2}r\partial_r$ on the tangent cone at infinity.
 
 A Riemannian manifold $(M^n,g)$ has quadratic curvature decay with derivatives if for some (hence any) point $p\in M$ and all $k\in\mathbb{N}_{0}$,
\begin{equation*}
A_{k}(g):=\sup_{x\in M}|(\nabla^{g})^{k}\operatorname{Rm}(g)|_{g}(x)d_{g}(p,\,x)^{2+k}<\infty,
\end{equation*}
where $\operatorname{Rm}(g)$ denotes the curvature of $g$ and $d_{g}(p,\,x)$ denotes the distance between $p$ and $x$ with respect to $g$.

Finally, to each expanding gradient Ricci soliton $(M^n,g,\nabla^gf)$, consider the one-parameter family $(\varphi_t)_{t>0}$ of diffeomorphisms generated by $-\frac{\nabla^gf}{t}$ and define as in Siepmann \cite{Siepmann} the following sets, for any $a\geq 0$:
\begin{equation}\label{def-M-a}
M_{a}:=\Big\{x\in M\,|\,\liminf_{t\to 0^{+}}tf(\varphi_t(x))>a\Big\}.
\end{equation}
These sets are invariant under the diffeomorphisms $(\varphi_t)_{t\in(0,1]}.$ 
%\todo[inline]{FS: I see that from the definition these sets are invariant under the diffeomorphisms $\varphi_t$. But taking for example a positively curved expander, then aside from the point where $\nabla f =0$ do not all points move strictly 'inwards' under $\varphi_t$? What am I missing here?}

 \begin{theo}[Map to the tangent cone for expanding Ricci solitons]\label{ac-str-exp}
Let $(M^n,g,\nabla^gf)$ be a complete expanding gradient Ricci soliton with quadratic curvature decay with derivatives. 
Then on each end of $M$,
\begin{enumerate}
\item the one-parameter family of functions $(tf\circ\varphi_{t})_{t>0}$ converges to a non-negative
continuous real-valued function $q(x):=\lim_{t\to0^{+}}tf(\varphi_{t}(x))$ pointwise on $M$ as $t\to 0^{+}$. Moreover, the convergence is locally smooth on $M_0$.\\[-1ex]
\item The metrics $(g(t))_{t>0}$ converge to a Riemannian metric $\tilde{g}_{c}$ in $C^{\infty}_{\operatorname{loc}}(M_{0})$ as $t\to0^{+}$ and $\tilde{g}_{c}=2\operatorname{Hess}_{\tilde{g}_{c}}q$.\\[-1ex]
\item The function $q$ is proper and satisfies on $M_0$, $|\nabla^{\tilde{g}_{c}}q|^{2}_{\tilde{g}_{c}}=q$ and $\nabla^{\tilde{g}_{c}}q=\nabla^gf$. In particular, there exists $c_0>0$ such that the level sets $q^{-1}\{\frac{c^2}{4}\}$ are compact connected smooth hypersurfaces for $c\geq c_0.$\\[-1ex]
\item Define $\rho:=2\sqrt{q}$ on $M_0$ and let $S:=\rho^{-1}(\{c\})$ for any $c$ large enough. 
Then there exists a diffeomorphism $\iota:(c,\,\infty)\times S\to M_{\frac{c^{2}}{4}}$ such that $g_{c}:=\iota^{*}\tilde{g}_{c}=dr^{2}+r^{2}\frac{g_{S}}{c^{2}}$ and $d\iota(r\partial_{r})=2\nabla^gf$, where
$r$ is the coordinate on the $(c,\,\infty)$-factor and $g_{S}$ is the restriction of $\tilde{g}_{c}$ to $S$. Moreover, along $M_{\frac{c^{2}}{4}}$, we have that
\begin{equation*}\label{northbeach}
|(\nabla^{g_{c}})^{k}(\iota^{*}g-g_{c}-2\operatorname{Ric}(g_{c}))|_{g_{c}}=O(r^{-4-k})\quad\textrm{for all $k\in\mathbb{N}_{0}$}.
\end{equation*}
\end{enumerate}
\end{theo}
The statement and the proof of Theorem \ref{ac-str-exp} is the combination of the statements and the proofs contained in \cite[Section $3.1$]{CDS}.

An expanding gradient Ricci soliton $(M^n,g,\nabla^gf)$ with quadratic curvature decay with derivatives is then said to be asymptotically conical, a definition legitimated by Theorem \ref{ac-str-exp}.

\textbf{From now on, $(M^n,g,\nabla^gf)$ denotes an asymptotically conical expanding gradient Ricci soliton whose potential function is normalized such that }
\begin{equation}\label{norm-f-pot-fct}
|\nabla^gf|_g^2+\R_g=f-\frac{n}{2}.
\end{equation}

This normalisation fixes the free parameter in the soliton identities from Lemma \ref{id-sol-egs} once and for all. Moreover, one observes that thanks to this choice, the potential function $f$ is strictly positive on $M$ and that $\Delta_ff=f$ on $M$.

\subsection{Spectral theory of the linearized operator}\label{sec-spec-l-f}

For a given $(M^n,g,\nabla^gf)$ expanding gradient Ricci soliton, define the weighted measure $d\mu_f:=e^{f}d\mu_g$ and denote the corresponding $L^2$-space induced by $d\mu_f$ by $L^2_f$. According to \cite[Lemma $5.1$]{Mun-Wan-Wei-Ric-Sol}, the following Hardy-type inequality holds:
\begin{equation*}
\int_M\left(\R_g+\frac{n}{2}\right)\varphi^2\,d\mu_f\leq \int_M|\nabla^g\varphi|^2_g\,d\mu_f,\quad\forall \varphi\in C_c^{\infty}(M).
\end{equation*}
%\todo[inline]{F: I assume $C_c^{\infty}(M) = C_\text{c}^{\infty}(M)$? I would prefer to use $C_\text{c}^{\infty}(M)$.}
Now, thanks to Lemma \ref{id-sol-egs}, $\R_g+\frac{n}{2}$ is a positive function on $M$. Since the curvature of $g$ is assumed to decay at infinity, we infer that $\inf_M\R_g+\frac{n}{2}>0$ which proves that the bottom of the $L^2_f$-spectrum of $-\Delta_f$ is bounded from below by $\inf_M\R_g+\frac{n}{2}>0$.

Let $L_f:=\Delta_f+2\Rm(g)\ast$ be the linearized operator acting on symmetric $2$-tensors defined by:
\begin{equation*}
\begin{split}
&L_fh:=\underbrace{\Delta_gh+\nabla^g_{\nabla^gf}h}_{=:\Delta_fh}+2\Rm(g)\ast h,\quad h\in S^2T^*M,\\
&(\Rm(g)\ast h)_{ij}:=\Rm(g)_{ipqj}h^{pq}.\\
\end{split}
\end{equation*}
Then $L_f$ is symmetric with respect to the $L^2$-weighted scalar product induced by the weighted measure $e^fd\mu_g$ and denoted by $\langle\cdot,\cdot\rangle_{L^2_f}$:
\begin{equation*}
\langle L_fh_1,h_2\rangle_{L^2_f}:=\int_M\langle L_fh_1,h_2\rangle_g\,e^fd\mu_g=\langle h_1,L_fh_2\rangle_{L^2_f},\quad \forall h_i\in C^{\infty}_0(S^2T^*M),\,i=1,2.
\end{equation*}
Recall that the $L^2_f$-spectrum of $L_f$, denoted by $\sigma(L_f)$ is discrete and is bounded from above: see \cite{Der-Sta-Egs}. Let $E_{-}$, $E_0$ and $E_{+}$ be the subspaces of $L^2_f$ spanned by the eigentensors of $L^2_f$ associated to negative, zero and positive eigenvalues of $L_f$. Define $\pi_{-}$, $\pi_0$ and $\pi_+$ to be the $L^2_f$-projections onto the corresponding subspaces. Then $E_+$ and $E_0$ are finite dimensional. Let $\lambda_+:=\min \sigma(L_f)\cap \RR_+$ and $\lambda_{-}:=\max \sigma(L_f)\cap \RR_{-}$. Denote $\sigma(L_f)\cap \RR_{-}=:(\lambda_i^-)_i$ and let $(h_i^-)_i$ be a Hilbert basis of $E_{-}$ such that $L_fh_i^-=\lambda_i^-h_i^-$ for all $i$. Then it can be shown that eigentensors of $L_f$ associated to an eigenvalue $\lambda$ are smooth and decay quantitatively with respect to $\lambda$ (see \cite[Corollary $2.10$]{Der-Asy-Com-Egs}):
\begin{lemma}\label{lem-dec-eigen-lic}
Let $(M^n,g,\nabla^gf)$, $n\geq 3$, be an asymptotically conical expanding gradient Ricci soliton. If $L_fh=\lambda h$, $\lambda\in \RR$, with $h\in L^2_f$ then for each $k\geq 0$, there exists $C_k>0$ such that on $M$,
\begin{equation}\label{dec-eigen-lic}
 |\nabla^{g,\,k}h|_g\leq C_k\|h\|_{L^2_f}\,f^{-\lambda+\frac{k-n}{2}}e^{-f}.
\end{equation}
%\todo[inline]{TO: there was a problem, I believe linked to the convention with a minus, $L_fh=-\lambda h$, which was inconsistent with $\Delta_{g}H=-2\Rm(g)\ast H+\lambda H+\left(\frac{1}{4}|X|^2_{g}+\frac{1}{2}\Delta_{g}f\right)H.$ for instance, I think that I corrected it, but to be read again.\\F: Do we want to replace 'bounded curvature operator' by 'bounded curvature'?}
\end{lemma}
\begin{proof}
By conjugating $L_fh=\lambda h$ with a suitable weight, notice that the tensor $H:=e^{\frac{f}{2}}h$ satisfies
\begin{equation}
\begin{split}\label{love-drift-lap-conj-prelim}
\Delta_{g}H&=-2\Rm(g)\ast H+\lambda H+\left(\frac{1}{4}|\nabla^gf|^2_{g}+\frac{1}{2}\Delta_{g}f\right)H.
\end{split}
\end{equation}
As $|\Rm(g)|_g$ is bounded on $M$ and since the last two terms on the righthand side of \eqref{love-drift-lap-conj-prelim} are nonnegative thanks to Lemma \ref{id-sol-egs} 
%\todo{F: Don't we get just bounded from below?}
, this implies that $|H|_g$ satisfies the following differential inequality in the weak sense: since $\Delta_g f = \R_g+\frac{n}{2}>0,$
\begin{equation}
\Delta_{g}|H|\geq -C|H|_g.\label{love-drift-lap-conj}
\end{equation}

Since $(M^{n},\,g)$ is a Riemannian manifold with Ricci curvature bounded from below,
the results of \cite{Sal-Cos-Uni-Ell} yield the following local Sobolev inequality:
\begin{equation}
\begin{split}\label{sob-inequ-loc}
\left(\frac{1}{\operatorname{vol}_{g} (B_{g}(x,\,r))}\int_{B_{g}(x,\,r)}|\varphi|^{\frac{2n}{n-2}}\,d\mu_g\right)^{\frac{n-2}{n}}\leq \frac{C(r_0)r^2}{\operatorname{vol}_{g}(B_{g}(x,\,r))}\int_{B_{g}(x,\,r)}|\nabla^g\varphi|_{g}^2\,d\mu_g
\end{split}
\end{equation}
for any $\varphi\in H^1_0(B_{g}(x,\,r))$ and for all $x\in M$ and $0<r<r_0$, where $r_0>0$ is some fixed positive radius.

Performing a DeGiorgi-Nash-Moser iteration applied to \eqref{love-drift-lap-conj} satisfied by $|H|_{g}$ ensures that there exists a positive constant $C>0$ such that for a given ball $B_g(x,2)\subset M$, 
\begin{equation*}
\begin{split}
\sup_{B_g(x,1)}|H|_{g}^2&\leq C\int_{B_g(x,2)}|H|_g^2\,d\mu_g= C\int_{B_g(x,2)}|h|_g^2\,e^fd\mu_g\leq C\|h\|_{L^2_f}^2.
\end{split}
\end{equation*}
Then, one obtains the foolowing a priori uniform exponential decay, namely
\begin{equation*}
|h(x)|_g\leq Ce^{-\frac{f(x)}{2}}\|h\|_{L^2_f},\qquad \textrm{for all $x\in M$}.
\end{equation*}
To get the full decay on $|h|_g$ and its covariant derivatives, one invokes \cite[Corollary $2.10$]{Der-Asy-Com-Egs}.
\end{proof}
Recall that if $h_0\in L^2_f$ with $h_0=\sum_{\lambda\in \sigma(L_f)}a_{\lambda} h_{\lambda}$, the heat semi-group associated to $L_f$ acting on $L^2_f$ is defined by:
\begin{equation*}
P_{\tau}h_0:=\sum_{\lambda\in \sigma(L_f)}e^{\lambda \tau}a_{\lambda} h_{\lambda},\quad \tau>0.
\end{equation*}
Here the series converges in $L^2_f$ for any $\tau\geq 0$ and in $C^{\infty}_{loc}(\RR_+\times M)$: see \cite[Chapter $10$]{Gri-Boo}. Observe that for each eigenfunction $\lambda\in\sigma(L_f)$ and $h_{\lambda}\in L^2_f$ an eigentensor associated to $\lambda$, the heat semi-group acting on $h_{\lambda}$ extends to $\RR$ by setting $P_{\tau}h_{\lambda}:=e^{\lambda \tau}h_{\lambda}$ for $\tau\in\RR$.
%\todo[inline]{F: For $\tau = 0$ we have convergence in $L^2_f$, no? A: Yes, I added it}
%\todo[inline]{recall this fact [Deruelle-$L^2$-stability] and recall decay of eigentensors}
%\todo[inline]{recall linearization of Ricci tensor}
\subsection{Set-up of the equation}\label{sec-set-up-eqn}
We first recall in this section the linearisation of the Ricci tensor at a given metric $g$:
\begin{lemma}\label{app-C-monster}
For $g$ a Riemannian metric and $h$ a symmetric $2$-tensor such that $g+h$ is a Riemannian metric, 
\begin{equation*}
	\begin{split}
-2\Ric(g+h)&=-2\Ric(g)+\Delta_{L,\,g}h-\Li_{B_{g}(h)^{\sharp}}g+Q_{g}(h),\\
B_g(h)&:=\div_gh-\frac{1}{2}d\tr_gh,\quad \mbox{(Bianchi $1$-form)}\\
\Delta_{L,g}(h)&:= \Delta_g h + 2\Rm(g)\ast h - \Ric(g)\circ h - h\circ \Ric(g), \quad\mbox{(Lichnerowicz Laplacian)}\\
Q_{g}(h)&=
\left((g+h)^{-1}-g^{-1}\right)\ast \nabla^{g,\,2}h+(g+h)^{-1}\ast (g+h)^{-1}\ast \nabla^{g}h\ast\nabla^{g}h\\
&\quad+(g+h)^{-1}\ast g^{-1}\ast g^{-1}\ast g\ast \Rm(g)\ast h\ast h.
	\end{split}
	\end{equation*}
	where we slightly abuse notation with $\circ$ denoting composition after applying the suitable musical isomorphism.

	%\todo[inline]{T: there's an additional term $\Li_{B_{g}(h)}h$ in the linearization which of course disappears when linearizing the Ricci-Deturck flow. Should we directly linearize that R-D tensor or add the Lie derivative? A: I introduced two notions of Bianchi $1$-form, one stands for the linearized version $B_g(h)$ and the other one stands for the non-linear version $\tilde{B}_g(\cdot)$ that is only used in the lemmata below. F: I think we have not defined $\Delta_{L,g}$ at this point, no? A: done}
Moreover, one can write $Q_g(h)=:R_0(g,h)+(\nabla^{g})^{\ast}R_1(g,h)$ and there exists $C>0$ such that if $|h|_g\leq 1$, then
	\begin{equation*}
	\begin{split}
	|R_0(g,h)|_{g}&\leq C\left(|\nabla^{g}h|^2_{g}+|\Rm(g)|_{g}|h|^2_{g}\right),\\
	|R_1(g,h)|_{g}&\leq C|h|_g|\nabla^{g}h|_{g}.
	\end{split}
	\end{equation*}
	\end{lemma}
	
We also recall below how to get a solution to the (normalized) Ricci-DeTurck flow:
	
\begin{lemma}\label{lemma-eqn-h-tau}
If $(h(\tau))_{\tau}$ is a smooth one-parameter family of $S^2T^*M$ such that for each $\tau$, $g+h(\tau)=:\tilde{g}(\tau)$ is a metric on $M$, which satisfies:
\begin{equation}\label{eqn-h-tau}
\begin{split}
\partial_{\tau}h&=L_fh+Q_{g}(h(\tau)),\\
Q_{g}(h(\tau))&=R_0(g,h(\tau))+(\nabla^{g})^{\ast}R_1(g,h(\tau)),
%+\Li_{\nabla^{g+h(\tau)}f-\nabla^gf}h\\
%R_0(g,h(\tau))&:=(g+h(\tau))^{-1}\ast (g+h(\tau))^{-1}\ast \nabla^{g}h(\tau)\ast\nabla^{g}h(\tau)+(g+h(\tau))^{-1}\ast \Rm(g)\ast h(\tau)\ast h(\tau),\\
%R_1(g,h(\tau))_{kl}^i&:=((g+h(\tau))^{ij}-g^{ij})\nabla^g_jh(\tau)_{kl},
\end{split}
\end{equation}
then $\tilde{g}(\tau):=g+h(\tau)$ satisfies the following DeTurck normalized Ricci flow:
\begin{equation}
\begin{split}\label{oh-mymy}
\partial_{\tau}\tilde{g}(\tau)&=-2\Ric(\tilde{g}(\tau))-\tilde{g}(\tau)+\Li_{\nabla^{g}f}\tilde{g}(\tau)+\Li_{\tilde{B}_{g}(\tilde{g}(\tau))^{\sharp}}\tilde{g}(\tau),\\
\tilde{B}_{g}(\tilde{g}(\tau))_k&:=\tilde{g}(\tau)^{ij}\left(\nabla^g_i\tilde{g}(\tau)_{jk}-\frac{1}{2}\nabla^g_k\tilde{g}(\tau)_{ij}\right).
\end{split}
\end{equation}
\end{lemma}

\begin{rk}
The authors have decided to consider the above system \eqref{oh-mymy} instead of the following:
\begin{equation*}
\begin{split}
\partial_{\tau}\tilde{g}(\tau)&=-2\Ric(\tilde{g}(\tau))-\tilde{g}(\tau)+\Li_{\nabla^{\tilde{g}(\tau)}f}\tilde{g}(\tau)+\Li_{\tilde{B}_{g}(\tilde{g}(\tau))^{\sharp}}\tilde{g}(\tau),\\
\tilde{B}_{g}(\tilde{g}(\tau))_k&:=\tilde{g}(\tau)^{ij}\left(\nabla^g_i\tilde{g}(\tau)_{jk}-\frac{1}{2}\nabla^g_k\tilde{g}(\tau)_{ij}\right)+\tilde{g}(\tau)\left(\nabla^{g}f-\nabla^{\tilde{g}(\tau)}f,\cdot\right)_k,
\end{split}
\end{equation*}
the difference being that that we fix the vector field once and for all, a crucial property in the K\"ahler setting. This is convenient as it avoids time‑dependent weights in the nonlinear estimates.
\end{rk}

%\todo[inline]{$P_\tau$ to define, maybe $P_\tau h_i = e^{\tau\lambda_i}h_i$ if $L_fh_i = \lambda_i$. A: let's define it in the previous section}

%Notice that the first equation is an ill-posed parabolic Cauchy problem.

Recall that if $(h(\tau))_{\tau\,\leq\,0}$ satisfies \eqref{oh-mymy} then the curve of metrics defined by $\overline{g}(t):=t\varphi_t^*(g+h(\log t))$, where $\partial_t\varphi_t=-t^{-1}\nabla^gf\circ\varphi_t$, $\varphi_t\big|_{t=1}=\Id_M$, satisfies for $t\in(0,1]$:
\begin{equation}
\begin{split}\label{oh-mymy-bis}
\partial_{t}\overline{g}(t)&=-2\Ric(\overline{g}(t))+\Li_{\tilde{B}_{g(t)}(\overline{g}(t))^{\sharp}}\overline{g}(t),\\
\tilde{B}_{g(t)}(\overline{g}(t))_k&=\overline{g}(t)^{ij}\left(\nabla^{g(t)}_i\overline{g}(t)_{jk}-\frac{1}{2}\nabla^{g(t)}_k\overline{g}(t)_{ij}\right).
\end{split}
\end{equation}
From \eqref{oh-mymy-bis} and Lemma \ref{app-C-monster}, we obtain the following lemma:
\begin{lemma}\label{lemma-pde-h-unnorm}
If $(h(\tau))_{\tau\,\leq\,0}$ is a smooth one-parameter family of $S^2T^*M$ such that for each $\tau$, $g+h(\tau)=:\tilde{g}(\tau)$ is a metric on $M$, which satisfies \eqref{eqn-h-tau}, then the family $(\overline{h}(t):=t\varphi_t^*h(\log t))_{t\in(0,1]}$ satisfies:
%\todo[inline]{TO: We haven't found the fixed point yet, it should be about $(\overline{H}(t):=t\varphi_t^*H(\log t))_{t\in(0,1]}$--notation used afterwards.}
\begin{equation*}
\begin{split}
\partial_{t}\overline{h}(t)&=\Delta_{L,g(t)}\overline{h}(t)+Q_{g(t)}(\overline{h}(t)),\\
Q_{g(t)}(\overline{h}(t))&=R_0(g(t),\overline{h}(t))+(\nabla^{g(t)})^{\ast}R_1(g(t),\overline{h}(t)).
%R_0(g(t),\overline{h}(t))&:=\overline{g}(t)^{-1}\ast \overline{g}(t)^{-1}\ast \nabla^{g(t)}\overline{h}(t)\ast\nabla^{g(t)}\overline{h}(t)+\overline{g}(t)^{-1}\ast \Rm(g(t))\ast \overline{h}(t)\ast \overline{h}(t),\\
%R_1(g(t),\overline{h}(t))_{kl}^i&:=(\overline{g}(t))^{ij}-g(t)^{ij})\nabla^{g(t)}_j\overline{h}(t)_{kl}.
\end{split}
\end{equation*}
\end{lemma}

\subsection{Function spaces}\label{sec-defn-fct-spa}
%\todo[inline]{is there a global  embedding of the kind $H^m_f\hookrightarrow C^{k}$ provided $m>C(n,k)$ ? It should as $H^m\hookrightarrow C^k$ for $m>n/2+k$ and the norm on $H^m_f$ dominates that of $H^m$ as $f$ is bounded from below, a fact that is not true for shrinkers !...}

Define $$H_f^1:=\left\{h\in H^1_{loc}(S^2T^*M)\,\bigg|\,\|h\|_{L^2_f}+\|\nabla^{g}h\|_{L^2_f}<\infty\right\},$$ endowed with the norm: $$\|h\|_{H^1_f}:=\|\pi_-h\|_{L^2_f}+\|\nabla^{g}\pi_-h\|_{L^2_f}+\|\pi_0h\|_{L^2_f}+\|\pi_+h\|_{L^2_f}.$$
\begin{lemma}\label{lemma-eq-norms}
The following norms on $H^1_{f}$,
\begin{equation*}
\begin{split}
N_0(h)&:=\|\nabla^{g}\pi_-h\|_{L^2_f}+\|\pi_0h\|_{L^2_f}+\|\pi_+h\|_{L^2_f},\\
N_1(h)&:=\|\nabla^{g}\pi_-h\|_{L^2_f}+\|\nabla^{g}\pi_0h\|_{L^2_f}+\|\nabla^{g}\pi_+h\|_{L^2_f},\\
N_2(h)&:=\int_M|\nabla^g\pi_-h|^2_g-\langle2\Rm(g)\ast \pi_{-}h,\pi_-h\rangle_g\,d\mu_f+\|\pi_0h\|_{L^2_f}+\|\pi_+h\|_{L^2_f},
\end{split}
\end{equation*}
are all equivalent to $\|\cdot\|_{H^1_{f}}$.

In particular, if $h\in L^2_f$, $h\in H^1_f$ if and only if $\sum_i\left(|\lambda_i^{-}|+1\right)\langle h,h_i^-\rangle_{L^2_f}^2<\infty.$
\end{lemma}
\begin{proof}
Let $h\in H^1_f$. First observe that $N_0(h)\leq \|h\|_{H^1_f}$ by definition of $N_0$. Then as the bottom of the $L^2_f$-spectrum of $-\Delta_f$ acting on functions is positive as recalled at the beginning of Section \ref{sec-spec-l-f}
%\todo{F: Where do we know this from? Is there a reference? A: answered}
, Kato's inequality shows that $\|\nabla^g\pi_{-}h\|_{L^2_f}\geq c\| \pi_{-}h\|_{L^2_f}$ for some positive constant $c$ independent of $h$. This shows that $N_0$ and $\|\cdot\|_{H^1_f}$ are equivalent. 
Regarding the equivalence of $N_1$ and $\|\cdot\|_{H^1_f}$, as any tensor in $E_+$ or $E_0$ decays to $0$ at infinity thanks to \cite[Corollary 6.17]{Der-Sta-Egs}, it ensures that  $\|\nabla\cdot\|_{L^2_f}$ is a norm on $E_+$ or $E_0$, a fact that shows that the norms $\|\cdot\|_{L^2_f}$ and $\|\nabla\cdot\|_{L^2_f}$ are equivalent on $E_+$ and $E_0$ as $E_+$ and $E_0$ are finite dimensional.

Finally, by definition of $N_2$ and since $|\Rm(g)|_g$ is bounded on $M$, there is a positive constant independent of $h\in H^1_f$ such that $N_2(h)\leq C\|h\|_{H^1_f}.$
The Poincar\'e inequality $\int_M|\nabla^g\pi_-h|^2_g-\langle2\Rm(g)\ast \pi_{-}h,\pi_-h\rangle_g\,d\mu_f\geq -\lambda_{-}\|\pi_{-}h\|_{L^2_f}^2$ together with the trivial inequality $ \int_M|\nabla^g\pi_-h|^2_g-\langle2\Rm(g)\ast \pi_{-}h,\pi_-h\rangle_g\,d\mu_f\geq \|\nabla^g\pi_{-}h\|^2_{L^2_f}-C\|\pi_{-}h\|_{L^2_f}^2$ where $C:=C_0\sup_M|\Rm(g)|_g$, $C_0=C_0(n)$, gives for $\varepsilon\in(0,1)$:
\begin{equation*}
\int_M|\nabla^g\pi_-h|^2_g-\langle2\Rm(g)\ast \pi_{-}h,\pi_-h\rangle_g\,d\mu_f\geq ((1-\varepsilon)(-\lambda_{-})-\varepsilon C)\|\pi_{-}h\|_{L^2_f}^2+\varepsilon \|\nabla^g\pi_{-}h\|^2_{L^2_f}.
\end{equation*}
This leads to the equivalence of $N_2$ and $\|\cdot\|_{H^1_f}$ provided $\varepsilon$ is chosen small enough.
\end{proof}
%\todo[inline]{to be proved with the help of a suitable Hardy inequality}
We next define the function spaces within which we will work.
\begin{itemize}
\item Pick $\delta\in(0,\lambda_+)$ and define:
%$$|||h|||^2:=\int_{-\infty}^{0}\|h(\tau)\|^2_{H^1_f}\,d\tau.$$
\begin{equation*}
H^1_{f,\,\delta}:=\left\{h\in H^1_{loc}(M\times\RR_-)\,\bigg|\,\sup_{\tau\,\leq\,0}e^{-\delta\tau}\|h(\tau)\|_{L^2_f}<\infty,\,\text{and}\,\int_{-\infty}^0e^{-2\delta\tau}\|h(\tau)\|^2_{H^1_f}\,d\tau<\infty\right\},
\end{equation*}
endowed with the corresponding norm $\|h\|_{H^1_{\delta}}$ defined by:
\begin{equation*}
\|h\|_{H^1_{\delta}}:=\sup_{\tau\,\leq\,0}e^{-\delta\tau}\|h(\tau)\|_{L^2_f}+\left(\int_{-\infty}^0e^{-2\delta\tau}\|h(\tau)\|^2_{H^1_f}\,d\tau\right)^{\frac{1}{2}}.
\end{equation*}
We remove the dependency on $f$ in the notation of the norm to keep it light here.\\
\item For a given space-time point $(x,t)\in M\times(0,1]$, define its parabolic neighborhood by $P_r(x,t):=B_{g(t)}(x,r)\times (t-r^2,t]$ for $r^2\in(0,t)$. 
 Let $(\varphi_{t})_{t<0}$ denote the flow of $-\frac{X}{t}$ such that $\varphi_{t}\big|_{t=1}=\Id_M$. 
 
 Define for a symmetric $2$-tensor $H$ on $M\times(-\infty,0]$ the following symmetric $2$-tensor $\overline{H}$ on $M\times(0,1]$:
\begin{equation}\label{self-sim-H}
\overline{H}(x,t):=(t\varphi_{t}^*H)(x,\log t),\qquad x\in M,\qquad t\in(0,1].
\end{equation}
\item %Define finally the following Schauder-type function spaces for $k\geq 0$, $\alpha\in[0,\frac{1}{2})$ and $\iota_0\in(0,\inj_{g}(M)/4)$:

For $k\geq 0$ and $\alpha\in\left(0,\frac{1}{2}\right)$, denote the standard $C^{2k,\,2\alpha}$-norm of tensors on $P_{r}(x,t)$ for $2r<\inj_{g(t)}(x)$, by $\|\cdot\|_{C^{2k,\,2\alpha}(P_{r}(x))}$:
\begin{equation*}
\begin{split}
 \|T\|_{C^{2k,\,2\alpha}(P_{r}(x,t))}&:=\sum_{i+2j\,\leq\, 2k}r^{i+2j}\cdot\sup_{P_{r}(x,t)}|\nabla^{g(t),\,i}(\partial_{\tau}^{(j)})T|_{g(t)}\\
&\quad\ +\sum_{i+2j\,=\,2k}r^{2k+2\alpha}\cdot\left[\nabla^{g(t),\,i}(\partial_{\tau}^{(j)})T\right]
_{C^{0,\,2\alpha}(P_{r}(x,t))},
\end{split}
\end{equation*}
where $[\,\cdot\,]_{C^{0,2\alpha}}$ denotes the $2\alpha$-semi-norm on tensors on $M\times(0,1]$ induced by $g(t)$.\\
\item Since $\inj_{g(t)}x\geq \iota_0\sqrt{t}$ for all $x\in M$ and $t>0$ due to the self-similarity of the solution $g(t)=t\varphi_t^*g$ and the fact that $(M^n,g)$ has positive injectivity radius, let $\iota_0\in(0,1)$ be such that if $r:=\iota_0\sqrt{t}$, $2r<\inj_{g(t)}(x)$. 

Define then for $k\geq 0$, $\alpha\in\left(0,\frac{1}{2}\right)$ and $r:=\iota_0\sqrt{t}$,
\begin{equation*}
C^{2k,\,2\alpha}_{\delta}:=\left\{h\in C^{k,2\alpha}_{loc}(M\times(-\infty,0])\,\bigg|\,\sup_{t\in(0,1],\, x\in M}t^{-\delta}\|h\|_{C^{2k,\,2\alpha}(P_{\iota_0\sqrt{t}}(x,t))}<\infty\right\},
\end{equation*}
endowed with the corresponding norm $\|\cdot\|_{C^{k,2\alpha}_{\delta}}$ defined by:
\begin{equation*}
\|h\|_{C^{2k,\,2\alpha}_{\delta}}:=\sup_{t\in(0,1],\, x\in M}t^{-\delta}\|h\|_{C^{2k,\,2\alpha}(P_{\iota_0\sqrt{t}}(x,t))}.
\end{equation*}

\item Finally, \textbf{consider from now on} $h_0\in \operatorname{Im} \pi_+$  and define the final function space of interest here for $\alpha\in(0,\frac{1}{2})$:

\begin{equation*}
X^{2,2\alpha}_{\delta}(h_0):=\left\{h\in C^{2,2\alpha}_{\delta}\cap H^1_{f,\,\delta}\,\bigg|\,\pi_+(h(0))=h_0 \right\},
\end{equation*}
%\todo[inline]{TO: should it be $\pi_+(h(0))=h_0$? A: oups, yes indeed}
endowed with the corresponding norm $\|h\|_X:=\|h\|_{C^{2,2\alpha}_{\delta}}+\|h\|_{H^1_{\delta}}.$ 

The normed space $(X^{2,2\alpha}_{\delta}(h_0), \|\cdot\|_X)$ is a closed subspace of the Banach space $C^{2,2\alpha}_{\delta}\cap H^1_{f,\,\delta}$. In particular, $(X^{2,2\alpha}_{\delta}(h_0), \|\cdot\|_X)$ is a Banach space.
\end{itemize}

\begin{rk}\label{rk-simple-est}
Observe that by definition of the function space $C^{2}_{\delta}$, if $h\in C^{2}_{\delta}$ then there exists a constant $C>0$ such that for each $j=0,\,1,\,2$, for all $x\in M$ and $\tau\leq 0$, 
\begin{equation}\label{simple-est}
|\nabla^{g,\,j}h(x,\tau)|_g\leq C\|h\|_{C^{2}_{\delta}}e^{\delta \tau}.
\end{equation}
Indeed, by self-similarity of $g(t)$, $P_{\iota_0\sqrt{t}}(x,t)=B_{g}(\varphi_t(x),\iota_0)\times((1-\iota_0^2)t,t]$ and by the scaling properties of the tensor norms induced by $g(t)$, one has $t^{\frac{j}{2}}|\nabla^{g(t),\,j}\overline{h}|_{g(t)}=\varphi_t^*|\nabla^{g,\,j}h(\log t)|_{g}$ so that $$\sup_{B_g(\varphi_t(x),\iota_0)\times((1-\iota_0^2)t,t]}\varphi_t^*|\nabla^{g,\,j}h(\log t)|_{g}\leq C\|h\|_{C^{2}_{\delta}}\,t^{\delta}.$$ In particular, if $t=e^{\tau}$ for $\tau\leq 0$, and $y\in M$, $|\nabla^{g,\,j}h(y,\tau)|_{g}\leq C\|h\|_{C^{2}_{\delta}}\,e^{\delta \tau},$ since $\varphi_t$ is a diffeomorphism for each $t>0$. This leads to the above claimed estimate.
\end{rk}
%\todo[inline]{one may need/choose different exponents $\delta$ in the above definition}

%We hope for the following unstable result to be true:
\section{Existence of an unstable manifold}\label{sec-exi-uns-mfd}
Let $(M^n,g,\nabla^gf)$ be an asymptotically conical expanding gradient Ricci soliton such that $\sigma(L_f)\,\cap\,\RR_+\neq\emptyset$, i.e.~such that it is unstable. Let $h_0\in E_+$ where $E_+$ is defined in Section \ref{sec-spec-l-f}.

Consider the following formal map:
\begin{equation}\label{eq:def T}
\begin{split}
T:&(h(\tau))_{\tau\,\leq \,0} \longrightarrow \underbrace{P_{\tau}h_0-\int_{\tau}^0P_{\tau-s}\pi_+Q_{g}(h(s))\,ds}_{=:H_+(\tau)}\\
&+\underbrace{\int_{-\infty}^{\tau}\pi_0Q_{g}(h(s))\,ds}_{=:H_0(\tau)}+\underbrace{\int_{-\infty}^{\tau}P_{\tau-s}\pi_{-}Q_{g}(h(s))\,ds}_{=:H_{-}(\tau)},\quad \tau\leq0.
\end{split}
\end{equation}
%\todo[inline]{F: Why do we need $\pi_+h_0$ if $h_0\in E_+$? Do we want to specify already here in what space $h$ is?}
%\todo[inline]{F: on page 7 we define the heat semigroup $P_\tau$ for any $\tau > 0$, but in (3.1) we use it for $\tau < 0$, i.e. we are evolving backwards? I assume this works since we only apply it for $P_\tau h_0$ and  $P_{\tau-s} \pi_+Q_g(h(s))$, so in both cases it is a finite sum, so we can formally evolve backwards in time? A: I added a remark p.7}
The families of tensors $H_{\pm}$ and $H_0$ solve the following formal "Cauchy'' problems:
\begin{equation*}
\begin{split}
\partial_{\tau}H_+&=L_fH_{+}+\pi_{+}Q_{g}(h(\tau)),\quad H_+(0)=h_0,\\
\partial_{\tau}H_{-}&=L_fH_{-}+\pi_{-}Q_{g}(h(\tau)),\quad \lim_{\tau\rightarrow-\infty}e^{-(1-\varepsilon)\lambda_-\tau}\|H_-(\tau)\|_{L^2_f}=0,\quad \forall \varepsilon\in(0,1),\\
\partial_{\tau}H_0&=L_fH_0+\pi_0Q_{g}(h(\tau)),\quad \lim_{\tau\rightarrow -\infty}\|H_{0}(\tau)\|_{L^2_f}=0.\\
\end{split}
\end{equation*}
In particular, $H(\tau):=H_+(\tau)+H_{0}(\tau)+H_{-}(\tau)=T(h(\tau))$ solves on $M\times (-\infty,0)$,
\begin{equation*}
\partial_{\tau}H(\tau)=L_fH(\tau)+Q_g(h(\tau)).
\end{equation*}

The main result of this section is the following unstable manifold theorem:
%\todo[inline]{A: explicit dependence on $\|h_0\|_{L^2_f}$}
\begin{theo}\label{theo-ancient-sol}
For $\delta_0\in(\lambda_+/2,\lambda_+)$ and $\alpha\in (0,1/2)$, there exists $\varepsilon>0$ and $C_0>0$ such that for every $h_0\in \operatorname{Im} \pi_+$ satisfying $\|h_0\|_{L^2_f}<\varepsilon$, there exists a unique fixed point $h$ of $T$ in the ball $B_X(0,C_0\,\varepsilon)\subset X^{2,2\alpha}_{\delta_0}(h_0)$, i.e. there exists a unique solution $(h(\tau))_{\tau\,\leq\, 0}$ to $\partial_{\tau}h=L_fh+Q_{g}(h)$ with $\pi_+(h(0))=h_0$ such that $\|h\|_X<C_0\,\varepsilon$.

Moreover, there exists $C'_0=C'_0(n,\delta_0,\alpha)>0$ such that $\|h\|_{X}\leq C_0\|h_0\|_{L^2_{f}}$ and there exists a positive constant $C=C(n,\varepsilon,\delta_0,\alpha)$ such that the solution $(h(\tau))_{\tau\,\leq\, 0}$ satisfies:
\begin{equation}
\begin{split}\label{oh-dear-ene-est}
\sup_{\tau\,\leq\, 0}e^{-\lambda_+\tau}\|h(\tau)-P_{\tau}h_0\|_{L^2_f}&\leq C\|h\|_{C^2_{\delta}}\|h\|_{H^1_{\delta}}.
\end{split}
\end{equation}
\end{theo}
\subsection{Energy estimates}\label{sec-ene-est}
Before we start estimating each piece of the map $T(h)$ defined in \eqref{eq:def T}, we state the following elementary but helpful lemma that clarifies in which function space the non-homogeneous data $Q_g(h)$ lives:
\begin{lemma}\label{lemma-str-fct-space-Q}
If $h\in B_X(0,\gamma)\subset X^{2,2\alpha}_{\delta}(h_0)$ for some $\delta>0$ and $\gamma\in(0,1)$ then $Q_g(h)$ lies in $C^{0,2\alpha}_{2\delta-1}$ and for each $s\leq 0$, $Q_g(h(s))\in L^2_f$. Moreover, it satisfies:
\begin{equation*}
\int_{-\infty}^0\|Q_g(h(s))\|^2_{L^2_f}e^{-4\delta s}\,ds\leq C\|h\|^2_{C^2_{\delta}}\|h\|^2_{H^1_{\delta}},
\end{equation*}
for some uniform positive constant $C=C(n,\gamma,\delta)$.
\end{lemma} 

\begin{proof}
By the very definition of $Q_g(h)$ and that of the H\"older spaces, we immediately get that $Q_g(h)\in C^{0,2\alpha}_{2\delta-1}$ (it is essentially due to the ``quadratic'' nonlinearities of $Q_g(h)$ and the boundedness of $\Rm(g)$ and its covariant derivatives on $M$).

Now, thanks to [\eqref{simple-est}, Remark \ref{rk-simple-est}], observe that there exists a constant $C$ such that for $s\leq 0$,
\begin{equation*}
\begin{split}
|Q_g(h(s))|_g&\leq C\left(|h(s)|^2_g+|\nabla^gh(s)|^2_g+|\nabla^{g,\,2}h(s)|_g|h(s)|_g\right)\\
&\leq Ce^{\delta s}\|h\|_{C^2_{\delta}}\left(|h(s)|_g+|\nabla^gh(s)|_g\right),
\end{split}
\end{equation*}
pointwise on $M$. This implies in particular that for $s\leq 0$,
\begin{equation*}
\|Q_gh(s)\|^2_{L^2_f}\leq Ce^{2\delta s}\|h\|^2_{C^2_{\delta}}\|h(s)\|^2_{H^1_f}.
\end{equation*}
This shows that $Q_g(h(s))\in L^2_f$ for each $s\leq 0$. Finally, the expected last estimate follows once the previous inequality is multiplied across by $e^{-4\delta s}$ and integrated on $(-\infty,0]$.

\end{proof}

\begin{lemma}[$L^2$-estimate on $H_0$]\label{H_0-L2-norm}
Let $\delta>0$ and $\gamma\in(0,1)$. Then there exists $C>0$ such that if $h\in B_X(0,\gamma)\subset X^{2,2\alpha}_{\delta}(h_0)$, the family of symmetric $2$-tensors $(H_0(\tau))_{\tau\,\leq\, 0}$ lies in $H^1_{\delta}$ and,
\begin{equation*}
\sup_{\tau\,\leq\,0}e^{-2\delta\tau}\|H_0(\tau)\|_{L^2_f}+\left(\int_{-\infty}^0e^{-2\delta\tau}\|H_0(\tau)\|^2_{H^1_f}\,d\tau\right)^{\frac{1}{2}}\leq C\|h\|_{C^2_{\delta}}\|h\|_{H^1_{\delta}}.
\end{equation*}
\end{lemma}

\begin{proof}
 Let $(h_i^0)_i$ be an orthonormal basis of $E_0$. Observe first that $H_0(\tau)$ is a well-defined element of $L^2_f$ for each $\tau\leq 0$: indeed, by Lemma \ref{lemma-str-fct-space-Q}, $\pi_0(Q_g(h(s))$ is well defined for each $s\leq 0$ and for each $x\in M$,
\begin{equation*}
|\pi_0(Q_g(h(s))|_g(x)\leq \|Q_g(h(s))\|_{L^2_f}\left(\sum_{i}|h^0_i(x)|_g\right).
\end{equation*}
  Therefore, since $\|h_i^0\|_{L^2_f}=1$ for each index $i$, one has for $\tau\leq 0$,
\begin{equation*}
\int_M\left|\int_{-\infty}^{\tau}\pi_0(Q_g(h(s))(x)\,ds\right|^2_g\,d\mu_f(x)\leq C\left(\int_{-\infty}^{\tau}\|Q_g(h(s))\|_{L^2_f}\,ds\right)^2,
\end{equation*}
where $C$ depends on $\dim E_0$ and where the righthand side is finite thanks to Lemma \ref{lemma-str-fct-space-Q} again.

We have then obtained thanks to Cauchy-Schwarz inequality:
\begin{equation*}
\begin{split}
\|H_0(\tau)\|_{L^2_f}^2&\leq Ce^{4\delta\tau}\left(\int_{-\infty}^{\tau}e^{-4\delta s}\|Q_g(h(s))\|^2_{L^2_f}\,ds\right)\\
&\leq Ce^{4\delta \tau}\|h\|^2_{C^2_{\delta}}\|h\|^2_{H^1_{\delta}},
\end{split}
\end{equation*}
where we have used Lemma \ref{lemma-str-fct-space-Q} in the last line.

In particular, 
\begin{equation*}
\sup_{\tau\,\leq\,0}e^{-2\delta\tau}\|H_0(\tau)\|_{L^2_f}+\left(\int_{-\infty}^0e^{-2\delta\tau}\|H_0(\tau)\|^2_{H^1_f}\,d\tau\right)^{\frac{1}{2}}\leq C\,\|h\|_{C^2_{\delta}}\|h\|_{H^1_{\delta}},
\end{equation*}
where $C$ is a positive constant that may vary from line to line.

\end{proof}

\begin{lemma}[$L^2$-estimate on $H_+$]\label{H_+-L2-norm}
Let $\delta\in(\lambda_+/2,\lambda_+)$ and $\gamma\in(0,1)$. Then there exists $C>0$ such that if $h\in B_X(0,\gamma)\subset X^{2,2\alpha}_{\delta}(h_0)$, the family of symmetric $2$-tensors $(H_+(\tau))_{\tau\,\leq\, 0}$ lies in $H^1_{\delta}$ and if $\tau\leq 0$,
\begin{equation*}
\begin{split}
\sup_{\tau\,\leq\, 0}e^{-\lambda_+\tau}\|H_+(\tau)-P_{\tau}h_0\|_{L^2_f}&\leq C\|h\|_{C^2_{\delta}}\|h\|_{H^1_{\delta}}.
\end{split}
\end{equation*}

In particular, $\|H_+(\tau)-P_{\tau}h_0\|_{H^1_{\delta}}\leq C\|h\|_{C^2_{\delta}}\|h\|_{H^1_{\delta}},$ for each $\tau\leq 0$.
\end{lemma}
%\begin{rk}\label{rk-not-sharp}
%In the conclusion of Lemma \ref{H_+-L2-norm}, one could refine the estimates as follows: if $\delta\in (0,\lambda_+)$ then for every $\delta'\in [\delta,\lambda_+)$, there exists $C>0$ such that:
%\begin{equation*}
%\begin{split}
%\sup_{\tau\,\leq\, 0}e^{-\delta'\tau}\|H_+(\tau)-P_{\tau}\pi_+h_0\|_{L^2_f}&\leq C\sup_{\tau\,\leq\, 0}e^{-\delta \tau}\|h(\tau)\|_{C^2}\left(\int_{-\infty}^0e^{-2\delta\tau}\|h(\tau)\|^2_{H^1_f}\,d\tau\right)^{\frac{1}{2}}.
%\end{split}
%\end{equation*}
%\end{rk}

\begin{proof}
Similar to the proof of Lemma \ref{H_0-L2-norm}, if $Q(s):=Q_g(h(s))$ and $Q_i^+(s):=\langle Q(s),h_i^+\rangle_{L^2_f}$, then, since $\lambda_i^+\geq \lambda^+$ for each index $i$,
\begin{equation*}
\begin{split}
e^{-2\lambda_+\tau}\|H_+(\tau)-P_{\tau}h_0\|_{L^2_f}^2&=\sum_i\left(\int_{\tau}^0Q_i^+(s)e^{(\tau-s)\lambda_i^+}\,ds\right)^2e^{-2\lambda_+\tau}\\
&\leq \sum_i\left(\int_{\tau}^0Q_i^+(s)^2e^{-4\delta s}\,ds\right)\left(\int_{\tau}^0e^{2(\tau-s)\lambda_++4\delta s}\,ds\right)e^{-2\lambda_+\tau}\\
&\leq C\|h\|^2_{C^2_{\delta}}\|h\|^2_{H^1_{\delta}}\left(\int_{\tau}^0e^{2(2\delta-\lambda_+) s}\,ds\right)\\
&\leq C\|h\|^2_{C^2_{\delta}}\|h\|^2_{H^1_{\delta}}\left(\int_{-\infty}^0e^{2(2\delta-\lambda_+) s}\,ds\right)\\
&\leq \frac{C}{2(2\delta-\lambda_+)}\|h\|^2_{C^2_{\delta}}\|h\|^2_{H^1_{\delta}},
\end{split}
\end{equation*}
which proves the expected estimate.

\end{proof}

\begin{lemma}[$L^2$-estimate on $H_-$]\label{H_--L2-norm}
Let $\delta>0$ and $\gamma\in(0,1)$. Then there exists $C>0$ such that if $h\in B_X(0,\gamma)\subset X^{2,2\alpha}_{\delta}(h_0)$, the family of symmetric $2$-tensors $(H_{-}(\tau))_{\tau\,\leq\, 0}$ lies in $H^1_{\delta}$ and if $\tau\leq 0$,
\begin{equation*}
\begin{split}
\sup_{\tau\,\leq\,0}e^{-2\delta\tau}\|H_-(\tau)\|_{L^2_f}&\leq C\|h\|_{C^2_{\delta}}\|h\|_{H^1_{\delta}},\\
\|H_-\|_{H^1_{\delta}}&\leq C\|h\|_{C^2_{\delta}}\|h\|_{H^1_{\delta}}.
%\int_{-\infty}^{0}e^{-2\delta s}\|H_-(s)\|^2_{H^1_f}\,ds&\leq C\|h\|_{C^2_{\delta}}^2\|h\|_{H^1_{\delta}}^2.
\end{split}
\end{equation*}
\end{lemma}

\begin{proof}
 Let $(h_i^{-})_i$ be an orthonormal basis of $E_{-}$. For ease of reading, define $Q(s):=Q_g(h(s))$ and $Q_i^-(s):=\langle Q(s),h_i^-\rangle_{L^2_f}$ which is well-defined by Lemma \ref{lemma-str-fct-space-Q} for each $s\leq 0$. Then by Fubini's theorem, one has for each index $i$,
\begin{equation*}
\begin{split}
\int_M\left(\int_{-\infty}^{\tau}|P_{\tau-s}(\pi_{-}Q(s))|_g(x)\,ds\right)\,|h_i^{-}(x)|_g\,d\mu_f(x)&=\int_{-\infty}^{\tau}\int_M|P_{\tau-s}(\pi_{-}Q(s))|_g(x)|h_i^{-}(x)|_g\,d\mu_f(x)ds\\
&\leq\int_{-\infty}^{\tau}\|P_{\tau-s}(\pi_{-}Q(s))\|_{L^2_f}\,ds\\
&\leq \int_{-\infty}^{\tau}\|(\pi_{-}Q(s))\|_{L^2_f}\,ds\\
&\leq  \int_{-\infty}^{\tau}\|Q(s)\|_{L^2_f}\,ds<\infty,
\end{split}
\end{equation*}
thanks to Lemma \ref{lemma-str-fct-space-Q}. Here we have used Cauchy Schwarz inequality in the second line together with the fact that the heat semi-group is contractive on $E_{-}$ in the third line. This proves that for each index $i$, the scalar product $\langle H_-(\tau),h_i^-\rangle_{L^2_f}$ is well-defined and equals $\int_{-\infty}^{\tau}\langle P_{\tau-s}(\pi_{-}Q(s)),h_i^{-}\rangle_{L^2_f}\,ds$ which in turn equals $\int_{-\infty}^{\tau}e^{\lambda_i^-(\tau-s)}\langle Q(s),h_i^{-}\rangle_{L^2_f}\,ds$ by symmetry of the heat semi-group and the fact that $P_{\tau-s}h_i^{-}=e^{\lambda_i^-(\tau-s)}h_i^-$ for $\tau-s\in\RR$.

 A further application of Cauchy Schwarz inequality gives us:
\begin{equation*}
\begin{split}
\|H_-(\tau)\|^2_{L^2_f}&=\sum_i\langle H_-(\tau),h_i^-\rangle_{L^2_f}^2=\sum_i\left(\int_{-\infty}^{\tau}Q_i^-(s)e^{(\tau-s)\lambda_i^-}\,ds\right)^2\\
&\leq \sum_i\left(\int_{-\infty}^{\tau}Q_i^-(s)^2\,ds\right)\left(\int_{-\infty}^{\tau}e^{2(\tau-s)\lambda_i^-}\,ds\right)\leq \frac{1}{2|\lambda_-|}\int_{-\infty}^{\tau}\|Q(s)\|_{L^2_f}^2\,ds\\
&\leq C\|h\|_{C^2_{\delta}}^2\|h\|_{H^1_{\delta}}^2e^{4\delta \tau},
\end{split}
\end{equation*}
where we have used Lemma \ref{lemma-str-fct-space-Q} in the last line.

A similar reasoning shows on the basis of Lemma \ref{lemma-eq-norms}
 that:
\begin{equation*}
\begin{split}
\|H_-\|_{H^1_{\delta}}^2&\leq C\sum_{i}|\lambda_i^-|\int_{-\infty}^0\left(\int_{-\infty}^{\tau}|Q_i^-(s)|e^{(\tau-s)\lambda_i^-}\,ds\right)^2e^{-2\delta\tau}\,d\tau\\
&\leq \sum_i|\lambda_i^-|\int_{-\infty}^0\left(\int_{-\infty}^{\tau}|Q_i^-(s)|^2e^{-4\delta s}\,ds\right)\left(\int_{-\infty}^{\tau}e^{2(\tau-s)\lambda_i^-+4\delta s}\,ds\right)e^{-2\delta\tau}\,d\tau\\
&\leq \left(\int_{-\infty}^0\left(\sum_i\frac{|\lambda_i^-|}{2|\lambda_i^-|+4\delta}|Q_i^-(s)|^2e^{-4\delta s}\right)\,ds\right)\int_{-\infty}^0e^{2\delta\tau}\,d\tau\\
&\leq C\|h\|_{C^2_{\delta}}^2\|h\|_{H^1_{\delta}}^2,
\end{split}
\end{equation*}
as expected.
\end{proof}

By combining the previous Lemmata \ref{H_0-L2-norm}, \ref{H_+-L2-norm} and \ref{H_--L2-norm}, one arrives at:
\begin{coro}\label{coro-l2-est}
Let $\delta\in(\lambda_+/2,\lambda_+)$ and $\gamma\in(0,1)$. Then there exists $C>0$ such that if $h\in B_X(0,\gamma)\subset X^{2,2\alpha}_{\delta}(h_0)$,
\begin{equation*}
\begin{split}
\|T(h)-P_{\tau}h_0\|_{H^1_{\delta}}\leq C\|h\|_{C^2_{\delta}}\|h\|_{H^1_{\delta}}.
\end{split}
\end{equation*}
Moreover,
\begin{equation*}
\begin{split}
\sup_{\tau\,\leq\, 0}e^{-\lambda_+\tau}\|T(h)(\tau)-P_{\tau}h_0\|_{L^2_f}&\leq C\|h\|_{C^2_{\delta}}\|h\|_{H^1_{\delta}}.
\end{split}
\end{equation*}

\end{coro}
%\begin{rk}The first statement of Corollary \ref{coro-l2-est} would only need $\delta\in(0,\lambda_+)$.
% The restriction $(\lambda_+/2,\lambda_+)$ comes from Lemma \ref{H_+-L2-norm} where a sharp estimate is obtained on the time decay of the component $H_+$.
%\end{rk}

\subsection{Pointwise estimates}\label{sec-pt-est}

%For a given space-time point $(x,t)\in M\times(0,1]$, define its parabolic neighborhood by $P_r(x,t):=B_{g(t)}(x,r)\times (t-r^2,t]$ for $r^2\in(0,t)$. 
For a given $C^2_{loc}$ family of symmetric $2$-tensors $(\overline{h}(t))_{t\in(0,1]}$, if $(\overline{H}(t))_{t\in(0,1]}$ is a $C^2_{loc}$ solution to $\partial_{t}\overline{H}(t)=\Delta_{L,g(t)}\overline{H}(t)+Q_{g(t)}(\overline{h}(t))$ on $M\times(0,1]$, parabolic Schauder estimates read in this setting:
\begin{equation*}
\|\overline{H}\|_{C^{2,2\alpha}(P_r(x,t))}\leq C\left(r^2\|Q_{g(\cdot)}(\overline{h}(\cdot))\|_{C^{0,2\alpha}(P_{2r}(x,t))}+\|\overline{H}\|_{C^{0}(P_{2r}(x,t))}\right),\quad r:=\iota_0\sqrt{t}.
\end{equation*}

%From now on, define $r:=\iota_0\sqrt{t}$ for a possibly smaller uniform constant $\iota_0$. With such a choice of $r$, 
Consider the supremum over $x\in M$ in the previous estimate while multiplying across by $t^{-\delta}$  and apply it to $(\overline{H}(t):=t\varphi_t^*H(\log t))_{t\in(0,1]}$ where $H(\tau):=T(h)(\tau)$, $\tau\leq 0$, with $T$ defined in \eqref{eq:def T}, satisfying $\partial_{t}\overline{H}(t)=\Delta_{L,g(t)}\overline{H}(t)+Q_{g(t)}(\bar h(t))$ by Lemma \ref{lemma-pde-h-unnorm}, to get for $\delta>0$:
\begin{equation}\label{Schauder-par}
\begin{split}
\|H\|_{C^{2,2\alpha}_{\delta}}&\leq C\left(\|Q_{g}(h)\|_{C^{0,2\alpha}_{\delta-1}}+\|H\|_{C^{0}_{\delta}}\right)\\
&\leq C\left(\|h\|_{C^{2,2\alpha}_{\delta}}^2+\|H\|_{C^{0}_{\delta}}\right),
\end{split}
\end{equation}
 and where $C$ is a uniform positive constant.
%\todo[inline]{F: I am still confused here about what we do with time. A: I've changed the definition of the spaces $C^{2,2\alpha}_{\delta}$ p.9-10 and put as a remark the basic but crucial estimate \ref{simple-est} that we use in all the estimates below}
Now, thanks to Corollary \ref{coro-inter-inequ-Appleton} applied at $t=1$ to $H(\tau)$ for a given $\tau\leq 0$:
\begin{equation*}
\begin{split}
\|H(\tau)\|_{L^{\infty}(g)}&\leq C\|H(\tau)\|_{L^2}^{\frac{2}{n+2}}\|\nabla^gH(\tau)\|_{C^0}^{\frac{n}{n+2}}\\
&\leq\beta \|\nabla^gH(\tau)\|_{C^0}+C_{\beta}\|H(\tau)\|_{L^2(g)},
\end{split}
\end{equation*}
for any $\beta>0$.
In particular, as the $L^2_f$ norm dominates the $L^2$ norm induced by $g$, multiplying the previous inequality by $e^{-\delta\tau}$ and considering the supremum over $\tau\leq 0$ gives us:
\begin{equation*}
\begin{split}
\|H\|_{C^{0}_{\delta}}&\leq\beta \|H\|_{C^1_{\delta}}+C_{\beta}\sup_{\tau\,\leq \,0}e^{-\delta\tau}\|H(\tau)\|_{L^2_f},
\end{split}
\end{equation*}
thanks to [\eqref{simple-est}, Remark \ref{rk-simple-est}].

Feeding back \eqref{Schauder-par} with the previous estimate by choosing $\beta>0$ small enough so that after absorption by the lefthand side, one gets the following proposition:
\begin{prop}\label{prop-pt-est}
Let $\delta>0$, $\alpha\in(0,1/2)$ and $\gamma\in(0,1)$. Then there exists $C>0$ such that if $h\in B_X(0,\gamma)\subset X^{2,2\alpha}_{\delta}(h_0)$, 
\begin{equation*}\label{Schauder-par-bis}
\|T(h)\|_{C^{2,2\alpha}_{\delta}}\leq C\left(\|h\|^2_{C^{2,2\alpha}_{\delta}}+\sup_{\tau\,\leq \,0}e^{-\delta\tau}\|T(h)(\tau)\|_{L^2_f}\right).
\end{equation*}
\end{prop}
\subsection{Proof of Theorem \ref{theo-ancient-sol}}\label{sec-pf-thm}

We are in a position to prove Theorem \ref{theo-ancient-sol}:
\begin{claim}\label{T-preserve}
Given $\delta_0\in(\lambda_+/2,\lambda_+)$ and $\alpha\in(0,1/2)$, there exists $\varepsilon_0>0$ and a positive constant $C_0$ such that if $\varepsilon\in(0,\varepsilon_0]$ and $h_0\in \operatorname{Im} \pi_+$ satisfies $\|h_0\|_{L^2_f}<\varepsilon$ then the map $T$ preserves $B_X(0,C_0\,\varepsilon)$.
\end{claim}

%\todo[inline]{F: I guess this works for $0<\varepsilon<\varepsilon_0$? A: I changed it accordingly}

\begin{proof}[Proof of Claim \ref{T-preserve}]
According to Corollary \ref{coro-l2-est}, one gets from the definition of the norm on $X^{2,2\alpha}_{\delta_0}(h_0)$:
Given $\delta_0\in(\lambda_+/2,\lambda_+)$ and $\gamma=\frac{1}{2}$, there exists $C>0$ such that for $\tau\leq 0$,
\begin{equation*}
\begin{split}
\|T(h)-P_{\tau}h_0\|_{H^1_{\delta_0}}\leq C\|h\|^2_{X},
\end{split}
\end{equation*}
so that by the triangular inequality,
\begin{equation}
\begin{split}\label{hello-quad-est-l2}
\|T(h)\|_{H^1_{\delta_0}}&\leq \|P_{\tau}h_0\|_{H^1_{\delta_0}}+C\|h\|^2_{X}\leq C\|h_0\|_{L^2_f}+C\|h\|^2_{X}.
\end{split}
\end{equation}
Here we have used that $\delta<\lambda_+$ to estimate $\|P_{\tau}h_0\|_{H^1_{\delta_0}}$ from above.

Now, Proposition \ref{prop-pt-est} lets us estimate as follows:
\begin{equation}
\begin{split}\label{hello-quad-est-pt}
\|T(h)\|_{C^{2,2\alpha}_{\delta_0}}&\leq C\left(\|h\|^2_{C^{2,2\alpha}_{\delta}}+\|T(h)\|_{H^1_{\delta_0}}\right)\\
&\leq C\left(\|h\|^2_{C^{2,2\alpha}_{\delta}}+\|h_0\|_{L^2_f}+\|h\|^2_{X}\right)\\
&\leq C\|h\|^2_X+C\|h_0\|_{L^2_f},
\end{split}
\end{equation}
where we have used \eqref{hello-quad-est-l2} in the second line and where $C$ is a uniform positive constant that may vary from line to line.
Combining \eqref{hello-quad-est-l2} and \eqref{hello-quad-est-pt} gives us the existence of a positive constant $C$ such that for a given $C_0>0$ and $\varepsilon>0$ such that $C_0\varepsilon\leq 1/2=\gamma$, $$T(B_X(0,C_0\,\varepsilon))\subset B_X(0,C\varepsilon (1+C_0^2\varepsilon)),$$
provided $\|h_0\|_{L^2_f}\leq \varepsilon$. 

Therefore, if $C_0:=2\,\max\{C,1\}$ and $\varepsilon\leq \varepsilon_0:=C_0^{-2},$ then $T(B_X(0,C_0\,\varepsilon))\subset B_X(0,C_0\,\varepsilon)$ as required.
\end{proof}

\begin{claim}\label{T-cont}
Given $\delta_0\in(\lambda_+/2,\lambda_+)$ and $\alpha\in(0,1/2)$, there exists $\varepsilon>0$ and a positive constant $C_0$ such that if $h_0\in \operatorname{Im} \pi_+$ satisfies $\|h_0\|_{L^2_f}<\varepsilon$ then the map $T$ is a contraction on $B_X(0,C_0\,\varepsilon)$.
\end{claim}

\begin{proof}[Proof of Claim \ref{T-cont}]
Adapting the proof of Claim \ref{T-preserve}, there exists a positive constant $C$ such that for $h_i\in B_X(0,C_0\,\varepsilon)$, $i=1,2$:
\begin{equation*}
\begin{split}
\|T(h_2)-T(h_1)\|_{H^1_{\delta_0}}+\|T(h_2)-T(h_1)\|_{C^{2,2\alpha}_{\delta_0}}&\leq C(\|h_1\|_X+\|h_2\|_X)\|h_2-h_1\|_{X}.
\end{split}
\end{equation*}

Combining the two previous estimates lead to the desired claim by choosing $C_0$ large enough compared to the above constant $C$ and defining $\varepsilon_0:=C_0^{-2}$ as in the proof of Claim \ref{T-preserve}.
\end{proof}

According to Claim \ref{T-cont} and Picard's fixed point theorem applied to $T:B_X(0,C_0\,\varepsilon)\rightarrow B_X(0,C_0\,\varepsilon)$, there exists a unique fixed point of $T$ in $B_X(0,C_0\,\varepsilon)$ denoted by $h$. Moreover, combining estimates \eqref{hello-quad-est-l2} and \eqref{hello-quad-est-pt} to $h$ gives by absorption the desired quantitative estimate on $\|h\|_X$ in terms of $\|h_0\|_{L^2_{f}}<\varepsilon$, provided $\varepsilon$ is chosen small enough. Finally, Corollary \ref{coro-l2-est} applied to $T(h)=h$ gives the very last required estimate \eqref{oh-dear-ene-est}. This ends the proof of Theorem \ref{theo-ancient-sol}.

\section{Improved decay in space}\label{sec-dec-spa}
The purpose of this section is to show that not only the solution provided by Theorem \ref{theo-ancient-sol} lies in $L^2_f$ but also decays exponentially fast at spatial infinity uniformly in time.

Recall that $f$ is well-defined up to a constant. In this section we assume that the normalisation of $f$ is such that \eqref{equ:3} holds with $f+\tfrac{n}{2}$ on the righthand side. This especially implies that $f>0$ on $M$.

%\todo[inline]{F: Is this true ans sufficient? Do we want to put this earlier? A: I put in section 2 just before 2.2 this normalisation once and for all. }

%\subsection{Exponential decay in space}
%\todo[inline]{say a word about regularity of the solution and show that decay in time is $e^{\lambda_+\tau}$}
\begin{prop}\label{prop-alm-sharp-rate}
Under the setting and notations of Theorem \ref{theo-ancient-sol}, for $\delta\in(\lambda_+/2,\lambda_+)$ and $\alpha\in(0,1/2)$, there exists $C>0$ and $\varepsilon>0$ such that if $\|h_0\|_{L^2_f}<\varepsilon$ for $\tau\leq 0$, then
\begin{equation*}
|h(\tau)|_g\leq C\|h_0\|_{L^2_f}\,e^{\lambda_+\tau}f^{-\delta-\frac{n}{2}}e^{-f}.
\end{equation*}
Moreover, for each $k\geq 1$ and $\delta\in(\lambda_+/2,\lambda_+)$, there exists $C_k>0$ such that for $\tau\leq 0$, 
\begin{equation*}
|\nabla^{g,\,k}h(\tau)|_g\leq C_k\|h_0\|_{L^2_f}\,e^{\lambda_+\tau}f^{-\delta+\frac{k-n}{2}}e^{-f}.
\end{equation*}
\end{prop}

\begin{proof}
\begin{enumerate}
\item According to Lemma \ref{lemma-pde-h-unnorm},  the associated solution $(\overline{h}(t):=t\varphi_t^*h(\log t))_{t\in(0,1)}$ satisfies the following quasilinear parabolic equation:
\begin{equation}\label{eqn-h-over}
\partial_{t}\overline{h}(t)=\Delta_{L,\,g(t)}\overline{h}(t)+Q_{g(t)}(\overline{h}(t)),\quad t\in(0,1),
\end{equation}
a DeGiorgi-Nash-Moser iteration applied to $|\overline{h}(t)|_{g(t)}$ ensures that there exists a positive constant $C$ such that for every $t\in(0,1)$ and $x\in M$,
\begin{equation*}
%\label{DeGNM}
\sup_{(y,s)\in B_{g(t)}(x,\,\sqrt{t})\times\left(\frac{t}{2},t\right)}|\bar h(y,s)|_{g(s)}\leq C\left(\fint_{B_{g(t)}(x,\,2\sqrt{t})\times\left(\frac{t}{4},t\right)}|\overline{h}(y,s)|_{g(s)}^2\,d\mu_{g(s)}(y)ds\right)^{\frac{1}{2}}.
\end{equation*}
%\todo[inline]{F: Shouldn't this be $\bar h$ on the RHS? Otherwise the diffeomorphism moves the center in time, no?}
In particular, since $\vol_{g(t)}B_{g(t)}(x,\,2\sqrt{t})\geq v_0\,t^{\frac{n}{2}}$  for some uniform $v_0>0$, there exists a positive constant $C$ such that for every $\tau<0$ and $x\in M$,
\begin{equation*}
\begin{split}
|h(x,\tau)|_{g}&\leq C\left(e^{-\left(\frac{n}{2}+1\right)\tau}\int_{M\times\left(\tau-2\log 2,\tau\right)}|h(y,\tau')|_{g}^2\,d\mu_{g}(y)e^{\left(\frac{n}{2}+1\right)\tau'}d\tau'\right)^{\frac{1}{2}}\\
&\leq C\left(e^{-\left(\frac{n}{2}+1\right)\tau}\int_{M\times\left(\tau-2\log 2,\tau\right)}|h(y,\tau')|_{g}^2\,e^{f(y)}d\mu_{g}(y)e^{\left(\frac{n}{2}+1\right)\tau'}d\tau'\right)^{\frac{1}{2}}\\
&\leq C\left(e^{-\left(\frac{n}{2}+1\right)\tau}\int_{\tau-2\log 2}^\tau\|h(\tau')\|^2_{L^2_f}\,e^{\left(\frac{n}{2}+1\right)\tau'}d\tau'\right)^{\frac{1}{2}},
\end{split}
\end{equation*}
%\todo[inline]{F: Where does the $2\log 2$ come from?}
where $C$ is a uniform positive constant $C$ that may vary from line to line. Here we have used the fact that $f$ is bounded from below in the second line.

Thanks to Theorem \ref{theo-ancient-sol} and the previous estimate, there exists a positive constant $C$ such that if $\|h_0\|_{L^2_f}<\varepsilon$ and $\tau\leq 0$,
\begin{equation}
\begin{split}\label{DeGNM}
\sup_{ M}|h(\tau)|_{g}\leq C\|h_0\|_{L^2_f} e^{\lambda_+\tau}.
\end{split}
\end{equation}
%\todo[inline]{F: Where do we get this estimate from?}
\item As the associated solution $(\overline{h}(t))_{t\in(0,1)}$ satisfies \eqref{eqn-h-over}, parabolic regularity ensures that $(\overline{h}(t))_{t\in(0,1)}$ is smooth on $M\times(0,1)$. Therefore, $(h(\tau))_{\tau<0}$ is smooth on $M\times(-\infty,0)$ and satisfies the following interior Schauder estimates for parabolic equations thanks to \eqref{DeGNM}:  for each integer $k\geq 0$, there exists $C_k$ such that
\begin{equation}\label{pt-cov-der-time}
\sup_M|\nabla^{g,\,k}h(\tau)|_g\leq C_k\|h_0\|_{L^2_f} \,e^{\lambda_+ \tau},\quad \tau\leq 0.
\end{equation}
Notice that the previous bounds are a sharp improvement of Remark \ref{rk-simple-est}.

\item In view of \eqref{pt-cov-der-time}, an application of Corollary \ref{coro-inter-inequ-Appleton} to $h$ localized at infinity gives the existence of $\beta_0\in(0,1)$ such that $|h(\tau)|_g\leq C\|h_0\|_{L^2_f}\,e^{\lambda_+ \tau}e^{-\left(\frac{1-\beta_0}{2}\right)f}.$

We invoke the following standard interpolation inequalities whose proof can be found for instance in \cite{Aubin} for $m\geq 2$, an arbitrary tensor $T$ and $x\in M$:
\begin{equation}
\begin{split}\label{aub-int-inequ-fixed-sca}
&\sup_{B_{g}\left(x,\frac{\iota_0}{2}\right)}|\nabla^{g}T|_{g}\leq C\sup_{B_{g}(x,\iota_0  )}|T|_{g}^{1-\frac{1}{m}}\cdot \sup_{B_{g}\left(x,\iota_0 \right)}\bigg(\sum_{k=0}^m|\nabla^{g,\,k}T|_{g}\bigg)^{\frac{1}{m}},
\end{split}
\end{equation}
where $\inj_{g}(x)\geq 2\iota_0 $ for some positive constant $\iota_0$ uniform in $x\in M$, a property that holds thanks to the asymptotics of $g$ provided by Theorem \ref{ac-str-exp}.

Repeated applications of \eqref{aub-int-inequ-fixed-sca} to $h(\tau)$ and its covariant derivatives for a given $\tau\leq 0$ shows that for each nonnegative integer $k$, there exists $\beta_0\in (0,1)$ (possibly smaller than the previous one) and $C_{k,\,\beta_0}>0$ such that

\begin{equation}\label{intermed-improved-cov-der-exp-dec}
|\nabla^{g,\,k}h(\tau)|_g\leq C_{k,\,\beta_0}\|h_0\|_{L^2_f}\,e^{\lambda_+ \tau}e^{-\left(\frac{1-\beta_0}{2}\right)f}.
\end{equation}

\item From the previous item, one gets that there exists $\beta_0\in (0,1)$, $C_{\beta_0}>0$ and $C_0=C_0(n)>0$ such that 
\begin{equation}\label{rough-subsol-h}
\left(\partial_{\tau}-\left(\Delta_f+C_0|\Rm(g)|_g\right)\right)|h|_g\leq C_{\beta_0}\|h_0\|_{L^2_f}\,e^{2\lambda_+ \tau}e^{-(1-\beta_0)f}.
\end{equation}
In particular, in order to absorb the zeroth order curvature term, a straightforward computation based on soliton identities (cf. Lemma \ref{id-sol-egs}) show that for $C>0$:
\begin{equation*}
\begin{split}
&\left(\partial_{\tau}-\Delta_{f+2Cf^{-1}}\right)\left(e^{-Cf^{-1}}|h|_g\right)\leq C_{\beta_0} \|h_0\|_{L^2_f} \,e^{2\lambda_+ \tau}e^{-(1-\beta_0)f-Cf^{-1}}\\
&\quad+\underbrace{\left(C_0\sup_M(f|\Rm(g)|)-C+\left(2+Cf^{-1}\right)Cf^{-1}\right)}_{=:V_C}f^{-1}e^{-Cf^{-1}}|h|_g,
\end{split}
\end{equation*}
where $V_C\leq -\frac{C}{2}$ outside a compact set $K_C$ depending on $C$, provided $C$ is uniformly large.
Now, a similar computation shows that for $\beta\in(0,1)$ fixed and $C>0$:
\begin{equation*}
\begin{split}
\left(\partial_{\tau}-\Delta_{f+2Cf^{-1}}\right)&\Big(e^{-\beta f+\lambda_+\tau}\Big)\\
&\geq\left(\lambda_++\beta\left((1-\beta)f+\beta\frac{n}{2}\right)+\beta^2\R_g-2C\beta f^{-1}\right)e^{-\beta f+\lambda_+\tau}\\
&\geq\underbrace{\left(\lambda_++\beta(1-\beta)f-2C\beta f^{-1}\right)}_{=:W_C}e^{-\beta f+\lambda_+\tau},
\end{split}
\end{equation*}
where $W_C$ is a proper function bounded from below that grows linearly with $f$.
%\todo[inline]{A: we have $\Delta_gf=\R_g+n/2$ and $|\nabla^gf|^2_g+\R_g=f-n/2$ now so that $|\nabla^gf|^2_g\leq f$ since $\R_g+n/2>0$, I have modified the estimates above accordingly. F: I think we can drop the $n/2$, no?}
Therefore, the previous computation applied to $1-\beta_0$ (instead of $\beta$) gives for $A>0$ and $C>0$,
\begin{equation*}
\begin{split}
\left(\partial_{\tau}-\Delta_{f+2Cf^{-1}}\right)\Big(e^{-Cf^{-1}}|h|_g&-A\|h_0\|_{L^2_f} e^{-(1-\beta_0) f+\lambda_+\tau}\Big)\\
&\leq V_Cf^{-1}\left(e^{-Cf^{-1}}|h|_g-A\|h_0\|_{L^2_f} e^{-(1-\beta_0) f+\lambda_+\tau}\right)\\
&\quad+\left(C_{\beta_0}\|h_0\|_{L^2_f}+A\|h_0\|_{L^2_f}\left(f^{-1} V_C-W_C\right)\right)e^{-(1-\beta_0) f+\lambda_+\tau}.
\end{split}
\end{equation*}
Let us apply the maximum principle to the function $e^{-Cf^{-1}}|h|_g-A\|h_0\|_{L^2_f} e^{-(1-\beta_0) f+\lambda_+\tau}$ on $M\times[\tau_i,0]$ where $(\tau_i)_i$ is an arbitrary decreasing sequence of times going to $-\infty$: if it attains its maximum at an interior point $(x_0,\tau_0)\in M\times(\tau_i,0]$ then either it is nonpositive or it must be attained in a compact set of the form $\{f\leq \tilde{C}(A^{-1}+1)\}$. By choosing $A>0$ large enough such that the maximum of $e^{-Cf^{-1}}|h|_g-A\|h_0\|_{L^2_f} e^{-(1-\beta) f+\lambda_+\tau}$ is nonpositive on $\{f\leq \tilde{C}(A^{-1}+1)\}\times[\tau_i,0]$, a choice that can be made independently of the sequence $\tau_i$ thanks to \eqref{intermed-improved-cov-der-exp-dec}, one arrives at:
\begin{equation*}
\begin{split}
\sup_{M\times[\tau_i,0]}\Big(e^{-Cf^{-1}}|h|_g&-A\|h_0\|_{L^2_f}e^{-(1-\beta_0) f+\lambda_+\tau}\Big)\\
&\leq \max\left\{0,\sup_{M\times\{\tau_i\}}\left(e^{-Cf^{-1}}|h|_g-A\|h_0\|_{L^2_f}e^{-(1-\beta_0) f+\lambda_+\tau}\right)\right\}\\
&\leq \max\left\{0,\sup_{M\times\{\tau_i\}}\left(e^{-Cf^{-1}}|h|_g\right)\right\}.
\end{split}
\end{equation*}
The righthand side of the previous estimate converges to $0$ as $\tau_i$ tends to $-\infty$ thanks to \eqref{intermed-improved-cov-der-exp-dec}. Therefore, we have proved that for $\beta_0\in(0,1)$ as in \eqref{intermed-improved-cov-der-exp-dec}, 
\begin{equation}\label{almost-sharp-dec-h-bis}
|h(\tau)|_g\leq C_{\beta_0}\|h_0\|_{L^2_f}\,e^{\lambda_+ \tau}e^{-(1-\beta_0)f},\quad \tau\leq 0.
\end{equation}
Similarly to the proof of \eqref{intermed-improved-cov-der-exp-dec}, one also obtains that for each nonnegative integer $k$, there exists $C_{\beta_0,\,k}>0$ such that,
\begin{equation*}
|\nabla^{g,\,k}h(\tau)|_g\leq C_{\beta_0,\,k}\|h_0\|_{L^2_f}\,e^{\lambda_+ \tau}e^{-\frac{3}{4}(1-\beta_0)f},\quad \tau\leq 0.
\end{equation*}
This is a strict improvement of \eqref{intermed-improved-cov-der-exp-dec}. 

A reasoning as above shows then that for each nonnegative integer $k$ and $\beta\in(0,1)$, there exists $C_{\beta,\,k}>0$ such that,
\begin{equation*}
|\nabla^{g,\,k}h(\tau)|_g\leq C_{\beta,\,k}\|h_0\|_{L^2_f}\,e^{\lambda_+ \tau}e^{-(1-\beta)f},\quad \tau\leq 0.
\end{equation*}
\item Let us finally show that for every $\gamma\in(0,1)$, there exists $C_{\gamma}>0$ such that
\begin{equation}\label{dec-exp-almost-sharp}
|h(\tau)|_g\leq C_{\gamma}\|h_0\|_{L^2_f}\,e^{\lambda_+ \tau}f^{-(1-\gamma)\lambda_+-\frac{n}{2}}e^{-f},\quad\tau\leq 0.
\end{equation}
%\todo[inline]{TO: no $\delta$ in the equation above? A: good catch! I removed it from the statement}
Similarly to the proof of \eqref{almost-sharp-dec-h-bis}, if $\delta\in(0,\lambda_+)$,
\begin{equation}
\begin{split}\label{est-last-1}
&\left(\partial_{\tau}-\Delta_{f}\right)\left(f^{-\delta-\frac{n}{2}}e^{- f+\lambda_+\tau}\right)\geq\left(\lambda_+-\delta-C(n,g,\delta) f^{-1}\right)f^{-\delta-\frac{n}{2}}e^{- f+\lambda_+\tau}.
\end{split}
\end{equation}
Similar to the proof of \eqref{rough-subsol-h}, one also gets that for every $\beta\in (0,1)$, there exists $C_{\beta}>0$ and $C_0>0$ such that 
\begin{equation}\label{est-last-2}
\left(\partial_{\tau}-\left(\Delta_f+C_0|\Rm(g)|_g\right)\right)|h|_g\leq C_{\beta}\|h_0\|_{L^2_f}\,e^{2\lambda_+ \tau}e^{-(1+\beta)f}.
\end{equation}
On the basis of the two previous estimates \eqref{est-last-1} and \eqref{est-last-2}, and by considering a function of the form $e^{-Cf^{-1}}\left(|h|_g-A\|h_0\|_{L^2_f}f^{-\delta-\frac{n}{2}}e^{- f+\lambda_+\tau}\right)$ with $C$ and $A$ large enough, one arrives at the expected estimate \eqref{dec-exp-almost-sharp} by applying the maximum principle.
\item The proof of the decay on higher covariant derivatives $\nabla^{g,\,k}h$ is adapted from the one on $|h|_g$ by induction on $k\geq 0$: one establishes a parabolic equation satisfied by $\nabla^{g,\,k}h$ by differentiating [\eqref{eqn-h-tau}, Lemma \ref{lemma-eqn-h-tau}] successively and then one applies the maximum principle in the spirit of the proof that led to \eqref{dec-exp-almost-sharp}.
\end{enumerate}

\end{proof}

\section{Geometry of the solutions}\label{sec-geo-sol}

In this section, we consider an unstable asymptotically conical expanding gradient Ricci soliton $(M,g,\nabla^gf)$ which enjoys some curvature constraints to be made precise below. The purpose of this section is to show that the solutions provided by Theorem \ref{theo-ancient-sol} satisfy the same curvature constraints.

\subsection{Prerequisites on expanding gradient Ricci solitons}\label{sec-prelim-egs-curv-con}

We recall the following basic estimates from \cite{Der-Asy-Com-Egs}:
\begin{lemma}\label{lower-bd-ricci-curv-non-flat-cone}
Let $(M^n,g,\nabla^gf)$ be an asymptotically conical expanding gradient Ricci soliton with positive scalar curvature. Then there exists $c>0$ such that on $M$:
\begin{equation*}
 \R_g\geq \frac{c}{f},\quad \liminf_{f\rightarrow+\infty}f\Ric(g)(\nabla^gf,\nabla^gf)>0.
 \end{equation*}
\end{lemma}
%\todo[inline]{F: Does $c$ depend on the asymptotic cone? A: it does indeed.}
\begin{proof}
Regarding the first lower bound, we refer to \cite{Der-Asy-Com-Egs}. The lower bound on the radial directions can be easily derived from the proof of \cite{Der-Asy-Com-Egs} that we reproduce here as it is rather short: recall that $\Delta_f\R_g=-\R_g-2|\Ric(g)|^2_g$ on $M$. Now, the second contracted Bianchi identity from Lemma \ref{id-sol-egs} gives $2\Ric(g)(\nabla^gf,\,\cdot)=-d\R_g$. Once it is combined with the previous equation satisfied by $\R_g$, one gets:
\begin{equation*}
2\Ric(g)(\nabla^gf,\nabla^gf)=-g(\nabla^g\R_g,\nabla^gf)=\Delta_g\R_g+2|\Ric(g)|^2_g+\R_g.
\end{equation*}
%\todo[inline]{TO: there's an expression for mass that looks a lot like this integrated over large spheres in the AE setting... Not sure if there's any link.}
As the convergence to the cone is assumed to be smooth outside the tip, $$|\Delta_g\R_g+2|\Ric(g)|^2_g|\leq Cf^{-2}.$$ Thanks to the previously established lower bound on $\R_g$, the desired asymptotic lower bound on the Ricci curvature along the radial direction can be derived.
\end{proof}

\subsection{Preservation of curvature constraints}

%Assuming the background expander $(M,g)$ enjoys some curvature constraints, the purpose of this section is to show that the solution provided by Theorem \ref{theo-ancient-sol} also satisfies the same curvature constraints.
The main result of this section is the following proposition:

\begin{prop}\label{prop-curv-con}
Let $(M^n,g,\nabla^gf)$ be an unstable asymptotically conical expanding gradient Ricci soliton. With the notations of Theorem \ref{theo-ancient-sol}, there exists $\varepsilon>0$ such that if $\|h_0\|_{L^2_f}<\varepsilon$, the following assertions hold:\\[-2ex]
\begin{enumerate}
\item Assume the asymptotic cone $(C_0,d_c)$ has positive scalar curvature. Then the solutions provided by Theorem \ref{theo-ancient-sol} have positive scalar curvature.\\[-1ex]
\item Assume $(M^n,g)$ has positive Ricci curvature and the link $(X,g_X)$ of the asymptotic cone satisfies $\Ric(g_X)>(n-2)g_X$. Then the solutions provided by Theorem \ref{theo-ancient-sol} have positive Ricci curvature.\\[-1ex]
%\item Assume that both $(M^n,g)$ and the asymptotic cone $(C_0,g_c)$ have positive isotropic curvature (PIC). Then the solutions provided by Theorem \ref{theo-ancient-sol} are PIC too.\\[-1ex]
\end{enumerate} 
 \end{prop}
% \todo[inline]{A: the result above on PIC might be empty if $n=4$ and the assumption on the expander to be PIC provided the cone is PIC might be redundant?}
\begin{rk}
Nonnegative Ricci curvature is not necessarily preserved along the Ricci flow as demonstrated in \cite{Max-Ricci}.  Proposition \ref{prop-curv-con} ensures that positive Ricci curvature is preserved provided the solution at stake is a perturbation of a background asymptotically conical (unstable) expanding gradient Ricci soliton with positive Ricci curvature. Notice finally that the Ricci curvature of such a soliton decays at different speeds whether one considers directions parallel to the vector field $\nabla^gf$ or not.
\end{rk}

%\begin{rk}
%Proposition \ref{prop-curv-con} could be extended to the curvature conditions ``positive curvature operator'', ``2-positive curvature operator'', ``PIC2'' and ``PIC1''. However, it has been shown in \cite{Der-Lam} that an expanding gradient Ricci soliton with nonnegative sectional curvature is strictly stable, i.e. with the notations of Section \ref{sec-spec-l-f}, $E_+$ and $E_0$ are reduced to $\{0\}$. Since the condition ``nonnegative sectional curvature'' is implied by  the conditions ``nonnegative curvature operator'' and  ``weakly PIC2'', this would make the statement of Proposition \ref{prop-curv-con} empty. More recently, Lavoyer \cite{Lav-pers-com} has proved that a weakly PIC1 expanding gradient Ricci soliton is strictly stable too. Therefore, the remaining curvature conditions that can satisfy an unstable gradient Ricci soliton are ``positive scalar curvature'', ``positive Ricci curvature'' and ``PIC''. 
%\end{rk}
\begin{proof}
According to \cite{Bam-Che}, $(M^n,g)$ has positive scalar curvature. Lemma \ref{lower-bd-ricci-curv-non-flat-cone} ensures that $f\R_g\geq c>0$ on $M$. Thanks to Proposition \ref{prop-alm-sharp-rate} and the fact that if $|h(\tau)|_g\leq 1/2$, for some $\tau\leq 0$, 
\begin{equation*}
|\R_{g+h(\tau)}-\R_g|\leq C\left(|\nabla^{g,\,2}h(\tau)|_g+|\Rm(g)|_g|h(\tau)|_g+|\nabla^gh(\tau)|^2_g\right),
\end{equation*}
one arrives at $f\R_{g+h(\tau)}\geq c/2$ on $M$ for $\tau\leq 0$ provided $\varepsilon$ from Theorem \ref{theo-ancient-sol} is chosen small enough with respect to $c$.

If $(M^n,g)$ has positive Ricci curvature, then Lemma \ref{lower-bd-ricci-curv-non-flat-cone} gives us that $f\Ric(g)(\nabla^gf,\nabla^gf)\geq c>0$ outside a sufficiently large compact set containing the zero set of $\nabla^gf$. By a similar reasoning as above, one can show that $f\Ric(g+h(\tau))(\nabla^gf,\nabla^gf)\geq c/2>0$ outside a sufficiently large compact set. The assumption $\Ric(g_X)>(n-2)g_X$ guarantees in turn that the Ricci curvature of $g$ will be positive in the directions orthogonal to $\nabla^gf$ outside a sufficiently large compact set. More precisely, by Theorem \ref{ac-str-exp}, there exists $c>0$ such that 
\begin{equation*}
\liminf_{f\rightarrow+\infty}f\inf_{v\perp\nabla^gf,\,|v|_g=1}\Ric(g)(v,v)\geq c>0.
\end{equation*}
%\todo[inline]{F: Do we need to say something here about the decay of $g$ to the cone metric?}
In particular, the same estimate holds for the Ricci curvature of $g+h(\tau)$ uniformly in $\tau\leq 0$ for a slightly worse constant $c/2$. By choosing $\varepsilon$ from Theorem \ref{theo-ancient-sol} small enough, one can ensure that the Ricci curvature of $g+h(\tau)$ is positive for all $\tau\leq 0$ on $M$ as required.

\end{proof}

\subsection{Convergence of the solutions at the initial time}\label{sec-conv-init-time}

In this section, we aim at proving that the solutions obtained in Theorem \ref{theo-ancient-sol} give rise to solutions to the Ricci flow starting at the same asymptotic cone (up to isometry) as that of the background expanding gradient Ricci soliton.

%We are in a good position to show that any solution to the Ricci flow obtained as a perturbation of a given unstable expanding gradient Ricci soliton
\begin{prop}
Let $(M^n,g_b,\nabla^{g_b}f)$ be an asymptotically conical unstable expanding gradient Ricci soliton and let $(g_b(t))_{t>0}$ be its corresponding self-similar solution to the Ricci flow. Let $(\varphi_t)_{t>0}$ be the flow generated by $-t^{-1}\nabla^{g_b}f$ such that $\varphi_t\big|_{t=1}=\Id_M$. Let $(g(\tau))_{\tau\,\leq\, 0}$ be a solution provided by Theorem \ref{theo-ancient-sol}. Let $\overline{g}(t)=t\varphi_t^*g(\log t)$, $t\in(0,1]$, be the corresponding solution to \eqref{oh-mymy-bis} and let 
\begin{equation}\label{corr-rf-sol}
g_{RF}(t):=\psi_t^*\overline{g}(t)=t(\varphi_t\circ\psi_t)^*g(\log t),\quad t\in(0,1],
\end{equation}
 where $(\psi_t)_{t\in(0,1]}$ is the flow generated by $-\tilde{B}_{g_b(t)}(\overline{g}(t))^{\sharp_{\overline{g}(t)}}$ such that $\psi_{t}\big|_{t=1}=\Id_M$ (where $B$ is as in \eqref{oh-mymy-bis}). 
Then, \\[-2ex]
 \begin{enumerate}
 \item $(\overline{g}(t))_{t\in(0,1]}$ is a solution to the DeTurck-Ricci flow \eqref{oh-mymy-bis} which converges locally smoothly to $g_c$ on $C_0\setminus\{0\}$ as $t$ tends to $0^+$.\\[-1ex]
 \item There exists a smooth diffeomorphism $\tilde{\psi}_0$ of $C_0\setminus\{0\}$ such that $(g_{RF}(t))_{t\in(0,1]}$ is a solution to the Ricci flow which converges locally smoothly to $\tilde{\psi}_0^*g_c$ on $C_0\setminus\{0\}$ as $t$ tends to $0^+$.\\[-1ex]
  \item 
%\item \todo[inline]{$\overline{g}(t)$ is not self-similar provided $h_0:=g(0)-g$ is not identically $0$.}

 %If $\Ric(g_X)>(n-2)g_X$ and $\Ric(g(t))>0$ then $\Ric(\overline{g}(t))>0$ for all $t\in(0,1]$. In particular, the distances $d_{\overline{g}(t)}$ satisfy:
 %\begin{equation*}
 %d_c(x,y)-C\sqrt{t}\leq d_{\overline{g}(t)}(x,y)\leq d_c(x,y),\quad x,y\in M.
 %\end{equation*}
If $p$ lies in the zero set of $\nabla^{g_b}f$, $(M,d_{\overline{g}(t)},p)_{t\in(0,1]}$ and $(M,d_{g_{RF}(t)},\psi_t^{-1}(p))_{t\in(0,1]}$ converge to $(C_0,d_c,o)$ as $t$ tends to $0^+$ in the pointed Gromov-Hausdorff topology.

Moreover, if $\Ric(g)\geq 0$ then $(M,d_{\overline{g}(t)})$ converges pointwise to $(C_0,d_c)$ as $t$ tends to $0^+$.
 \end{enumerate}
% \todo[inline]{TO: I don't really understand the ``pointwise'' convergence of distance if we don't have the same  topology, is that GH?}
\end{prop}

\begin{proof}
Recall that $\overline{g}(t):=t\varphi_t^*g(\log t)=t\varphi_t^*g_b+\overline{h}(t)=g_b(t)+\overline{h}(t)$ where $(g_b(t))_{t\in(0,1]}$ denotes the self-similar solution to Ricci flow associated to the expanding gradient Ricci soliton $(M^n,g_b,\nabla^{g_b}f)$. 
%\todo[inline]{F: Are we using $g(t)$ both for the expanding solution and the our constructed RDTF?}
Now, Theorem \ref{ac-str-exp} ensures that $g_b(t)$ converges locally smoothly to $g_c$ on $C_0\setminus\{0\}$ (up to pulling back by the diffeomorphism $\iota$). As an intermediate conclusion, $\overline{g}(t)$ converges locally smoothly to $g_c$ on $C_0\setminus\{0\}$ as $t$ tends to $0^+$ since $\overline{h}(t)$ converges locally smoothly to $0$ on $C_0\setminus\{0\}$ (up to pulling back by the diffeomorphism $\iota$) by Proposition \ref{prop-alm-sharp-rate} as $t$ tends to $0^+$.\\

Regarding the second assertion, define $f_t:=t\varphi_t^*f$ so that $f_t$ converges locally smoothly to $r^2/4$ on $C_0\setminus\{o\}$ by Theorem \ref{ac-str-exp} (up to pulling back by the diffeomorphism $\iota$). Now, let us compute the time derivative of $\psi_t^*f_t$ as follows:
\begin{equation*}
\begin{split}
\partial_{t}\psi_{t}^*f_t&=-\psi_{t}^*\tilde{B}_{g_b(t)}(\overline{g}(t))^{\sharp_{\overline{g}(t)}}\cdot f_t+\psi_t^*\partial_tf_t\\
&=-\psi_{t}^*\tilde{B}_{g_b(t)}(\overline{g}(t))^{\sharp_{\overline{g}(t)}}\cdot f_t+(\varphi_t\circ\psi_t)^*\left(f-|\nabla^{g_b}f|^2_{g_b}\right)\\
&=-\psi_{t}^*\tilde{B}_{g_b(t)}(\overline{g}(t))^{\sharp_{\overline{g}(t)}}\cdot f_t+(\varphi_t\circ\psi_t)^*\left(\R_{g_b}+\frac{n}{2}\right)
,\quad t\in (0,1],
\end{split}
\end{equation*}
%\todo[inline]{F: How are we getting from the first to the second line? Where does the $n/2$ come from?}
%\todo[inline]{TO: Should it be $=-\psi_{t}^*\tilde{B}_{g(t)}(\overline{g}(t))+(\varphi_t\circ\psi_t)^*(\R_g+\frac{n}{2})$?}
where we have used Lemma \ref{id-sol-egs} in the last line. Therefore, Proposition \ref{prop-alm-sharp-rate} together with the quadratic curvature decay of $g$ gives,
\begin{equation*}
\begin{split}
|\partial_{t}\psi_{t}^*f_t|&\leq C\varepsilon t^{\lambda_+-\frac{1}{2}}\left(\frac{\psi_t^*f_t}{t}\right)^{-\delta+\frac{1-n}{2}}e^{-\frac{\psi_t^*f_t}{t}}\psi_t^*\underbrace{|\nabla^{\overline{g}(t)}f_t|_{\overline{g}(t)}}_{\leq \,C'\sqrt{f_t}}+C\\
&\leq CC'\varepsilon t^{\lambda_+}\left(\frac{\psi_t^*f_t}{t}\right)^{-\delta+\frac{2-n}{2}}e^{-\frac{\psi_t^*f_t}{t}}+C\leq C''
,\quad t\in (0,1],
\end{split}
\end{equation*}
where $C'$ is a positive constant independent of $t\in(0,1]$. 

%Therefore, $\psi_t^*f_t$ converges uniformly locally to a continuous function $f_0$ on $M$ and $\tilde{M}_{2a}:=\{x\in M\,|\, f_0(x)>2a\}\subset M_a\subset \tilde{M}_{\frac{a}{2}}$ where $M_a$ is defined in \eqref{def-M-a} for $a\geq 2C''$ where $C''$ is as in estimate \eqref{proper-fct-dist-f_t}.

%Observe that a similar estimate holds for $f_t$ and $\psi_t^*f$ in place of $\psi_t^*f_t$.

 In particular, 
\begin{equation}\label{proper-fct-dist-f_t}
f-C'(1-t)\leq \psi_{t}^*f_t\leq f+C'(1-t),\quad \text{on $M\times(0,1]$.}
\end{equation}
Observe that a similar estimate holds for $f_t$ in place of $\psi_t^*f_t$, more precisely, there exists $C>0$ such that for $0<s\leq t\leq 1$, $f_s-C(t-s)\leq f_t\leq f_s+C(t-s)$ on $M$. Therefore there exists $C>0$ such that for $t\in(0,1]$,
\begin{equation}\label{inequ-f_t-psi_t}
|\psi_t^*f_t-f_t|\leq C(1-t),\quad\text{on $M$}.
\end{equation}
As $\psi_t$ is a diffeomorphism of $M$ for each $t>0$, the previous estimate \eqref{inequ-f_t-psi_t} implies also:
\begin{equation}\label{inequ-f_t-psi_t-inverse}
|(\psi_t^{-1})^{*}f_t-f_t|\leq C(1-t),\quad\text{on $M$}.
\end{equation}

We claim that $(\psi_t)_t$ converges smoothly to a map $\psi_0$ uniformly locally on $M_c$ for some large but fixed $c>0$. 

Indeed, choose $c_0>0$  large enough so that $\psi_t(x)\in M_{c_0}$ for all $t\in(0,1]$ and $x\in M_{2c_0}$, a choice that can be made thanks to the previous estimates on $f_t$ and $\psi_t^*f_t$. Then Proposition \ref{prop-alm-sharp-rate} ensures that there exists $C>0$ such that for $t\in(0,1]$ and $x\in M_{2c_0}$,
\begin{equation*}
\begin{split}
|\partial_t\psi_t(x)|_{g_b(t)}&=|\tilde{B}_{g_b(t)}(\overline{g}(t))^{\sharp_{\overline{g}(t)}}|_{g_b(t)}(\psi_t(x))=t^{-\frac{1}{2}}\varphi_t^*|\tilde{B}_{g_b}(\overline{g})^{\sharp_{\overline{g}}}|_{g_b}(\psi_t(x))\leq C\varepsilon t^{\lambda_+-\frac{1}{2}},
\end{split}
\end{equation*}
where we have used the scaling properties of the Bianchi $1$-form and those of the norm induced by $g_b(t)$ in the second equality. Since $t\rightarrow t^{\lambda_+-\frac{1}{2}}$ is integrable on $(0,1]$, and since $g_b(t)$ is uniformly comparable to $\tilde{g}_c$ for $t\in(0,1]$ where $\tilde{g}_c$ is a cone metric isometric to $g_c$ as introduced in Theorem \ref{ac-str-exp}, one gets that $|\partial_t\psi_t(x)|_{\tilde{g}_c}\leq C\varepsilon t^{\lambda_+-\frac{1}{2}}$ for some positive uniform constant $C$ and 
\begin{equation*}
|\partial_td_{\tilde{g}_c}(\psi_{t}(x_0),\psi_s(x_0))\leq C\varepsilon\, t^{\lambda_+-\frac{1}{2}},\quad x_0\in M_{2c_0},
\end{equation*}
so that after integration between $t\in (0,1]$ and $s\in (0,1]$ (since this function vanishes at $t=s$!):
\begin{equation*}
d_{\tilde{g}_c}(\psi_{t}(x_0),\psi_s(x_0))\leq C\varepsilon \,t^{\lambda_++\frac{1}{2}},\quad t\geq s.
\end{equation*}
Therefore, $\psi_0(x_0):=\lim_{t\rightarrow 0^+}\psi_t(x_0)$ exists and lies in  $M_{c_0}$.

%Indeed, as $C_0\setminus\{0\}$ can be identified with $M_0$ (with the notations of Theorem \ref{ac-str-exp}), and since $|\psi_t^*\tilde{f}_t-\tilde{f}_t|\leq Ct$ for $t\in (0,1]$ due to \eqref{proper-fct-dist-f_t} and the remark right after it, if $x\in M_0$, then $\liminf_{t\rightarrow 0^+}  \psi_t^*\tilde{f}_t(x)>0$. The convergence of $\tilde{f}_t$ being locally uniform, one gets that $r(\psi_0(x_0))>0$ as expected.
%\todo[inline]{A: can someone check the scaling below please? F: Looks fine to me. In the second line should it be $\tilde{B}_{g_b}(g(\log t))$?}
Now, invoking the scaling properties of the Bianchi $1$-form and those of the norm induced by $g_b(t)$ again, one gets for $k\geq 1$,
\begin{equation}
\begin{split}\label{high-cov-der-g_b-psi}
 |\nabla^{g_b(t),\,k}\partial_t\psi_t|_{g_b(t)}(x)&=|\nabla^{g_b(t),\,k}\tilde{B}_{g_b(t)}(\overline{g}(t))^{\sharp_{\overline{g}(t)}}|_{g_b(t)}(\psi_t(x))\\
 &=t^{-\frac{k+1}{2}}\varphi_t^*|\nabla^{g_b,\,k}\tilde{B}_{g_b}(g(\log t))^{\sharp_{\overline{g}}}|_{g_b}(\psi_t(x))\\
 &\leq C_k\varepsilon t^{\lambda_+-\frac{k+1}{2}}\left(\frac{\psi_t^*f_t}{t}\right)^{-\delta+\frac{k+1-n}{2}}e^{-\frac{\psi_t^*f_t}{t}}\\
 &\leq C_ke^{-\frac{C(c_0)}{t}},
 \end{split}
 \end{equation}
 provided $x\in M_{2c_0}$ and where $C(c_0)$ is a positive constant. Here we have made use of Proposition \ref{prop-alm-sharp-rate} in the penultimate line together with \eqref{inequ-f_t-psi_t}
 in the last line.

As $g_b(t)$ converges smoothly locally to $\tilde{g}_c$ on $M_{2c_0}$ by Theorem \ref{ac-str-exp}, the estimates \eqref{high-cov-der-g_b-psi} can be converted for each $k\geq 0$ into a uniform bound on $|\nabla^{\tilde{g}_c,\,k}\partial_t\psi_t|_{\tilde{g}_c}(x)$ for $x\in M_{2c_0}$ and $t\in(0,1]$. This shows that the map $\psi_0$ is smooth.

In order to show that $\psi_0$ is a local diffeomorphism, it suffices to show that its differential is invertible. To do so, we observe that on $M_{2c_0}$, there exists $C>0$ such that for $t\in (0,1]$, 
\begin{equation*}
\begin{split}
 |\partial_t\psi_t^*g_b(t)|_{\psi_t^*g_b(t)}&=\left|-\psi_t^*\Li_{\tilde{B}_{g_b(t)}(\overline{g}(t))^{\sharp_{\overline{g}(t)}}}g_b(t)-2\psi_t^*\Ric(g_b(t))\right|_{\psi_t^*g_b(t)}\\
 &\leq C\varepsilon t^{\lambda_+-1}\left(\frac{\psi_t^*f_t}{t}\right)^{-\delta+\frac{2-n}{2}}e^{-\frac{\psi_t^*f_t}{t}}+\frac{C(g_b)}{t+f_t}\\
 &\leq C(c_0).
 \end{split}
 \end{equation*}
 Here we have invoked Proposition \ref{prop-alm-sharp-rate} with $k=2$ in the second line together with the quadratic Ricci curvature decay of $g_b$.
 
 By integrating in time the previous estimate in the form $\partial_t\psi_t^*g_b(t)\leq C(c_0)\psi_t^*g_b(t)$ on $M_{2c_0}\times(0,1]$, one arrives at the following coercive estimate for $0< t\leq 1$ for some positive uniform positive constant $C(c_0)$: $e^{-C(c_0)(1-t)}g_b\leq \psi_t^*g_b(t)$.

Sending $t$ to $0^+$ leads to $e^{-C(c_0)}g_b\leq \psi_0^*\tilde{g}_c$ on $M_{2c_0}$ which gives in particular that $\psi_0$ is an immersion.

Let us show that $\psi_0$ is injective as a map defined on $M_{2c_0}$. Assume there exist two points $x$ and $x'$ in $M_{2c_0}$ such that $\psi_0(x)=\psi_0(x')\in M_{c_0}$. The very same reasoning applied to the inverse of $\psi_t^{-1}$ restricted to $M_{c_0}$ with values into $M_{\frac{c_0}{2}}$ by \eqref{inequ-f_t-psi_t-inverse}  implies the existence of a smooth map $\vartheta_0$ defined on $M_{c_0}$ satisfying $\vartheta_0 \circ \psi_0=\Id.$ This forces $x=x'\in M_{2c_0}$ i.e. $\psi_0$ is a smooth diffeomorphism onto its image as required.

Let us now show that $g_{RF}(t)$ converges to $\psi_0^*\tilde{g}_c$ smoothly locally on $M_{2c_0}$ as $t$ tends to $0^+$.

By Theorem \ref{ac-str-exp} and the convergence properties of $\psi_t$ towards $\psi_0$ previously established, one can infer that for each compact subset $K$ of $M_{2c_0}$, there exists $C_K>0$ such that on $K\times(0,1]$, $C_K^{-1}\psi_0^{*}\tilde{g}_c\leq \psi_t^*g_b(t)\leq C_K\psi_0^{*}\tilde{g}_c$ and $ C_K^{-1}\psi_0^{*}\tilde{g}_c\leq \psi_t^*\tilde{g}_c\leq C_K\psi_0^{*}\tilde{g}_c$.
 The triangular inequality then guarantees that:
 \begin{equation}
 \begin{split}\label{lovely-tri-inequ}
 |g_{RF}(t)-\psi_0^*\tilde{g}_c|_{\psi_0^*\tilde{g}_c}&\leq |g_{RF}(t)-\psi_t^*g_b(t)|_{\psi_0^*\tilde{g}_c}+|\psi_t^*g_b(t)-\psi_t^*\tilde{g}_c|_{\psi_0^*\tilde{g}_c}+|\psi_t^*\tilde{g}_c-\psi_0^*\tilde{g}_c|_{\psi_0^*\tilde{g}_c} \\
 &\leq C|\psi_t^*\overline{h}(t)|_{\psi_t^*g_b(t)}+C|\psi_t^*g_b(t)-\psi_t^*\tilde{g}_c|_{\psi_t^*\tilde{g}_c}+|\psi_t^*\tilde{g}_c-\psi_0^*\tilde{g}_c|_{\psi_0^*\tilde{g}_c}\\
 &\leq C\varepsilon\, t^{\lambda_+}+C\psi_t^*|g_b(t)-\tilde{g}_c|_{\tilde{g}_c}+|\psi_t^*\tilde{g}_c-\psi_0^*\tilde{g}_c|_{\psi_0^*\tilde{g}_c},
 \end{split}
 \end{equation}
where we have used  Proposition \ref{prop-alm-sharp-rate} with $k=0$. Theorem \ref{ac-str-exp} ensures that the second term on the righthand side of the previous estimate \eqref{lovely-tri-inequ} converges to $0$ as $t\rightarrow 0^+$. The last term on the righthand side of  \eqref{lovely-tri-inequ} converges to $0$ as $t\rightarrow 0^+$ by smooth convergence of $\psi_t$ to the map $\psi_0$. Local smooth convergence of $g_{RF}(t)$ to $\psi_0^*\tilde{g}_c$ as $t\rightarrow 0^+$ can be proved along the same lines.
%\begin{equation*}
%|\nabla^{\psi_t^*g_b(t),\,k}\psi_t^*\overline{h}(t)|_{\psi_t^*g_b(t)}\leq C_k\varepsilon\, t^{\lambda_+}\tilde{f}^{-\delta+\frac{k-n}{2}}e^{-\tilde{f}},\quad t\in(0,1].
%\end{equation*}
%As an intermediate conclusion, we infer that $g_{RF}(t)-\psi_t^*g_b(t)$ converges locally uniformly smoothly to $0$ as $t$ tends to $0^+$.

As for the last assertion regarding the pointwise convergence of the distances induced by $\overline{g}(t)$, observe that Proposition \ref{prop-alm-sharp-rate} ensures that $|h(\tau)|_g\leq C_0\varepsilon \,e^{\lambda_+\tau}$ on $M\times(-\infty,0]$. Therefore, by definition of $g(t)$ and $\overline{g}(t)$, one gets for $t\in(0,1]$:
\begin{equation}\label{dist-est-g-bar-b}
\left(1-\sqrt{C_0\varepsilon}\, t^{\frac{\lambda_+}{2}}\right)d_{g_b(t)}\leq d_{\overline{g}(t)}\leq \left(1+\sqrt{C_0\varepsilon}\, t^{\frac{\lambda_+}{2}}\right)d_{g_b(t)}.
\end{equation}
Now, if $d_{g_b(t)}$ converges pointwise to $d_c$ as $t$ tends to $0$, the previous estimate ensures that  $d_{\overline{g}(t)}$ converges pointwise to $d_c$ as $t$ tends to $0$. 

If $\Ric(g)\geq 0$, \cite{Sim-3d} shows that $d_{g_b(t)}$ converges pointwise to $d_c$ as $t$ tends to $0$ which by the previous conclusion shows the desired result. 

Finally, if $p$ lies in the zero set of $\nabla^gf$, the pointed metric spaces $(M,d_{g_b(t)},p)$ for each $t\in(0,1]$ are isometric to $ (M,\sqrt{t} d_{g_b},p)$ which converges to $(C_0,d_c,o)$ in the pointed Gromov-Hausdorff topology as a consequence of Theorem \ref{ac-str-exp}. Now, \eqref{dist-est-g-bar-b} implies directly that $(M,d_{\overline{g}(t)},p)_{t\in(0,1]}$ converge to $(C_0,d_c,o)$ as $t$ tends to $0^+$ in the pointed Gromov-Hausdorff topology. Observe that $(M,d_{g_{RF}(t)},\psi_t^{-1}(p))$ is isometric to $(M,d_{\overline{g}(t)}, p)$ for each $t\in(0,1]$, a fact that implies that the limit of $(M,d_{g_{RF}(t)},\psi_t^{-1}(p))$ as $t\rightarrow 0^+$ is $(C_0,d_c,o)$ in the pointed Gromov-Hausdorff topology.
%and since $\psi_t(p)$ stays within bounded distance from $p$ for $t\in(0,1]$ by estimate \eqref{inequ-f_t-psi_t} (since $f_t(p)=tf(p)$ is bounded, \eqref{inequ-f_t-psi_t} forces $f_t(\psi_t(p))$ and $d_{\overline g(t)}(p,\psi_t(p))\leq C\sqrt{f_t(\psi_t(p))}$  to be uniformly bounded), the limit of $(M,d_{g_{RF}(t)},p)$ as $t\rightarrow 0^+$ is $(C_0,d_c,o)$ in the Gromov-Hausdorff topology.
\end{proof}

\section{Non self-similarity of the solutions}\label{sec-non-self-sim-sol}
In this section, we aim at proving items \eqref{item-grad} and \eqref{item-non-grad} from Theorem \ref{theo-main}.

\subsection{A PDE for the weighted divergence of a vector field}
Before we prove the non-self-similarity of the solutions provided by Theorem \ref{theo-ancient-sol}, we establish the linear equation satisfied by the weighted divergence of a vector field satisfying a Poisson equation involving the weighted Laplacian. Similar results were derived in the setting of shrinking gradient Ricci solitons in \cite{Col-Min-Rig-Cyl}. See also \cite{Bre-Rot} for similar equations in the setting of steady gradient Ricci solitons. Finally, the results below (specifically Lemmata \ref{lem-equiv-gauge-eqn} and \ref{lem-bochner-lie-der}) are very close to \cite[Section 6.2]{Bam-Che}, however, we present short proofs for consistency.

%\todo[inline]{F: Add definition of $\div_f$. A: see below}
Recall that the weighted divergence of a vector field $Y$ on a triple $(M^n,g,\nabla^gf)$ is defined by:
\begin{equation*}
\div_fY:=\div_gY+g(\nabla^gf,Y).
\end{equation*}
\begin{lemma}\label{lem-equiv-gauge-eqn}
Let $(M^n,g,\nabla^gf)$ be an expanding gradient Ricci soliton. Let $Y$ be a $C^{3}_{loc}$ vector field on $M$ such that 
\begin{equation}
\Delta_fY-\frac{Y}{2}+\nabla^g(\div_fY)=Z,\label{eq-quasi-kill}
\end{equation}
for some $C^1_{loc}$ vector field $Z$.

Then $h:=\Li_Y(g)$ and $\div_fY$ satisfy the following linear elliptic equations:
\begin{equation}
\begin{split}
L_fh&=\Li_{Z-\nabla^g\div_fY}g,\\
\Delta_f\div_fY-\frac{\div_fY}{2}&=\frac{\div_fZ}{2}.\label{decouple-div-eqn}
\end{split}
\end{equation}

\end{lemma}

\begin{proof}
Consider the flow $(\phi_t)_t$ generated by the vector field $Y$ and the family of metrics $g(t):=\phi_t^{\ast}g$. 
%Equation $(2.31)$ of \cite[Chapter $2$]{Cho-Lu-Ni-Boo} gives the variation of the Ricci curvature at the initial time.
%\todo[inline]{F: Don't we have this recorded in Lemma \ref{app-C-monster}? A: I fixed it}
Thanks to Lemma \ref{app-C-monster},
\begin{equation*}
\begin{split}
\frac{\partial}{\partial t}(-2\Ric(g(t)))\Big|_{t=0}&=\Delta_{L,g}(h)-\Li_{B_g(h)}g,\quad\mbox{where $h:=\frac{\partial}{\partial t}g(t)\Big|_{t=0}$}.
\end{split}
\end{equation*}
Now, in this case $h=\Li_Y(g)$, so the Bochner formula for vector fields implies
\begin{equation}\label{eq-div-li}
\div_g(h)^{\sharp}-\frac{1}{2}\nabla^g\tr_g(h)=\Delta_g Y+\Ric(g)(Y).
\end{equation}
 Using equation \eqref{eq-quasi-kill} satisfied by the vector field $Y$ and the soliton equation,
\begin{equation}
\div_g(h)^{\sharp}-\frac{1}{2}\nabla^g\tr_g(h)=\Delta_g Y+\nabla^g_Y\nabla^g f-\frac{Y}{2}=[Y,\nabla^g f]-\nabla^g(\div_fY)+Z.\label{formula-Bianchi-gauge-lie-der}
\end{equation}
We can conclude :
\begin{equation*}
-2\Li_Y(\Ric(g))=\Delta_{L,g}(\Li_Y(g))-\Li_{[Y,\nabla^g f]-\nabla^g(\div_fY)+Z}(g).
\end{equation*}
With $2\Ric(g)=\Li_{\nabla^g f}(g)-g$ from the soliton equation \eqref{egs-eqn}, one has
\begin{equation*}
-\Li_Y(\Li_{\nabla^g f}(g)-g)=\Delta_{L,g}(\Li_Y(g))-\Li_{[Y,\nabla^g f]-\nabla^g(\div_fY)+Z}(g).
\end{equation*}
Since $[\Li_Y,\Li_{\nabla^gf}](g)=\Li_{[Y,\nabla^gf]}(g)$,
%\todo[inline]{F: Does $[\Li_Y(g),\Li_{\nabla^gf}(g)]$ make sense? Should this not better be $[\Li_Y,\Li_{\nabla^gf}](g) = \Li_Y(\Li_{\nabla^gf}(g)) - \Li_{\nabla^gf}(\Li_Y(g)) = \Li_{[Y,\nabla^gf]}(g)$?}
\begin{equation}
\Li_Y(g)-\Li_{\nabla^g f}(\Li_Y(g))=\Delta_{L,g}(\Li_Y(g))-\Li_{Z-\nabla^g(\div_fY)}(g).\label{eqn-li-evo-eqn-form}
\end{equation}
Now, 
%Lemma \ref{lemm-link-lie-der-cov-der} applied to $T:=\Li_Y(g)$ and $X:=\nabla^gf$ gives:
\begin{equation}
\begin{split}
\Li_{\nabla^gf}(\Li_Y(g))&=\nabla^g_{\nabla^gf}\Li_Y(g)+\Li_Y(g)\circ \nabla^{g,2}f+\nabla^{g,2}f\circ \Li_Y(g)\\
&=\nabla^g_{\nabla^gf}\Li_Y(g)+\Li_Y(g)+\Li_Y(g)\circ \Ric(g)+\Ric(g)\circ \Li_Y(g).\label{rel-li-der-nabla-lie-der-V}
\end{split}
\end{equation}
%\todo[inline]{F: Have we define what $\circ$ is? I guess contraction?}

Here we have used \eqref{egs-eqn} once more in the second line.
In particular, with the help of \eqref{eqn-li-evo-eqn-form} and \eqref{rel-li-der-nabla-lie-der-V} by recalling that $h:=\Li_Y(g)$:
\begin{equation}
\begin{split}\label{evo-eqn-lie-der}
\Delta_fh+2\Rm(g)\ast h&=\Delta_{L,g}h+\Ric(g)\circ h+h\circ \Ric(g)+\nabla^g_{\nabla^gf}h\\
&=h-\Li_{\nabla^g f}h+\Ric(g)\circ h+h\circ \Ric(g)+\nabla^g_{\nabla^gf}h\\
&\quad+\Li_{Z-\nabla^g(\div_fY)}(g)\\
&=\Li_{Z-\nabla^g(\div_fY)}(g).
\end{split}
\end{equation}
%\todo[inline]{F: Where does the first identity come from?}
This proves the first required equation satisfied by $h$.

By tracing \eqref{evo-eqn-lie-der}, 
\begin{equation*}
\Delta_f(\div_g Y)+\left\langle\Ric(g),\Li_Y(g)\right\rangle_g=\div_gZ-\Delta_g(\div_fY),
\end{equation*}
which implies by using the soliton equation \eqref{egs-eqn},
\begin{equation}
\Delta_f(\div_g Y)-\div_gY=\div_gZ-\Delta_g(\div_fY)-\frac{1}{2}\left\langle\Li_{\nabla^gf}(g),\Li_Y(g)\right\rangle_g.\label{intermed-evo-eqn-div}
\end{equation}
Let us now compute the elliptic equation satisfied by the scalar product $g( Y,\nabla^gf)$: 
\begin{equation}
\begin{split}\label{evo-eqn-scal-prod}
\Delta_f(g(Y,\nabla^gf))=&g(\Delta_fY,\nabla^gf)+2\left\langle\nabla^gY,\nabla^{g,2}f\right\rangle_g+g(Y,\Delta_f\nabla^gf)\\
=&\frac{1}{2}g(Y,\nabla^gf)+g(Z-\nabla^g(\div_fY),\nabla^gf)\\
&+\frac{1}{2}\left\langle\Li_{\nabla^gf}(g),\Li_Y(g)\right\rangle_g+\frac{1}{2}g(Y,\nabla^gf).
\end{split}
\end{equation}
Here we have used \eqref{eq-quasi-kill} in the second line together with the fact that $\nabla^{g,2}f$ is symmetric, as well as equation \eqref{eigen-vec-nabla^gf}.
%\todo[inline]{F: Where does the last term come from?}

We conclude the proof of \eqref{decouple-div-eqn} by adding \eqref{evo-eqn-scal-prod} to \eqref{intermed-evo-eqn-div}.
\end{proof}

We continue this section by deriving a Bochner identity for vector fields in the setting of expanding gradient Ricci solitons:
\begin{lemma}\label{lem-bochner-lie-der}
Let $(M^n,g,\nabla^gf)$ be an expanding gradient Ricci soliton. Let $Y$ be a $C^{2}_{loc}$ vector field on $M$. Then,
\begin{equation*}
\left(\div_f\Li_Y(g)\right)^{\sharp}=\Delta_fY-\frac{Y}{2}+\nabla^g\div_fY.
\end{equation*}
\end{lemma}

\begin{proof}
Recall the following (unweighted) Bochner identity (see \eqref{eq-div-li}):
\begin{equation*}
\begin{split}
\left(\div_g\Li_Y(g)\right)^{\sharp}&=\frac{1}{2}\nabla^g\tr_g\Li_Y(g)+\Delta_gY+\Ric(g)(Y)\\
&=\nabla^g\div_gY+\Delta_gY-\frac{Y}{2}+\nabla^{g,\,2}f(Y),
\end{split}
\end{equation*}
where we have used \eqref{egs-eqn} in the last line.

Now,
\begin{equation*}
\begin{split}
\Li_Y(g)(\nabla^gf,\cdot)&=g\left(\nabla^g_{\nabla^gf}Y,\,\cdot\right)+g\left(\nabla^gf,\nabla^g_{\cdot}Y\right),\\
d\div_fY(\cdot)&=d\div_gY(\cdot)+g(\nabla^g_{\cdot}\nabla^gf,Y)+g(\nabla^gf,\nabla^g_{\cdot}Y).
\end{split}
\end{equation*}
The two previous set of identities combined with the (unweighted) Bochner identity leads to the desired result thanks to the symmetry of the hessian of $f$.
\end{proof}
The following weighted $L^2$-estimate is a Korn-type identity:
\begin{lemma}\label{korn-type-id}
Let $(M^n,g,\nabla^gf)$ be a complete expanding gradient Ricci soliton. Let $Y$ be a vector field such that $Y\in L^2_f$ and $\Li_Yg\in L^2_f$. Then $\nabla^gY \in L^2_f$, $\div_fY\in L^2_f$ and
\begin{equation}\label{korn-type-id-est}
\| \Li_Yg\|^2_{L^2_{f}}=2\|\nabla^{g}Y\|^2_{L^2_{f}}+\|Y\|^2_{L^2_{f}}+2\| \div_{f}Y\|^2_{L^2_{f}}.
\end{equation}
\end{lemma}
\begin{proof}
Applying Lemma \ref{lem-bochner-lie-der} to the vector field $Y$ gives:
\begin{equation*}
\Delta_f|Y|^2_g=2|\nabla^gY|^2_g+2\left\langle \left(\div_f\Li_Y(g)\right)^{\sharp}+\frac{Y}{2}-\nabla^g\div_fY,Y\right\rangle_g.
\end{equation*}
Let $\varphi(x):=\psi(d_g(p,\cdot)/R)$, $R\geq 1$, be a cut-off function such that $\varphi\equiv 1$ on $B_g(p,R)$, $\varphi\equiv 0$ on $M\setminus B_g(p,2R)$ with $|\nabla^g\varphi|_g\leq c/R$ almost everywhere, for some uniform positive constant $c$. Then,
\begin{equation*}
\begin{split}
\int_M\varphi^2\Delta_{f}|Y|^2_{g}\,d\mu_{f}
&=\int_M\varphi^2\bigg(2|\nabla^{g}Y|^2_{g}+|Y|_{g}^2+2\langle -\nabla^{g}\div_{f}Y,Y\rangle_{g}\\
&\qquad\qquad +2\langle \left(\div_f\Li_Y(g)\right)^{\sharp},Y\rangle_{g}\bigg)\,e^{f}d\mu_g.
\end{split}
\end{equation*}
An integration by parts shows:
\begin{equation}\label{IBP-top-notch}
\begin{split}
\int_M-2\langle Y,\nabla_{\nabla^g\varphi^2}^gY\rangle_{g}\,d\mu_{f}
&=\int_M\varphi^2\left(2|\nabla^{g}Y|^2_{g}+|Y|_{g}^2+2|\div_fY|^2-|\Li_Yg|^2_g\right)\,e^fd\mu_g\\
&\quad+\int_M\div_{f}Y\langle \nabla^g\varphi^2,Y\rangle_g-2 \Li_Yg(Y,\nabla^g\varphi^2)\,e^{f}d\mu_g.
\end{split}
\end{equation}
Therefore,  there exists $C>0$ independent of $R\geq 1$ such that:
\begin{equation*}
\begin{split}
\int_M&\left(|Y|^2_g+|\nabla^g\div_fY|^2+|\nabla^gY|^2_g\right)\varphi^2\,e^fd\mu_g\leq \\&C\int_M \left(|\nabla^g\varphi|_g|Y|_g\right)\varphi\left(|\div_fY|+|\Li_Yg|_g+|\nabla^gY|_g\right)+\varphi^2|\Li_Yg|_g^2\,e^fd\mu_g,
\end{split}
\end{equation*}
which after applying Young's inequality to the righthand side of the above estimate (and absorbing the terms involving $|\div_fY|$ and $|\nabla^gY|_g$ by the lefthand side) leads to,
\begin{equation}
\begin{split}\label{prelim-inte-lemma-prop}
\int_M\left(|Y|^2_g+|\nabla^g\div_fY|^2+|\nabla^gY|^2_g\right)\varphi^2\,e^fd\mu_g&\leq C\int_M \left(|\nabla^g\varphi|^2_g|Y|^2_g\right)\\
&\qquad\qquad +\varphi^2|\Li_Yg|_g^2\,e^fd\mu_g\\
&\leq C\left(\frac{c^2}{R^2}\|Y\|^2_{L^2_f}+\|\Li_Yg\|^2_{L^2_f}\right)\\
&\leq C\left(c^2\|Y\|^2_{L^2_f}+\|\Li_Yg\|^2_{L^2_f}\right).
\end{split}
\end{equation}
As the righthand side of the above estimate \eqref{prelim-inte-lemma-prop} is independent of $R\geq 1$, this shows that $\div_fY$ and $\nabla^gY$ lie in $L^2_f$.

Going back to \eqref{IBP-top-notch}, and sending $R$ to $+\infty$ leads to the expected identity. 
%\todo[inline]{the IBP's above can be done by using a cut-off fct, by invoking Claim \ref{claim-ene-div-X-h_0} and a pointwise estimate on $L_Xg_0$, what do you think? I can add more details if you want, T: could say $\nabla X\leq C' e^{-\gamma f_0}$ as well}
\end{proof}

\begin{lemma}\label{lemma-div-free}
Assume $L_{f}h=\lambda h$ with $\lambda\geq -1/2$ where $h$ is a symmetric $2$-tensor in $L^2_{f}$. Then the symmetric $2$-tensor $h$ is weighted divergence free, i.e. $\div_{f}h=0$ on $M$.

In particular, if $h\in E_+\oplus E_0$ then $\div_{f}h=0$.
\end{lemma}

\begin{proof}[Proof of Lemma \ref{lemma-div-free}]
Recall the following Bochner-type formula going back to Lichnerowicz adapted here to the setting of Ricci solitons whose proof can be found for instance in \cite[Lemma 6.13]{Bam-Che}:
\begin{equation*}
\div_{f}(L_{f}h)=\Delta_{f}(\div_{f}h)-\frac{1}{2}\div_{f}h.
\end{equation*}

In particular, as $L_{f}h=\lambda h$, then, $\Delta_{f}\alpha=\left(\frac{1}{2}+\lambda\right)\alpha$ where $\alpha:=\div_{f}h$ is a $1$-form. Therefore, the function $|\alpha|_{g}^2$ satisfies $\Delta_{f}|\alpha|_{g}^2\geq (1+2\lambda)|\alpha|_{g}^2$ on $M$.  Thanks to [\eqref{dec-eigen-lic}, Lemma \ref{lem-dec-eigen-lic}], $|\nabla^{g}h|_{g}=O(f^{\frac{1-n}{2}-\lambda}e^{-f})$ so that $|\alpha|_{g}=o(1)$ as $f\rightarrow +\infty$. The use of the maximum principle allows us to conclude that $\alpha\equiv 0$ since $\lambda\geq -1/2$.
%\todo[inline]{F: Why can we apply the maximum principle here? $\alpha_0$ is a $1$-form, not a scalar? A: we apply the max ppe to the norm $|\alpha|_{g_0}$}
Finally, if $h\in E_+\oplus E_0$, $h$ is a finite linear combination of eigentensors associated to nonnegative eigenvalues of $L_{f}$ which leads to $\div_{f}h=0$ thanks to the previous result.
\end{proof}

We end this section with a suitable maximum principle for subsolutions of the weighted laplacian that grows sub-linearly at infinity:
\begin{lemma}\label{lemma-max-ppe-sub-lin}
Let $(M^n,g,\nabla^gf)$ be a complete expanding gradient Ricci soliton such that $f$ is proper. Let $u:M\rightarrow\RR$ be a $C^2_{loc}$ function satisfying on $M$:
\begin{equation*}
\Delta_fu\geq \frac{1}{2}u-C_0e^{-\gamma f},
\end{equation*}
for some $\gamma\in[0,1]$ and $C_0\geq 0$.
Assume $\max_{ \{f\,=\, r_i^2\}}u=o(r_i)$ for some diverging sequence of positive reals $(r_i)_i$. Then $u\leq 3C_0\,e^{-\gamma f}$ on $M$.
\end{lemma}
%\todo[inline]{F:What is the constant $a>0$ for? It does not seem to make sense here.}
\begin{proof}
Thanks to Lemma \ref{id-sol-egs} and the normalization of $f$ chosen in \eqref{norm-f-pot-fct}, one gets that on $M$, 
\begin{equation*}
\Delta_fe^{-\gamma f}=\left(-\gamma\Delta_ff+\gamma^2|\nabla^gf|^2_g\right)e^{-\gamma f}=\left(-\gamma \Delta_gf-\gamma(1-\gamma)|\nabla^gf|^2_g\right)e^{-\gamma f}\leq 0,
\end{equation*} 
as $\Delta_gf>0$ by Lemma \ref{id-sol-egs}.

Therefore, if $\varepsilon>0$ and $A>2C_0$ (say $A:=3C_0$), 
\begin{equation}
\begin{split}\label{u-eps-f-sqrt}
\Delta_{f}\left(u-A e^{-\gamma f}-\varepsilon f^{\frac{1}{2}}\right) &\geq \frac{1}{2}\left(u-A e^{-\gamma f}-\varepsilon f^{\frac{1}{2}}\right)+\left(\frac{A}{2}-C_0\right)e^{-\gamma f}\\&>\frac{1}{2}\left(u-A e^{-\gamma f}-\varepsilon f^{\frac{1}{2}}\right).
\end{split}
\end{equation}

Observe that for each index $i$, $\{f\,\leq\, r_i^2\}$ is a compact subset of $M$ thanks to the properness of $f$. The maximum principle applied to \eqref{u-eps-f-sqrt} on $\{f\,\leq\, r_i^2\}$  gives that 
\begin{equation*}
\max_{ \{f\,\leq\, a_i^2\}}\left(u-A e^{-\gamma f}-\varepsilon f^{\frac{1}{2}}\right)=\max\left\{0\,,\,\max_{ \{f\,=\,  r_i^2\}}\left(u-A e^{-\gamma f}-\varepsilon f^{\frac{1}{2}}\right)\right\}.
\end{equation*}

This leads to the required upper bound on $u$ on $M$ by taking the limit $i \to \infty$ and then $\varepsilon \to 0$. This ends the proof of the lemma.
\end{proof}

\subsection{Proof of item \ref{item-non-grad}, Theorem \ref{theo-main}}

Assume that $(M,g_1,X_1)$ is an expanding Ricci soliton with $g_1:=g_0+h(0)$. By Corollary \ref{cor-lin-growth}, there exists $\chi_1>0$ such that
\begin{equation}\label{eqn-sol-g_0}
\Ric(g_1)+\frac{g_1}{2}=\frac{1}{2}\Li_{X_1}g_1,\quad |X_1|_{g_0}\leq \chi_1f_0^{\frac{1}{2}},
\end{equation}

Define from now on $h:=h(0)$ and let us linearize the soliton equation \eqref{eqn-sol-g_0} with respect to that satisfied by $(M,g_0,X_0:=\nabla^{g_0}f_0)$ by using the fact that $\Li_{X_0}h=\nabla^{g_0}_{X_0}h+\nabla^{g_0,\,2}f_0\circ h+h\circ\nabla^{g_0,\,2}f_0$:
\begin{equation*}
\begin{split}
0&=-2\Ric(g_1)-g_1+\Li_{X_1}g_1\\
&=-2\Ric(g_0)-g_0+\Li_{X_0}g_0+\Delta_{L,\,g_0}h-\Li_{B_{g_0}(h)^{\sharp}}g_0-h+\Li_{X_0}h+Q_{g_0}(h)+\Li_{X_1-X_0}g_1\\
&=\Delta_{g_0}h+2\Rm(g_0)\ast h-\Ric(g_0)\circ h-h\circ \Ric(g_0)-\Li_{B_{g_0}(h)^{\sharp}}g_0\\
&\quad-h+\nabla^{g_0}_{X_0}h+\nabla^{g_0,\,2}f_0\circ h+h\circ\nabla^{g_0,\,2}f_0+Q_{g_0}(h)+\Li_{X_1-X_0}g_1\\
&=L_{f_0}h-\Li_{B_{g_0}(h)^{\sharp}}g_0+Q_{g_0}(h)+\Li_{X_1-X_0}g_1,
\end{split}
\end{equation*}
where $B_{g_0}(h)$ and the non-linear term $Q_{g_0}(h)$ are introduced in Lemma \ref{app-C-monster}. Here we have used the soliton equation \eqref{eqn-sol-g_0} satisfied by $(g_0,X_0)$ twice: in the third and last equality.
%\todo[inline]{F: How do we get from $\Delta_{L,g_0}$ here to $L_{f_0}h_0$? What happens to the Ricci terms (compare the equations after (6.2))? How do we get the $\nabla_{\nabla f}h$ term? Where does $B_{g_0}(h_0)$ come from? Looks ok now, but there were two sign errors.}

Therefore,
\begin{equation}
\begin{split}\label{hello-eqn-non-linear}
L_{f_0}h&=\Li_{B_{g_0}(h)^{\sharp}}g_0-\Li_{X_1-X_0}(g_1)-Q_{g_0}(h)\\
&=\Li_{B_{g_0}(h)^{\sharp}-(X_1-X_0)}g_0-\Li_{X_1-X_0}h-Q_{g_0}(h).
\end{split}
\end{equation}
%where $X:=\nabla^{g_0}f_0-\nabla^{g}f.$
Before proving item \ref{item-non-grad}, Theorem \ref{theo-main}, we establish the following intermediate lemmata.

\begin{lemma}\label{rough-asy-X1-growth}
For each $k\geq 1$, there exists $C_k=C(n,k,g_0,\|h_0\|_{L^2_{f_0}},\chi_1)>0$ such that:
\begin{equation*}
\sup_Mf_0^{\frac{k-1}{2}}|\nabla^{g_0,\,k}X_1|_{g_0}\leq C_k.
\end{equation*}
\end{lemma}
\begin{proof}[Proof of Lemma \ref{rough-asy-X1-growth}]
Since $\Delta_{g_1}X_1+\Ric(g_1)(X_1,\cdot)^{\sharp}=0$ by [\eqref{equ:2bis}, Lemma \ref{id-sol-es}],
\begin{equation}\label{cov-der-F}
\begin{split}
\Delta_{g_0}X_1+\Ric(g_0)(X_1)&=\left(\Delta_{g_0}-\Delta_{g_1}\right)X_1+\left(\Ric(g_0)-\Ric(g_1)\right)(X_1)=:F,\\  
\sup_Mf_0^{\frac{k}{2}+1}|\nabla^{g_0,\,k}F|&<\infty,\quad \text{for all $k\geq 0$.}
\end{split}
\end{equation} 
Here we have used the exponential decay of $h$ and its derivatives provided by Proposition \ref{prop-alm-sharp-rate} to infer the decay on $F$.
Standard interior elliptic estimates applied to \eqref{cov-der-F} ensure the expected estimates on $|\nabla^{g_0,\,k}X_1|_{g_0}$ for each $k\geq 1$.
\end{proof}

\begin{lemma}\label{rough-asy-X1}
There exists a Killing field $\tilde{K}$ on the cone $C_0$ orthogonal to $r\partial_r$ and a sequence $r_i\rightarrow+\infty$ such that:
$$\lim_{r_i\rightarrow +\infty}r_i^{-1}\sup_{\{r_i^2/4\,\leq\, f_0\,\leq\, r_i^2/2\}}|X_1-X_0-\iota_*\tilde{K}|_{g_0}=0.$$
\end{lemma}
\begin{proof}[Proof of Lemma \ref{rough-asy-X1}]
We claim that the one-parameter family of vector fields $X_1(t):=\varphi_{t}^*X_1$ converges (sequentially) locally smoothly on $C_0$ as $t$ tends to $0$.

Indeed, as $|X_1|_{g_0}=O(f_0^{1/2})$, one has
$|X_1(t)|_{g_0(t)}=\sqrt{t}\varphi_t^*|X_1|_{g_0}\leq C(t\varphi_t^*f_0)^{1/2}\leq Cf_0^{1/2}$ for $0<t\leq1$ and  by Lemma \ref{rough-asy-X1-growth}, $|\nabla^{g_0(t),\,k}X_1(t)|_{g_0(t)}=t^{\frac{1-k}{2}}\varphi_t^*|\nabla^{g_0,\,k}X_1|_{g_0}\leq C_k(t\varphi_t^*f_0)^{\frac{1-k}{2}}\leq C(a,k)$ on $M_a$ for a given $a>0$ and $k\geq 1$. Here we are using the same notation as introduced before Theorem \ref{ac-str-exp}.

These facts show the aforementioned smooth local sequential convergence of $\iota^*X_1(t)$ on the cone $C_0$ as $t$ tends to $0^+$ by an application of the Arzela-Ascoli theorem.

 Let us denote $\tilde{X}_1$ any such limit and observe that the soliton equation implies that $\Li_{\tilde{X}_1}g_c=g_c=\Li_{r\partial_r}g_c$, i.e. the vector field $\tilde{X}_1-r\partial_r=:\tilde{K}$ is Killing. By Proposition \ref{claim-cone-killing-orth}, $\tilde{K}$ is necessarily orthogonal to $r\partial_r$. 

%\todo{TO: the claim below, see where we want to put it?},  

 Let us finally show that $\limsup_{r_i\rightarrow +\infty}r_i^{-1}\sup_{\{r_i^2/4\,\leq\, f_0\,\leq\, r_i^2/2\}}|X_1-X_0-\iota_*\tilde{K}|_{g_0}=0$ for some diverging sequence of positive radii $(r_i)_i$.
 
 Observe that since $\varphi_t^*X_0=X_0$ for all $t>0$ and $[X_0,\iota_*\tilde{K}]=0$,
% \begin{equation*}
 %\varphi_t^*[X_1,X_0]=\varphi_t^*(\nabla^{g_0}_{X_1}X_0-\nabla^{g_0}_{X_0}X_1)=\nabla^{g_0(t)}_{X_1(t)}X_0(t)-\nabla^{g_0(t)}_{X_0(t)}X_1(t)=\nabla^{g_0(t)}_{X_1(t)}X_0-\nabla^{g_0(t)}_{X_0}X_1(t),
% \end{equation*}
% any sequential limit of $|\varphi_t^*[X_1,X_0]|_{g_c}$ as $t$ goes to $0^+$ must be $0$ as $[\tilde{X_1},r\partial_r]=0$.
 %In particular, $|\partial_tX_1(t)|_{g_0(t)}=\frac{1}{\sqrt{t}}|[X_1,X_0]|_{g_0}$
  \begin{equation*}
 \sqrt{t}\varphi_t^*|X_1-X_0-\iota_*\tilde{K}|_{g_0}=|X_1(t)-X_0-\iota_*\tilde{K}|_{g_0(t)},
 \end{equation*}
so that on a fixed compact set of the form $\{1\leq f_0\leq 2\}$, there exists $t_i$ tending to $0$ such that $\sqrt{t_i}\varphi_{t_i}^*|X_1-X_0-\iota_*\tilde{K}|_{g_0}$ converges uniformly to $0$ on $\{1\leq f_0\leq 2\}$. By Theorem \ref{ac-str-exp}, it means that $r_i^{-1}|X_1-X_0-\iota_*\tilde{K}|_{g_0}$ converges to $0$ uniformly on $\{r_i^2/4\leq f_0\leq r_i^2/2\}$ if $r_i:=2 t_i^{-\frac{1}{2}}$.
\end{proof}

Let $K$ be a vector field given by Proposition \ref{can-kill-vec-field}
 and define the vector field on $M$, $$X:=X_1-X_0-K.$$  A priori, one has $g_0(X,X_0)=o(f_0)$ sequentially only thanks to Lemma \ref{rough-asy-X1} and Proposition \ref{can-kill-vec-field}.
\begin{lemma}\label{Lemma-improved-ortho}
One has the improved asymptotic orthogonality behavior:
\begin{equation*}
g_0(X,X_0)=O(f_0^{\frac{1}{2}}).
\end{equation*}
\end{lemma}

\begin{proof}[Proof of Lemma \ref{Lemma-improved-ortho}]
Let us derive a non-homogeneous first-order linear ordinary differential equation satisfied by $g_0(X,X_0)$ as follows:
\begin{equation}
\label{ode-scal-prod-X-X0}
\begin{split}
X_0\cdot g_0(X,X_0)&=X_0\cdot g_0(X_1-X_0,X_0)-X_0\cdot g_0(K,X_0)\\
&=g_0(\nabla^{g_0}_{X_0}(X_1-X_0),X_0)+g_0(X_1-X_0,\nabla^{g_0}_{X_0}X_0) +X_0\cdot\div_{g_0}K\\
%&=g_0([X_0,X_1-X_0],X_0)+g_0(\nabla^{g_0}_{X_1-X_0}X_0,X_0)+g_0(X_1-X_0,\nabla^{g_0}_{X_0}X_0)\\
%&\quad +X_0\cdot\div_{g_0}K\\
&=\frac{1}{2}\Li_{X_1-X_0}g_0(X_0,X_0)+g_0(X_1-X_0,\nabla^{g_0}_{X_0}X_0)+X_0\cdot\div_{g_0}K\\
&=\frac{1}{2}\Li_{X_1-X_0}g_0(X_0,X_0)+\frac{1}{2}g_0(X_1-X_0,X_0)\\
&\quad+\Ric(g_0)(X_1-X_0,X_0) +X_0\cdot\div_{g_0}K\\
&=\frac{1}{2}g_0(X,X_0)+\frac{1}{2}\Li_{X_1-X_0}g_0(X_0,X_0)-\frac{1}{2}\div_{g_0}K\\
&\quad-\frac{1}{2}g_0(X_1-X_0,\nabla^{g_0}\R_{g_0}) +X_0\cdot\div_{g_0}K.
\end{split}
\end{equation}
Here we have used the property that $\div_{f_0}K=0$ from Proposition \ref{can-kill-vec-field} in the second and last equalities, together with the soliton equation \eqref{egs-eqn} satisfied by $g_0$ and $X_0$ in the forth equality. Finally, we have made use of the soliton identity [\eqref{equ:2}, Lemma \ref{id-sol-egs}] in the last line.

According to Proposition \ref{can-kill-vec-field}, $X_0\cdot\div_{g_0}K=O(f_0^{1-\frac{n}{2}}e^{-f_0})$ and $\div_{g_0}K=O(f_0^{-\frac{n}{2}}e^{-f_0})$. 

Equation \eqref{hello-eqn-non-linear} and Proposition \ref{prop-alm-sharp-rate} ensure that $|\Li_{X_1-X_0}g_0(X_0,X_0)|_{g_0}=O(e^{-\gamma f_0})$ for all $\gamma\in(0,1)$. Finally, the asymptotic properties of $g_0$ and $X_0$ together with the assumption on the linear growth of $X_1$ give us the worst decay $g_0(X_1-X_0,\nabla^{g_0}\R_{g_0})=O(f_0^{-1})$. Therefore, integrating \eqref{ode-scal-prod-X-X0} along the flow lines of $X_0$ gives us the required growth at infinity of $g_0(X,X_0)$. 
\end{proof}

\begin{lemma}\label{lemma-ell-pde-scal-prod}
The scalar product $S:=g_0(X,X_0)$ satisfies the following non-homogeneous elliptic partial differential equation:
\begin{equation*}
\Delta_{g_0}S=g_0(X,\nabla^{g_0}\R_{g_0})+\div_{g_0}\Li_Xg_0(X_0)-g_0(\nabla^{g_0}\div_{g_0}X,X_0)+\frac{1}{2}g_0(\Li_Xg_0,\Li_{X_0}g_0).
\end{equation*}
In particular, $\sup_Mf_0^{\frac{k-1}{2}}|\nabla^{g_0,\,k}S|_{g_0}<\infty$ for all $k\geq 0$.
\end{lemma}
\begin{proof}[Proof of Lemma \ref{lemma-ell-pde-scal-prod}]
A direct computation shows that:
\begin{equation*}
\begin{split}
\Delta_{g_0}S&=g_0(\Delta_{g_0}X,X_0)+2\langle\nabla^{g_0}X,\nabla^{g_0}X_0\rangle_{g_0}+g_0(X,\Delta_{g_0}X_0)\\
&=g_0(\Delta_{g_0}X,X_0)+\frac{1}{2}\langle\Li_Xg_0,\Li_{X_0}g_0\rangle_{g_0}+g_0(X,\Delta_{g_0}X_0)\\
&=g_0(\Delta_{g_0}X,X_0)+\frac{1}{2}\langle\Li_Xg_0,\Li_{X_0}g_0\rangle_{g_0}-\Ric(g_0)(X,X_0),
\end{split}
\end{equation*}
where we have made use of the gradientness of $X_0$ (with respect to $g_0$) in the second line together with the soliton identity $\Delta_{g_0}X_0+\Ric(g_0)(X_0)=0$ from Lemma \ref{id-sol-es} in the last line.

Recall the following Bochner type formula that holds for an arbitrary vector field $V\in C^3_{loc}$:
\begin{equation*}
\div_{g_0}\Li_Vg_0-g_0(\nabla^{g_0}\div_{g_0}V,\cdot)=g_0(\Delta_{g_0}V,\cdot)+\Ric(g_0)(V,\cdot),
\end{equation*}
which, once applied to $V:=X$, leads to,
\begin{equation*}
\begin{split}
\Delta_{g_0}S&=-2\Ric(g_0)(X,X_0)+\div_{g_0}\Li_Xg_0(X_0)-g_0(\nabla^{g_0}\div_{g_0}X,X_0)+\frac{1}{2}\langle\Li_Xg_0,\Li_{X_0}g_0\rangle_{g_0}\\
&=\underbrace{g_0(X,\nabla^{g_0}\R_{g_0})+\div_{g_0}\Li_Xg_0(X_0)-g_0(\nabla^{g_0}\div_{g_0}X,X_0)+\frac{1}{2}\langle\Li_Xg_0,\Li_{X_0}g_0\rangle_{g_0}}_{=:F},
\end{split}
\end{equation*}
where we have made use of the soliton identity \eqref{equ:2} from Lemma \ref{id-sol-egs} in the last line.

Now, if we denote by $F$ the righthand side of the previous estimate, \eqref{hello-eqn-non-linear} together with Proposition \ref{prop-alm-sharp-rate}  combined with Proposition \ref{can-kill-vec-field} let us get that $$\sup_Mf_0^{1+\frac{k}{2}}|\nabla^{g_0,\,k}F|_{g_0}<\infty,\quad k\geq 0.$$ Observe that the worse decay comes from the first term $g_0(X,\nabla^{g_0}\R_{g_0}) $ a priori. Interior estimates for elliptic equations of the form $\Delta_{g_0}S=F$ with $g_0$ asymptotically conical gives us the required estimates on the covariant derivatives of $S$ since $S$ grows at most linearly by Lemma \ref{Lemma-improved-ortho}.
\end{proof}
The next lemma shows that $X$ is bounded on $M$ together with its covariant derivatives. Notice that the estimates below are qualitative, i.e. they are not quantified in terms of the smallness of $\|h_0\|_{L^2_f}$ yet.
\begin{lemma}\label{lemma-evo-eqn-X}
The vector field $X$ satisfies:
\begin{equation*}
\begin{split}
\sup_M&f_0^{\frac{k}{2}+1}\left|\nabla^{g_0,\,k}\left(\Delta_{g_0}X+\Ric(g_0)(X)\right)\right|_{g_0}<\infty,\quad k\geq 0,\\
\Delta_{f_0}|X|_{g_0}-\frac{|X|_{g_0}}{2}&\geq-\left(|\nabla^{g_0}(\div_{f_0}X)|+|\div_{g_0}\left(\Li_{X_1-X_0}g_0\right)|_{g_0}+|\Li_{X_1-X_0}g_0|_{g_0}|X_0|_{g_0}\right).
\end{split}
\end{equation*}
In particular, for all $k\geq 0$, $\sup_M|f_0^{\frac{k}{2}}\nabla^{g_0,\,k}X|_{g_0}<\infty$.
\end{lemma}
\begin{proof}[Proof of Lemma \ref{lemma-evo-eqn-X}] 

Thanks to the Bochner identity, 
\begin{equation}\label{ell-eqn-X}
g_0(\Delta_{g_0}X,\cdot)+\Ric(g_0)(X,\cdot)=\div_{g_0}\left(\Li_{X}g_0\right)-d\left(\div_{g_0}X\right).
\end{equation}
The $k$-th covariant derivatives of the righthand side of \eqref{ell-eqn-X} decays faster than polynomially thanks to \eqref{hello-eqn-non-linear}, Proposition \ref{prop-alm-sharp-rate} and Proposition \ref{can-kill-vec-field}. This implies the expected first estimate.

Therefore, as $\Delta_{f_0}K-\frac{K}{2}=0$ and $\div_{f_0}K=0$ thanks to Proposition \ref{can-kill-vec-field}, Lemma \ref{lem-bochner-lie-der} gives:
\begin{equation}
\begin{split}\label{beauty-X-ell-wei-eqn}
g_0\left(\Delta_{f_0}X-\frac{X}{2}+\nabla^{g_0}(\div_{f_0}X),\cdot\right)&=g_0\left(\Delta_{f_0}(X_1-X_0)-\frac{X_1-X_0}{2},\cdot\right)\\
&\quad +d(\div_{f_0}(X_1-X_0))\\
&=\div_{f_0}\Li_{X_1-X_0}g_0.
\end{split}
\end{equation}
Therefore,
\begin{equation*}
\begin{split}
\Delta_{f_0}&|X|^2_{g_0}-|X|^2_{g_0}\geq \\
&2|\nabla^{g_0}X|^2_{g_0}-2\left(|\nabla^{g_0}(\div_{f_0}X)|+|\div_{g_0}\left(\Li_{X_1-X_0}g_0\right)|_{g_0}\right)|X|_{g_0}+2\Li_{X_1-X_0}g_0(X_0,X)\\
&\geq 2|\nabla^{g_0}X|^2_{g_0}-2\left(|\nabla^{g_0}(\div_{f_0}X)|+|\div_{g_0}\left(\Li_{X_1-X_0}g_0\right)|_{g_0}+|\Li_{X_1-X_0}g_0|_{g_0}|X_0|_{g_0}\right)|X|_{g_0}.
\end{split}
\end{equation*}
This in turn implies by Kato's inequality:
\begin{equation*}
\begin{split}
\Delta_{f_0}|X|_{g_0}-\frac{1}{2}|X|_{g_0}&\geq  -\left(|\nabla^{g_0}(\div_{f_0}X)|+|\div_{g_0}\left(\Li_{X_1-X_0}g_0\right)|_{g_0}+|\Li_{X_1-X_0}g_0|_{g_0}|X_0|_{g_0}\right)\\
&\geq -C,
\end{split}
\end{equation*}
thanks to Lemma \ref{lemma-ell-pde-scal-prod} applied to $k=1$. Invoking Lemma \ref{rough-asy-X1}, i.e. that $\max_{ \{f_0\,=\, r_i^2/2\}}|X|_{g_0}=o(r_i)$ with and some diverging sequence of positive reals $(r_i)_i$ lets us get that $|X|_{g_0}=O(1)$ thanks to Lemma \ref{lemma-max-ppe-sub-lin} applied to $u:=|X|_{g_0}-2C$ and $C_0:=0$.
By the first series of estimates proved on the covariant derivatives of $\Delta_{g_0}X+\Ric(g_0)(X)$ together with the boundedness of $X$ lead to $\sup_M|f_0^{\frac{k}{2}}\nabla^{g_0,\,k}X|_{g_0}<\infty$ for all $k\geq 0$ by interior elliptic estimates. 
\end{proof}

\begin{lemma}\label{lemma-pt-dec-rhs}
For each $\gamma\in(0,1)$, there exists $C=C(n,\gamma,g_0,\chi_1)>0$ such that if $\|h_0\|_{L^2_{f_0}}\leq 1$, one has the following pointwise estimate:
\begin{equation*}
\left|\div_{f_0}\left(\Li_{X}g_0\right)\right|+\left|\nabla^{g_0}\div_{f_0}\left(\Li_{X}g_0\right)\right|+\left|\div_{f_0}\div_{f_0}\left(\Li_{X}g_0\right)\right|\leq C\|h_0\|_{L^2_{f_0}} e^{-\gamma f_0}.
\end{equation*}
\end{lemma}
%Notice that the previous claim is qualitative, i.e. the main estimate is not quantified in terms of the smallness of $\|h_0\|_{L^2_f}$ yet.
\begin{proof}[Proof of Lemma \ref{lemma-pt-dec-rhs}]
According to \eqref{hello-eqn-non-linear}, one gets:
\begin{equation*}
\begin{split}
\div_{f_0}\left(\Li_{X_1-X_0}g_0\right)&=\div_{f_0}\left(-L_{f_0}h+\Li_{B_{g_0}(h)^{\sharp}}g_0-\Li_{X_1-X_0}h-Q_{g_0}(h)\right).
\end{split}
\end{equation*}
Now, by Proposition \ref{can-kill-vec-field} and Lemma \ref{lem-bochner-lie-der}, $\div_{f_0}\left(\Li_{X}g_0\right)=\div_{f_0}\left(\Li_{X_1-X_0}g_0\right)$.
Then repeated applications of Proposition \ref{prop-alm-sharp-rate} give us the desired estimate as the vector field $X_0$ (respectively $X_1$) and its covariant derivatives have polynomial growth thanks to the boundedness of the curvature of $g_0$ via the soliton equation (respectively thanks to Lemma \ref{rough-asy-X1-growth}).
\end{proof}

The next lemma consists in deriving a uniform decay on the weighted divergence of $X$ (or equivalently of the orthogonal projection of $X$ in the radial direction $X_0$):
\begin{lemma}\label{lemma-pt-est-wei-div}
For each $\gamma\in(0,1)$, there exists $C=C(n,\gamma,g_0,\chi_1)>0$ such that if $\|h_0\|_{L^2_{f_0}}\leq 1,$ then on $M$, $$|\div_{f_0}X|+|\nabla^{g_0}\div_{f_0}X|_{g_0}\leq C\|h_0\|_{L^2_{f_0}}^{\frac{1}{2}}e^{-\gamma f_0}.$$
\end{lemma}
\begin{proof}[Proof of Lemma \ref{lemma-pt-est-wei-div}]
Observe first that the one-parameter family of vector fields $\iota^*X(t):=t^{-\frac{1}{2}}\iota^*\varphi_{t}^*X$, $t>0$, subconverges locally smoothly on $C_0$ as $t$ tends to $0$.

Indeed, $|X(t)|_{g_0(t)}=\varphi_t^*|X|_{g_0}\leq C$ for $t>0$ and $|\nabla^{g_0(t),\,k}X(t)|_{g_0(t)}=t^{-\frac{k}{2}}\varphi_t^*|\nabla^{g_0,\,k}X|_{g_0}\leq C_k(t\varphi_t^*f_0)^{-\frac{k}{2}}\leq C(a,k)$ on $M_a$ for a given $a>0$ by Lemma \ref{lemma-evo-eqn-X}. Here we are using the same notations as those introduced before Theorem \ref{ac-str-exp}.

These facts show the aforementioned smooth local sequential convergence of $\iota^*X(t)$ as $t$ tends to $0^+$ by an application of the Arzela-Ascoli theorem.

Denote by $\tilde{X}_0$ such a limit vector field on $C_0$ and observe that it induces a Killing vector field as for $\gamma\in(0,1)$,
\begin{equation}\label{decay-lie-der}
|\Li_{X(t)}g_0(t)|_{g_0(t)}=t^{-\frac{1}{2}}\varphi_t^*|\Li_{X}g_0|_{g_0}\leq Ct^{-\frac{1}{2}}\varphi_t^*\left(e^{-\gamma f_0}\right)\leq Ct^{-\frac{1}{2}}e^{-\frac{C(a,\gamma)}{t}},\end{equation}
on $M_a$ for some positive constant $C(a,\gamma)$ thanks to \eqref{hello-eqn-non-linear}, Proposition \ref{prop-alm-sharp-rate} and Proposition \ref{can-kill-vec-field}. By Proposition \ref{claim-cone-killing-orth}, $\tilde{X}_0$ is orthogonal to $r\partial_r$. In particular, any sequential limit of $g(t)(X(t),X_0(t))$ is $0$ so that this implies that $g_0(X,X_0)=o(f_0^{1/2})$ as $g(t)(X(t),X_0(t))=t^{1/2}\varphi_t^*g_0(X,X_0)$.

Thanks to Lemma \ref{lem-bochner-lie-der}, one can apply Lemma \ref{lem-equiv-gauge-eqn} with $g_0(Z,\cdot):=\div_{f_0}\Li_{X_1-X_0}g_0$ and invoke
Lemma \ref{lemma-pt-dec-rhs} for a given $\gamma\in(0,1)$ to give the existence of a positive constant $C$ such that on $M$:
\begin{equation}
\begin{split}\label{decouple-div-eqn-bis}
\left|\Delta_{f_0}\div_{f_0}X-\frac{\div_{f_0}X}{2}\right|=\left|\Delta_{f_0}\div_{f_0}(X_1-X_0)-\frac{\div_{f_0}(X_1-X_0)}{2}\right|&\leq C\|h_0\|_{L^2_{f_0}}e^{-\gamma f_0},
\end{split}
\end{equation}
where we have used Proposition \ref{can-kill-vec-field} in the first equality. Therefore, the use of Lemma \ref{lemma-max-ppe-sub-lin} applied to $u:=\div_{f_0}X$ lets us get the required exponential decay on $\div_{f_0}X$ on $M$. A similar reasoning applied to $-\div_{f_0}X$ leads to the required exponential decay on $|\div_{f_0}X|$ on $M$.

In order to derive the decay on the gradient of $\div_{f_0}X$, we invoke the following basic interpolation inequality whose proof can be found for instance in \cite{Aubin} for $m\geq 2$ and $x\in M$:
\begin{equation}
\begin{split}\label{inter-ine}
&\sup_{B_{g_0}\left(x,\frac{\iota}{2} f_0(x)^{\frac{1}{2}}\right)}f_0(x)^{\frac{1}{2}}|\nabla^{g_0}\div_{f_0}X|_{g_0}\\
&\leq C\sup_{B_{g_0}(x,\iota f_0(x)^{\frac{1}{2}} )}|\div_{f_0}X|_{g_0}^{1-\frac{1}{m}}\cdot \sup_{B_{g_0}\left(x,\iota f_0(x)^{\frac{1}{2}}\right)}\bigg(\sum_{k=0}^mf_0(x)^{\frac{k}{2}}|\nabla^{g_0,\,k}\div_{f_0}X|_{g_0}\bigg)^{\frac{1}{m}},
\end{split}
\end{equation}
where $\inj_{g_0}(x)\geq \iota f_0(x)^{\frac{1}{2}}$ from some positive constant $\iota$ uniform in $x\in M$, a property that holds thanks to the asymptotics of $g_0$ provided by Theorem \ref{ac-str-exp}.

Now, by a repeated application of Lemmata \ref{rough-asy-X1-growth}  and \ref{lemma-evo-eqn-X}, covariant derivatives of $\div_{f_0}X=\div_{f_0}(X_1-X_0)$ are growing at most polynomially in $f_0$. This fact combined with the previous exponential decay on $\div_{f_0}X$ and \eqref{inter-ine} appled to $m=2$ leads to the desired quantitative exponential decay on the gradient of $\div_{f_0}X$.

 This ends the proof of the lemma.
%\todo[inline]{F:Why do we have $\Delta_{f_0} f_0^{1/2} \leq \tfrac{1}{2} f_0^{1/2}$? I don't completely see it from Lemma \ref{id-sol-egs}. A: it is now included in Lemma \ref{lemma-max-ppe-sub-lin}}
\end{proof}

The next lemma derives a quantitative energy estimate on $X$ independent of the growth of $|X_1|_{g_0}$ (quantified by the constant $\chi_1$), observe that this is the first instance where we use the assumption $\ker L_{f_0}\cap L^2_{f_0}=\{0\}$:
\begin{lemma}\label{lemma-ene-X-h_0-0}
The vector field $X\in L^2_{f_0}$ and there exist $C=C(n,g_0)>0$, $\varepsilon(n,g_0)>0$ such that if $\|h_0\|_{L^2_{f_0}}\leq \varepsilon(n,g_0)$,
\begin{equation*}
\|X\|_{L^2_{f_0}}+\|\nabla^{g_0}X\|_{L^2_{f_0}}\leq C\|h_0\|_{L^2_{f_0}}.
\end{equation*}
\end{lemma}
\begin{proof}[Proof of Lemma \ref{lemma-ene-X-h_0-0}]
Let us first show that $X\in L^2_{f_0}$. Applying Lemma \ref{lem-bochner-lie-der} to $Y:=X$ gives:
\begin{equation*}
\Delta_{f_0}|X|_{g_0}-\frac{1}{2}|X|_{g_0}\geq -|\nabla^{g_0}\div_{f_0}X|_{g_0}-|\div_{f_0}\Li_Xg_0|_{g_0}.
\end{equation*}
Now, invoking Lemmata \ref{lemma-pt-dec-rhs} and \ref{lemma-pt-est-wei-div} let us prove that for all $\gamma\in (0,1)$, there exists $C>0$ such that on $M$:
\begin{equation}\label{exp-dec-rhs-X-0}
|\nabla^{g_0}\div_{f_0}X|_{g_0}+|\div_{f_0}\Li_Xg_0|_{g_0}\leq Ce^{-\gamma f_0}.
\end{equation}

As $|X|_{g_0}$ is bounded on $M$ thanks to Lemma \ref{lemma-evo-eqn-X}
, one can apply Lemma \ref{lemma-max-ppe-sub-lin} to $u:=|X|_{g_0}$, which once combined with the exponential decay \eqref{exp-dec-rhs-X-0} for a given $\gamma\in(0,1)$ leads to a bound of the form $|X|_{g_0}\leq 3C e^{-\gamma f_0}$ on $M$. This proves in particular that $X\in L^2_{f_0}$ if we choose $\gamma\in(1/2,1)$.

Regarding the quantitative $L^2_{f_0}$-bounds on $X$, one can apply Lemma \ref{korn-type-id} to $Y:=X$:
\begin{equation}\label{oh-boy-korn}
\| \Li_Xg_0\|^2_{L^2_{f_0}}=2\|\nabla^{g_0}X\|^2_{L^2_{f_0}}+\|X\|^2_{L^2_{f_0}}+2\| \div_{f_0}X\|^2_{L^2_{f_0}},
\end{equation}
an energy identity that leaves us with estimating $|\Li_Xg_0|_{g_0}$ pointwise from above.

Since $\ker L_{f_0}\cap L^2_{f_0}=\{0\}$ by assumption, $\Li_Kg_0=0$ so that by \eqref{hello-eqn-non-linear}, one has pointwise on $M$:
\begin{equation*}
\begin{split}
|\Li_Xg_0|_{g_0}&\leq |L_{f_0}h|_{g_0}+C\left(|\nabla^{g_0,\,2}h|_{g_0}+|\nabla^{g_0}h|_{g_0}+|h|_{g_0}\right)+|\Li_Xh|_{g_0}\\
&\leq|L_{f_0}h|_{g_0}+C\left(|\nabla^{g_0,\,2}h|_{g_0}+(1+|X|_{g_0})|\nabla^{g_0}h|_{g_0}+(1+|\nabla^{g_0}X|_{g_0})|h|_{g_0}\right),
\end{split}
\end{equation*}
where we have used the identity $\Li_Xh=\nabla^{g_0}_Xh+h\ast \nabla^{g_0}X$ in the last line together with Lemma \ref{app-C-monster}.
Squaring and integrating the previous pointwise estimate gives:
\begin{equation}
\begin{split}\label{rough-L2-lie-der-0}
\|\Li_Xg_0\|_{L^2_{f_0}}&\leq\|L_{f_0}h\|_{L^2_{f_0}}+C\left(\sum_{k=0}^2\|\nabla^{g_0,\,k}h\|_{L^2_{f_0}}+\sum_{k=0}^1\|\nabla^{g_0,\,1-k}h\|_{C^0}\|\nabla^{g_0,\,k}X\|_{L^2_{f_0}}\right)\\
&\leq C\left(\|h_0\|_{L^2_{f_0}}+\|L_{f_0}h\|_{L^2_{f_0}}\right)+C\varepsilon\left(\|X\|_{L^2_{f_0}}+\|\nabla^{g_0}X\|_{L^2_{f_0}}\right)\\
&\leq C\|h_0\|_{L^2_{f_0}}+C\varepsilon\left(\|X\|_{L^2_{f_0}}+\|\nabla^{g_0}X\|_{L^2_{f_0}}\right),
\end{split}
\end{equation}
where we have invoked Proposition \ref{prop-alm-sharp-rate} in the last inequality and where $C=C(n,g_0)>0$ may vary from line to line. By choosing $0<\varepsilon\leq \varepsilon(n,g_0)$ small enough, an absorption principle obtained by feeding back the  energy identity \eqref{oh-boy-korn} with the help of the above integral estimate \eqref{rough-L2-lie-der-0} lead to the required estimate on $\|X\|_{L^2_{f_0}}$.
\end{proof}
We are now in a position to prove item \ref{item-non-grad}, Theorem \ref{theo-main}:
\begin{proof}[Proof of item \ref{item-non-grad}, Theorem \ref{theo-main}]
We are in a good position to conclude that $h_0\equiv 0$ provided $\|h_0\|_{L^2_f}\leq \varepsilon(n,g_0)$ is chosen sufficiently small: thanks to \eqref{hello-eqn-non-linear} and $\Li_Kg_0=0$,
\begin{equation*}
\begin{split}
\lambda_+\|h_0\|_{L^2_f}^2&\leq \langle L_{f_0}h_0,h_0\rangle_{L^2_{f_0}}=\langle L_{f_0}h,h_0\rangle_{L^2_{f_0}}\\
&=\langle \Li_{B_{g_0}(h)^{\sharp}-(X_1-X_0)}g_0-\Li_{X_1-X_0}h-Q_{g_0}(h),h_0\rangle_{L^2_{f_0}}=:I-II-III.
\end{split}
\end{equation*}

Observe that $I$, $II$ and $III$ are well-defined in the Lebesgue sense thanks to Lemma \ref{lem-dec-eigen-lic} (applied to $h:=h_0$), Proposition \ref{prop-alm-sharp-rate} and Lemma \ref{lemma-evo-eqn-X}.

Now, for $R>0$ large enough, $\{f_0=R\}$ is a compact connected smooth hypersurface and let $d\sigma_{f_0}$ denote the weighted measure on $\{f_0=R\}$ induced by $d\mu_{f_0}$. Then an integration by parts shows that:
\begin{equation}
\begin{split}\label{est-I}
I=\lim_{R\rightarrow+\infty}\Bigg(-\int_{\{f_0\,\leq\, R\}} &2\div_{f_0}h_0( B_{g_0}(h)^{\sharp}-(X_1-X_0))\,d\mu_{f_0}\\
&+\int_{\{f_0=R\}}h_0\left( B_{g_0}(h)^{\sharp}-(X_1-X_0),\frac{X_0}{|X_0|_{g_0}}\right)\,d\sigma_{f_0}\Bigg)=0,
\end{split}
\end{equation}
thanks to Lemma \ref{lemma-div-free}. Notice that the limit of the second boundary term vanishes as the it is of magnitude $f_0^{\frac{n-1}{2}}e^{f_0}|h_0( B_{g_0}(h)^{\sharp}-(X_1-X_0),X_0/|X_0|_{g_0})=O(f_0^{-\lambda_+})=o(1).$
%\todo[inline]{F: Shouldn't this be $O(f_0^{-\lambda_+/2})$? A: I think we agreed on $O(f_0^{-\lambda_+})$ last time we spoke }

As for the integral $II$, observe that pointwise on $M$:
\begin{equation}\label{pt-id-lie-der}
\begin{split}
\langle\Li_{X_1-X_0}h,h_0\rangle_{g_0}&=\langle\Li_{X}h,h_0\rangle_{g_0}+\langle\Li_{K}h,h_0\rangle_{g_0}\\
&=\langle\Li_{X}h,h_0\rangle_{g_0}+K\cdot \langle h,h_0\rangle_{g_0}-\langle h,\Li_{K}h_0\rangle_{g_0}.
\end{split}
\end{equation}
Here we have used that $\Li_Kg_0=0$ in the second line.

Now, on the one hand,
\begin{equation*}%\label{pt-id-lie-der-bis}
\begin{split}
\left|\langle\Li_{X}h,h_0\rangle_{g_0}\right|&\leq \left(|\nabla^{g_0}h|_{g_0}|X|_{g_0}+|h|_{g_0}|\nabla^{g_0}X|_{g_0}\right)|h_0|_{g_0},
\end{split}
\end{equation*}
 pointwise on $M$. Here we have invoked the formula $\Li_YT(U,V)=\nabla^{g_0}_YT(U,V)+T(\nabla^{g_0}_UY,V)+T(U,\nabla^{g_0}_VY)$ that holds for any vector field $U$, $V$ and $Y$ and an arbitrary $(0,2)$-tensor $T$ on $M$.

On the other hand, since $K$ is a Killing field with respect to $g_0$ and as $\div_{f_0}K=0$ by Proposition \ref{can-kill-vec-field}, one infers that $K\cdot f_0=\langle K,X_0\rangle_{g_0}=0$, i.e. the flow of $K$ preserves $f_0$. This implies that $\Li_Kh_0$ belongs to $E_+$ so that one gets:
\begin{equation}
\begin{split}\label{simplified-exp}
\langle h,\Li_{K}h_0\rangle_{L^2_{f_0}}&=\langle h_0,\Li_{K}h_0\rangle_{L^2_{f_0}}\\
&=\int_M\frac{K}{2}\cdot |h_0|^2_{g_0}\,d\mu_{f_0}\\
&=-\frac{1}{2}\int_M(\div_{f_0}K)|h_0|^2_{g_0}\,d\mu_{f_0}=0,
\end{split}
\end{equation}
since $\div_{f_0}K=0$. Here we have integrating by parts in the last line.

A similar integration by parts shows that
\begin{equation}
\begin{split}\label{simplified-exp-bis}
\int_MK\cdot \langle h,h_0\rangle_{g_0}\,d\mu_{f_0}&=-\int_M(\div_{f_0}K)\langle h_0,h\rangle_{g_0}\,d\mu_{f_0}=0.
\end{split}
\end{equation}

Therefore, thanks to \eqref{pt-id-lie-der}, \eqref{simplified-exp} and \eqref{simplified-exp-bis}: 
\begin{equation}
\begin{split}\label{est-II}
|II|&=\left|\langle\Li_{X}h,h_0\rangle_{L^2_{f_0}}\right|\leq \int_M\left(|\nabla^{g_0}h|_{g_0}|X|_{g_0}+|h|_{g_0}|\nabla^{g_0}X|_{g_0}\right)|h_0|_{g_0}\,d\mu_{f_0}\\
&\leq C\left(\|h\|_{C^0}+\|\nabla^{g_0}h\|_{C^0}\right)\left(\|X\|_{L^2_{f_0}}^2+\|\nabla^{g_0}X\|_{L^2_{f_0}}^2+\|h_0\|_{L^2_{f_0}}^2\right)\\
&\leq C\|h_0\|_{L^2_{f_0}}\left(\|X\|_{L^2_{f_0}}^2+\|\nabla^{g_0}X\|_{L^2_{f_0}}^2+\|h_0\|_{L^2_{f_0}}^2\right)\\
&\leq C\|h_0\|_{L^2_{f_0}}^3,
\end{split}
\end{equation}
 where $C=C(n,g_0)>0$ is a constant that may vary from line to line. Here, Proposition \ref{prop-alm-sharp-rate} is being used in the third inequality and Lemma \ref{lemma-ene-X-h_0-0} has been used in the last line.
 
 A similar estimate holds on $III$ thanks to the definition of $Q_{g_0}(h)$  and the use of Proposition \ref{prop-alm-sharp-rate}: there exists $C=C(n,g_0)>0$ such that,
 \begin{equation}\label{est-III}
 |III|\leq C\|h_0\|_{L^2_{f_0}}^3.
 \end{equation}
Therefore, the combination of \eqref{est-I}, \eqref{est-II} and \eqref{est-III} ensures the existence of  $C=C(n,g_0)>0$ such that if $\|h_0\|_{L^2_f}\leq \varepsilon(n,g_0)$,
 \begin{equation*}
 \lambda_+\|h_0\|_{L^2_f}^2\leq C\|h_0\|_{L^2_f}^3,
 \end{equation*}
which implies $h_0=0$ provided $\varepsilon(n,g_0)$ is chosen small enough compared to $C^{-1}\cdot\lambda_+$. 

Now, this implies that $g_1=g_0$ so that \eqref{eqn-sol-g_0} ensures that the vector field $X_1-X_0$ is a Killing field (with respect to the metric $g_0$). Thanks to Lemma \ref{lem-bochner-lie-der}, $\Delta_{f_0}Y-\frac{1}{2}Y+\nabla^{g_0}\div_{f_0}Y=0$ on $M$ with $Y:=X_1-X_0$. Therefore, Lemma \ref{lem-equiv-gauge-eqn} guarantees that $\div_{f_0}X$ (since $\div_{f_0}K=0$) lies in the kernel of $\Delta_{f_0}-\frac{1}{2}$. As $\div_{f_0}X$ is exponentially decaying on $M$ by Lemma \ref{lemma-pt-est-wei-div}, Lemma \ref{lemma-max-ppe-sub-lin} applied to $C_0:=0$ leads to $\div_{f_0}X=0$ on $M$ so that $\Delta_{f_0}X=X/2$ by invoking Lemma \ref{lem-bochner-lie-der} again. Moreover, Lemma \ref{rough-asy-X1} and Lemma \ref{lemma-max-ppe-sub-lin} applied to $u:=|X|_{g_0}$ and $C_0:=0$ lead to $X=0$ on $M$, i.e. $X_1-X_0=K$. Finally, since $\div_{f_0}K=0$ by Proposition \ref{can-kill-vec-field}, the previously established fact $\div_{f_0}(X_1-X_0)=0$ together with the fact that $X_1-X_0$ is Killing (implying that $\div_{g_0}(X_1-X_0)=0$ on $M$) lead us to prove that $K$ (and hence $X_1-X_0$) is orthogonal to $X_0$ pointwise on $M$.
\end{proof}

\begin{rk}\label{rk-non-self-sim}
If $(M^n,g_0,\nabla^{g_0}f_0)$ is an asymptotically conical unstable expanding gradient Ricci soliton which does not satisfy the nullity condition $\ker L_{f_0}\cap L^2_{f_0}=\{0\}$ a priori, one can still prove the following statement by inspecting the above proof of item \ref{item-non-grad}, Theorem \ref{theo-main}:
Given $\chi_1>0$, there exists $\varepsilon=\varepsilon(n,g_0,\chi_1)>0$ such that if $(g(\tau))_{\tau\,\leq\, 0}$ is a solution provided by Theorem \ref{theo-ancient-sol} satisfying $\|h_0\|_{L^2_f}\leq \varepsilon$ then $g+h(0)$ is not an expanding Ricci soliton $(M^n,g_1,X_1)$ such that $|X_1|_{g_0}\leq \chi_1f_0^{1/2}$ unless $h(0)\equiv0$, i.e. $g(\tau)=g_0$ for all $\tau\leq 0$ and $X_1-\nabla^{g_0}f_0$ is a Killing field orthogonal to $\nabla^{g_0}f_0$.
\end{rk}

\begin{conj}\label{conj-uniform-eps}
The constant $\varepsilon$ in Remark \ref{rk-non-self-sim} can be 
chosen independent of $\chi_1$; that is, item \eqref{item-non-grad} 
of Theorem \ref{theo-main} holds without the assumption 
$\ker L_{f_0}\cap L^2_{f_0}=\{0\}$, with $\varepsilon=\varepsilon(n,g_0)$.
\end{conj}

A natural strategy is to compare the two solitons via the time-dependent diffeomorphism
\[
\Psi_t:=\varphi_t^1\circ\psi_t^{-1}\circ\varphi_{1/t}^0,\qquad t\in(0,1],
\]
where $(\varphi_t^i)_{t>0}$ denotes the flow of $-X_i/t$ for $i=0,1$ and $(\psi_t)_{t\in(0,1]}$ is the DeTurck diffeomorphism of \eqref{corr-rf-sol}. One then has  $\Psi_t^{\ast}g_1=g_0+h(\log t)$. Applying the proof of item \eqref{item-non-grad} to $\Psi_t^{\ast}g_1$ in place of $g_1$ suggests that the conjecture would follow from a good enough estimate of $|\Psi_t^{\ast}f_0-f_0|$ as $f_0\to\infty$, uniformly in  $t\in(0,1]$ and independent of $\chi_1$.

\subsection{Proof of item \eqref{item-grad}, Theorem \ref{theo-main}}
The proof of item \eqref{item-grad}, Theorem \ref{theo-main} follows that of item \eqref{item-non-grad}, Theorem \ref{theo-main} very closely and it is simpler in that there is no need of considering a canonical extension of a Killing vector field given by Proposition \ref{can-kill-vec-field}. For the above reasons, we will be sketchy. 
\begin{proof}[Proof of item \eqref{item-grad}, Theorem \ref{theo-main}]
With the notations of Theorem \ref{theo-main}, let $(M^n,g_1,X_1=\nabla^{g_1}f_1)$ be an expanding gradient Ricci soliton with $g_1=g_0+h(0)=:g_0+h$, $\pi_+(h)=h_0$ and $\|h_0\|_{L^2_{f_0}}\leq \varepsilon$. Define $X:=X_1-X_0$ so that:
\begin{equation}
\begin{split}\label{hello-eqn-non-linear-bis}
L_{f_0}h&=\Li_{B_{g_0}(h)^{\sharp}-X}g_0-\Li_{X}h-Q_{g_0}(h),
\end{split}
\end{equation}
similarly to \eqref{hello-eqn-non-linear}.

Now, thanks to the soliton equation both satisfied by $(g_1,\nabla^{g_1}f_1)$ and $(g_0,\nabla^{g_0}f_0)$,
\begin{equation*}
\begin{split}
\nabla^{g_0}X&=\nabla^{g_1,\,2}f_1-\nabla^{g_0,\,2}f_0+\left(\nabla^{g_0}-\nabla^{g_1}\right)\nabla^{g_1}f_1\\
&=\Ric(g_1)-\Ric(g_0)+\frac{h}{2}+g_1^{-1}\ast \nabla^{g_0}h\ast \nabla^{g_1}f_1.
\end{split}
\end{equation*}
In particular, an application of Lemma \ref{app-C-monster} together with Proposition \ref{prop-alm-sharp-rate} shows that pointwise on $M$:
\begin{equation*}
\begin{split}
|\nabla^{g_0}X|_{g_0}&\leq C\left(|\nabla^{g_0,\,2}h|_{g_0}+|h|_{g_0}+|\nabla^{g_0}h|_{g_0}|X_1|_{g_0}\right)\\
&\leq C\left(|\nabla^{g_0,\,2}h|_{g_0}+|h|_{g_0}+|\nabla^{g_0}h|_{g_0}\left(|X|_{g_0}+f_0^{1/2}\right)\right),
\end{split}
\end{equation*}
This shows that the function $\log(|X|^2+1)^{1/2}$ has an exponentially decaying (with respect to $f_0$) gradient on $M$, i.e. there exists $C=C(n,g_0)>0$ such that on $M$, $|\nabla^{g_0}\log(|X|^2+1)^{1/2}|_{g_0}\leq Ce^{-f_0/2}$. In particular, the mean value inequality shows that $|X|$ is bounded on $M$.

Then an adaptation of the proof of Lemmata \ref{lemma-pt-dec-rhs} and \ref{lemma-pt-est-wei-div} lets us prove that for all $\gamma\in (0,1)$, there exists $C>0$ such that on $M$,
\begin{equation}
\begin{split}\label{exp-decay-grad-case-qual}
\left|\div_{f_0}\left(\Li_{X}g_0\right)\right|+\left|\nabla^{g_0}\div_{f_0}\left(\Li_{X}g_0\right)\right|+\left|\div_{f_0}\div_{f_0}\left(\Li_{X}g_0\right)\right|&\leq Ce^{-\gamma f_0},\\
|\div_{f_0}X|+|\nabla^{g_0}\div_{f_0}X|_{g_0}&\leq Ce^{-\gamma f_0},
\end{split}
\end{equation}
where $C$ here may depend on $\sup_M|X|_{g_0}$.

The estimates \eqref{exp-decay-grad-case-qual} allow to adapt the proof almost verbatim of Lemma \ref{lemma-ene-X-h_0-0} to get the existence of $C=C(n,g_0)>0$ such that for $\|h_0\|_{L^2_{f_0}}\leq \varepsilon(n,g_0)$,
\begin{equation}\label{l2-est-quan-X}
\|X\|_{L^2_{f_0}}+\|\nabla^{g_0}X\|_{L^2_{f_0}}\leq C\|h_0\|_{L^2_{f_0}}.
\end{equation}
From \eqref{l2-est-quan-X}, one can proceed as in the rest of the proof of  item \eqref{item-non-grad}, Theorem \ref{theo-main} to prove that $h_0=0$, $g_1=g_0$ and $X=0$ as required. 
%Notice that in order to handle $II$, one can skip \eqref{simplified-exp} by estimating the term $\langle\Li_{X}h,h_0\rangle_{g_0}$ directly by $C(|\nabla^{g_0}X|_{g_0}|h|_{g_0}|h_0|_{g_0}+|X|_{g_0}|\nabla^{g_0}h|_{g_0}|h_0|_{g_0})$ and invoke Proposition \ref{prop-alm-sharp-rate} to handle (the covariant derivatives of) $h$ pointwise.
\end{proof}

\section{Proof of Theorem \ref{theo-main}, item \eqref{item-mod-spa}}\label{sec-item-mod-spa}
For $h_1^+$ and $h_2^+$ in $E_+$ with $\|h_i^+\|_{L^2_f}<\varepsilon$, $i=1,2$, let $(g_i(t))_{t\in(0,1]}$ be the solutions to the Ricci flow associated to each $h_i^+$ as in [Theorem \ref{theo-main}, item \eqref{item-exi}]. Let $\varphi$ be a smooth diffeomorphism of $M$ such that for every $t\in(0,1]$, $g_2(t)=\varphi^*g_1(t)$. 

Now, evaluating at $t=1$ the previous equality gives:
\begin{equation}\label{fund-id}
g_0+h_2(0)=\varphi^*(g_0+h_1(0)),\quad \pi_+(h_i(0))=h_i^+,\quad i=1,2.
\end{equation} 
From now on, define $h_i:=h_i(0)$, $i=1,2$ and define the tensor 
\begin{equation}\label{eqn-defn-h-0-bis}
h_0:=h_2-\varphi^*h_1=\varphi^*g_0-g_0,
\end{equation}
Define the difference of the soliton vector fields as follows:
\begin{equation}\label{eqn-defn-X-bis}
X:=\nabla^{\varphi^*g_0}(\varphi^*f_0)-\nabla^{g_0}f_0=\varphi^*X_0-X_0,
\end{equation}
where $X_0:=\nabla^{g_0}f_0$.
Assume further that $\varphi$ preserves the zero set of the background soliton vector field $X_0$, i.e. if $E:=\{X_0=0\}$, then $\varphi(E)\subset E$. Recall that $E$ is a non-empty compact set by properness of $f_0$ and Lemma \ref{id-sol-egs}.

Finally, define the one-parameter family of symmetric $2$-tensors:
\begin{equation*}
h_0(\tau):=h_2(\tau)-\varphi^*h_1(\tau),\quad\tau\leq 0.
\end{equation*}
Observe that by \eqref{eqn-defn-h-0-bis}, $h_0(0)=h_2-\varphi^*h_1=h_0$. 

The following lemma serves as a suitable replacement of Lemma \ref{lem-dec-eigen-lic} and its proof makes heavy use of the assumption on $\varphi$, i.e. $\varphi(E)\subset E$.
\begin{lemma}\label{lemma-rough-equiv}
For each $k\geq 0$, there exists $C_k=C(n,k,g_0)>0$ and $\varepsilon_k=\varepsilon(n,k,g_0)>0$ such that if $\varepsilon\leq \varepsilon_k$, one has on $M\times(-\infty,0]$,
\begin{equation*}
|\nabla^{g_0,\,l}h_0(\tau)|_{g_0}+|\nabla^{g_0,\,l}\varphi^*h_1(\tau)|_{g_0}\leq C_k\,\varepsilon e^{\lambda_+\tau}e^{-f_0+C_0\sqrt{f_0}},\quad l\leq k,
\end{equation*}
and,
\begin{equation*}
|X|_{g_0}\leq C_0\varepsilon,\quad |\nabla^{g_0,\,l+1}X|_{g_0}\leq C_k\,\varepsilon e^{-f_0+C_0\sqrt{f_0}},\quad l\leq k.
\end{equation*}
\end{lemma}

\begin{proof}[Proof of Lemma \ref{lemma-rough-equiv}]
Regarding the first series of estimates, observe first that it suffices to prove that for all $k\geq 0$, there exists $C_k=C(n,k,g_0)>0$  and $\varepsilon_k=\varepsilon(n,k,g_0)>0$ such that if $\varepsilon\leq \varepsilon_k$, one has on $M\times(-\infty,0]$,
\begin{equation}\label{1-st-series-cov-der-h0}
|\nabla^{g_0,\,l}h_0(\tau)|_{g_0}\leq C_k\,\varepsilon e^{\lambda_+\tau} e^{-f_0+C_0\sqrt{f_0}},\quad l\leq k.
\end{equation}
Indeed, the triangular inequality together with \eqref{1-st-series-cov-der-h0} gives 
\begin{equation*}
|\nabla^{g_0,\,l}\varphi^*h_2(\tau)|_{g_0}\leq C_k\,\varepsilon e^{\lambda_+\tau} e^{-f_0+C_0\sqrt{f_0}},\quad l\leq k,
\end{equation*}
thanks to Proposition \ref{prop-alm-sharp-rate} applied to $h_1$.

Let us first prove a weaker version of \eqref{1-st-series-cov-der-h0} for at $\tau=0$: for all $k\geq 0$, there exists $C_k=C(n,k,g_0)>0$  and $\varepsilon_k=\varepsilon(n,k,g_0)>0$ such that if $\varepsilon\leq \varepsilon_k$, one has on $M$,
\begin{equation}\label{1-st-series-cov-der-h0-weak}
|\nabla^{g_0,\,l}h_0|_{g_0}\leq C_k\,\varepsilon e^{-(1-C_0\varepsilon)f_0},\quad l\leq k.
\end{equation}

Indeed, thanks to Proposition \ref{prop-alm-sharp-rate} applied both to $h_1$ and $h_2$ (at $\tau=0$), there exists $C=C(n,g_0)>0$ such that,
\begin{equation*}
h_0=\varphi^*g_0-g_0\leq C\|h_2\|_{L^2_{f_0}}g_0+C\|h_1\|_{L^2_{f_0}}\varphi^*g_0\leq C\varepsilon\left( g_0+\varphi^*g_0\right),
\end{equation*}
so that if $\varepsilon$ is small enough, $\varphi^*g_0\leq \frac{1+C\varepsilon}{1-C\varepsilon}g_0$ on $M$. Similarly, one gets $\varphi^*g_0\geq \frac{1-C\varepsilon}{1+C\varepsilon}g_0$ on $M$. As an intermediate result, we therefore get 
\begin{equation}\label{1-st-est-varphi}
|\varphi^*g_0-g_0|_{g_0}\leq C_0\,\varepsilon,\quad\text{ on $M$},
\end{equation}
 for some positive constant $C_0$. In particular, for all points $x$ and $y$ in $M$, 
 \begin{equation}\label{2-nd-est-varphi}
 \sqrt{1-C_0\varepsilon}\,d_{g_0}(x,y)\leq d_{\varphi^*g_0}(x,y)\leq \sqrt{1+C_0\varepsilon}\,d_{g_0}(x,y).
 \end{equation}
  Now, a Taylor expansion with integral remainder implies on $M$, 
  \begin{equation}\label{basic-pot-fct-f_0}
  \left|f_0(x)-\frac{d_{g_0}(p,x)^2}{4}\right|\leq C(1+d_{g_0}(p,x)), \quad \forall\,x\in M,
  \end{equation}
   for some positive constant $C$ and where $p$ is a critical point of $f_0$. This estimate uses the soliton equation \eqref{egs-eqn} and the quadratic Ricci curvature decay only: see \cite{Che-Der} for instance. By assumption, if $p$ is a critical point of $f_0$, $\varphi(p)$ is too. Therefore, evaluating the previous estimate \eqref{basic-pot-fct-f_0} at $\varphi(p)$ and $\varphi(x)$ with the help of \eqref{2-nd-est-varphi} gives:
 \begin{equation}
 \begin{split}\label{lower-bd-dist-varphi-f0}
 \varphi^*f_0(x)&\geq \frac{(1-C_0\varepsilon)\,d_{g_0}(p,x)^2}{4}-C( \sqrt{1+C_0\varepsilon}\,d_{g_0}(p,x)+1)\\
 &\geq (1-C'_0\varepsilon)f_0(x)-C'_0(\varepsilon),
 \quad \forall\,x\in M,
 \end{split}
 \end{equation}
%\todo[inline]{F: Would it be better to keep $\sqrt{f_0}$ on the RHS?}
 where $C'_0=C'_0(n,g_0)$ is a positive constant and where we have used \eqref{basic-pot-fct-f_0} together with Young's inequality in the last line.

Applying Proposition \ref{prop-alm-sharp-rate} both to $h_1$ and $h_2$ again and invoking \eqref{lower-bd-dist-varphi-f0} give:
\begin{equation*}
\begin{split}
|h_0|_{g_0}&=|\varphi^*g_0-g_0|_{g_0}\leq C\left(\|h_2\|_{L^2_{f_0}}f_0^{-\frac{n}{2}-\delta}e^{-f_0}+\|h_1\|_{L^2_{f_0}}\varphi^*f_0^{-\frac{n}{2}-\delta}e^{-\varphi^*f_0}\right)\\
&\leq C_0\varepsilon e^{-(1-C'_0\varepsilon)f_0}.
\end{split}
\end{equation*} 
%\todo[inline]{F: I assume this should give \eqref{1-st-series-cov-der-h0} for $k=0$ at $\tau =0$. But this is not really the claimed estimate, no? A: agreed, I added at the beginnign of the proof "a weaker version of ..."}
 Let us prove \eqref{1-st-series-cov-der-h0-weak} in the case $k=1$ only, higher covariant derivatives of $h_0$ can be handled similarly by induction on $k\geq 1$.
 
Indeed, the definition of $h_0$ and the triangular inequality give:
\begin{equation}
\begin{split}\label{cov-der-h0-1}
|\nabla^{g_0}h_0|_{g_0}&=|\nabla^{g_0}h_2-\nabla^{g_0}\varphi^*h_1|_{g_0}\\
&\leq |\nabla^{g_0}h_2|_{g_0}+|\nabla^{g_0}\varphi^*h_1|_{g_0}\\
&\leq |\nabla^{g_0}h_2|_{g_0}+|(\nabla^{g_0}-\nabla^{\varphi^*g_0})\varphi^*h_1|_{g_0}+|\nabla^{\varphi^*g_0}\varphi^*h_1|_{g_0}\\
&\leq |\nabla^{g_0}h_2|_{g_0}+|\varphi^*g_0^{-1}\ast \nabla^{g_0}(\varphi^*g_0-g_0)\ast \varphi^*h_1|_{g_0}+|\nabla^{\varphi^*g_0}\varphi^*h_1|_{g_0}\\
&\leq |\nabla^{g_0}h_2|_{g_0}+C|\nabla^{g_0}h_0|_{g_0} \varphi^*|h_1|_{g_0}+C\varphi^*|\nabla^{g_0}h_1|_{g_0}\\
&\leq C\varepsilon |\nabla^{g_0}h_0|_{g_0} +|\nabla^{g_0}h_2|_{g_0}+C\varphi^*|\nabla^{g_0}h_1|_{g_0},
\end{split}
\end{equation}
where $C=C(n,g_0)>0$ and where we have used the previous bound on $h_0$ together with \eqref{1-st-est-varphi} in the penultimate and last line. 

If $\varepsilon$ is chosen small enough (this choice made here reveals the dependency of $\varepsilon$ on the number $k$ of covariant derivatives of $h_0$), the previous estimate together with Proposition \ref{prop-alm-sharp-rate} lead to a weaker version of  the desired pointwise bound on $|\nabla^{g_0}h_0|_{g_0}$.

Now, observe that thanks to the soliton equation both satisfied by $(\varphi^*g_0,\nabla^{\varphi^*g_0}\varphi^*f_0)$ and $(g_0,\nabla^{g_0}f_0)$,
\begin{equation*}
\begin{split}
\nabla^{g_0}X&=\nabla^{\varphi^*g_0,\,2}\varphi^*f_0-\nabla^{g_0,\,2}f_0+\left(\nabla^{g_0}-\nabla^{\varphi^*g_0}\right)\nabla^{\varphi^*g_0}\varphi^*f_0\\
&=\Ric(\varphi^*g_0)-\Ric(g_0)+\frac{h_0}{2}+(\varphi^*g_0)^{-1}\ast \nabla^{g_0}h_0\ast \nabla^{\varphi^*g_0}\varphi^*f_0.
\end{split}
\end{equation*}

In particular, an application of Lemma \ref{app-C-monster} shows that pointwise on $M$:
\begin{equation*}
|\nabla^{g_0}X|_{g_0}\leq C\left(|\nabla^{g_0,\,2}h_0|_{g_0}+|h_0|_{g_0}+|\nabla^{g_0}h_0|_{g_0}\left(|\varphi^*(\nabla^{g_0}f_0)|_{g_0}+1\right)\right).
\end{equation*}
As $|\nabla^{g_0,\,k}h_0|_{g_0}\leq C_k\varepsilon e^{-(1-C_0\varepsilon)f_0}$ on $M$ for $k\geq 0$ thanks to the first series of estimates \eqref{1-st-series-cov-der-h0-weak} we proved above, there exists $C=C(n,g_0)>0$ such that,
\begin{equation*}
|\nabla^{g_0}X|_{g_0}\leq C\varepsilon e^{-(1-C_0\varepsilon) f_0}.
\end{equation*}
This proves a weaker version of the second series of estimates.

The mean value inequality together with the previously established estimate on $|\nabla^{g_0}X|_{g_0}$ give that $|X|_{g_0}$ is bounded on $M$. Moreover, by assumption on $\varphi$, $X(p)=0$ if $p$ lies in the zero set of $X_0=\nabla^{g_0}f_0$. In particular, one has pointwise on $M$,
\begin{equation*}
|X|_{g_0}\leq C\varepsilon,
\end{equation*}
for some positive constant $C=C(n,g_0)$.

This bound implies: 
\begin{equation*}
\begin{split}
|\nabla^{g_0}(\varphi^*f_0-f_0)|_{g_0}&\leq |(\nabla^{g_0}-\nabla^{\varphi^*g_0})\varphi^*f_0|_{g_0}+|X|_{g_0}\\
&\leq C|h_0|_{g_0}\varphi^*|X_0|_{g_0}+|X|_{g_0}\\
&\leq C\varepsilon.
\end{split}
\end{equation*}
A further application of the mean value inequality implies that $\varphi^*f_0-f_0$ grows at most linearly in the distance from a fixed point $p$ in the zero set of $X_0$ uniformly: there exists $C=C(n,g_0)>0$ such that $|\varphi^*f_0-f_0|\leq C\varepsilon\sqrt{f_0}$ on $M$ by the quadratic growth of $f_0$ recalled in \eqref{2-nd-est-varphi}.
%: observe that one could even get $|\varphi^*f_0-f_0|\leq C\varepsilon \sqrt{f_0}$ on $M$.
%but this does not improve the required estimates for higher covariant derivatives of $h_0(\tau)$ because of the polynomial powers $f_0^k$ in Proposition \ref{prop-alm-sharp-rate}.
Feeding back this estimate to the proof of \eqref{1-st-series-cov-der-h0-weak} that led to weaker versions of the desired estimates, this proves the expected estimates both on $h_0$ and $\nabla^{g_0}X$ and their covariant derivatives.

In order to prove \eqref{1-st-series-cov-der-h0} for any $\tau\geq 0$, it suffices to invoke Proposition \ref{prop-alm-sharp-rate} and the previously established bounds \eqref{1-st-series-cov-der-h0} at $\tau=0$. Indeed, for $k=0$, 
\begin{equation*}
\begin{split}
|h_0(\tau)|_{g_0}\leq |h_2(\tau)|_{g_0}+|\varphi^*h_1(\tau)|_{g_0}\leq |h_2(\tau)|_{g_0}+C\varphi^*|h_1(\tau)|_{g_0} \leq C_0\varepsilon e^{\lambda_+\tau}e^{-f_0+C_0\sqrt{f_0}}.
\end{split}
\end{equation*}
For $k=1$, one invokes a similar reasoning as that used to derive \eqref{cov-der-h0-1}:
\begin{equation*}
\begin{split}
|\nabla^{g_0}h_0(\tau)|_{g_0}&=|\nabla^{g_0}h_2(\tau)-\nabla^{g_0}\varphi^*h_1(\tau)|_{g_0}\\
&\leq |\nabla^{g_0}h_2(\tau)|_{g_0}+|\nabla^{g_0}\varphi^*h_1(\tau)|_{g_0}\\
&\leq |\nabla^{g_0}h_2(\tau)|_{g_0}+|(\nabla^{g_0}-\nabla^{\varphi^*g_0})\varphi^*h_1(\tau)|_{g_0}+|\nabla^{\varphi^*g_0}\varphi^*h_1(\tau)|_{g_0}\\
&\leq |\nabla^{g_0}h_2(\tau)|_{g_0}+C|\nabla^{g_0}h_0|_{g_0} \varphi^*|h_1(\tau)|_{g_0}+C\varphi^*|\nabla^{g_0}h_1(\tau)|_{g_0}\\
&\leq C_1\varepsilon e^{\lambda_+\tau}e^{-f_0+C_0\sqrt{f_0}},
\end{split}
\end{equation*}
as required.
 Higher covariant derivatives of $h_0(\tau)$ can be handled similarly by induction on $k\geq 1$.
\end{proof}

The following lemma establishes a system of elliptic-parabolic equations satisfied by the triple $(h_0(\tau))_{\tau\,\leq\,0}, h_0,X)$:
\begin{lemma}\label{lem-hello-eqn-non-linear-bis}
The following system holds:
\begin{equation*}
\begin{split}
\partial_{\tau}h_0(\tau)&=L_{f_0}h_0(\tau)+\tilde{Q}_{g_0,\,f_0}(\tau),\quad\tau\leq 0, \\
L_{f_0}h_0&=\Li_{B_{g_0}(h_0)^{\sharp}-X}g_0-\Li_{X}h_0-Q_{g_0}(h_0),\\
\Delta_{f_0}X-\frac{X}{2}+\nabla^{g_0}\div_{f_0}X&=\div_{f_0}\Li_{X}g_0.
\end{split}
\end{equation*}
Here, $\tilde{Q}_{g_0,\,f_0}(\tau)$ satisfies pointwise on $M\times(-\infty,0]$ for all $\gamma\in[0,1)$,
\begin{equation*}
\begin{split}
|\tilde{Q}_{g_0,\,f_0}(\tau)|_{g_0}&\leq C(n,g_0,\gamma)\varepsilon e^{\lambda_+\tau}\sum_{k=0}^2|\nabla^{g_0,\,k}h_0(\tau)|_{g_0}e^{-\gamma f_0}+C(n,g_0,\gamma)\varepsilon e^{\lambda_+\tau} |X|_{g_0}e^{-\gamma f_0},\\
|\nabla^{g_0}\tilde{Q}_{g_0,\,f_0}(\tau)|_{g_0}&\leq C(n,g_0,\gamma)\varepsilon e^{\lambda_+\tau}\sum_{k=0}^3|\nabla^{g_0,\,k}h_0(\tau)|_{g_0}e^{-\gamma f_0}\\
&\quad+C(n,g_0,\gamma)\varepsilon e^{\lambda_+\tau}\left(|X|_{g_0}+|\nabla^{g_0}X|_{g_0}\right)e^{-\gamma f_0}.
\end{split}
\end{equation*}
\end{lemma}
\begin{rk}
The third equation is exactly Lemma \ref{lem-bochner-lie-der} applied to $Y:=X$. 
\end{rk}
\begin{proof}
Recall that each family $(h_i(\tau))_{\tau\,\leq\, 0}$, $i=1,2$, satisfies 
\begin{equation}
\partial_{\tau}h_i(\tau)=L_{f_0}h_i(\tau)+Q_{g_0}(h_i(\tau)),\quad \tau\leq 0.
\end{equation}
In particular, if $L_{\varphi^*f_0}$ denotes the operator $\Delta_{\varphi^*g_0}+\nabla^{\varphi^*g_0}_{\varphi^*X_0}$,
\begin{equation*}
\begin{split}
\partial_{\tau}\varphi^*h_1(\tau)&=\varphi^*\partial_{\tau}h_1(\tau)\\
&=L_{\varphi^*f_0}\varphi^*h_1(\tau)+Q_{\varphi^*g_0}(\varphi^*h_1(\tau))\\
&=L_{f_0}\varphi^*h_1(\tau)+\left(L_{\varphi^*f_0}\varphi^*h_1(\tau)-L_{f_0}\varphi^*h_1(\tau)\right)+Q_{\varphi^*g_0}(\varphi^*h_1(\tau)).
\end{split}
\end{equation*}
Therefore,
\begin{equation}
\begin{split}\label{defn-Q-tilde}
\partial_{\tau}h_0(\tau)&=L_{f_0}h_0(\tau)+\tilde{Q}_{g_0,\,f_0}(\tau),\\
\tilde{Q}_{g_0,\,f_0}(\tau)&:=\left(Q_{g_0}(h_2(\tau))-Q_{g_0}(\varphi^*h_1(\tau))\right)\\
&\quad-\left(L_{\varphi^*f_0}-L_{f_0}\right)\varphi^*h_1(\tau)+\left(Q_{g_0}(\varphi^*h_1(\tau))-Q_{\varphi^*g_0}(\varphi^*h_1(\tau))\right).
\end{split}
\end{equation}
Inspecting each term in the definition of $Q_{g_0}(\cdot)$ recalled schematically in Lemma \ref{app-C-monster}, one gets for $\tau\leq 0$:
\begin{equation*}
\begin{split}
&|Q_{g_0}(h_2(\tau))-Q_{g_0}(\varphi^*h_1(\tau))|_{g_0}\leq C(n,g_0)|\nabla^{g_0,\,2}h_0(\tau)|\left(|h_2(\tau)|_{g_0}+|\varphi^*h_1(\tau)|_{g_0}\right)\\
&\quad+C(n,g_0)|\nabla^{g_0}h_0(\tau)|_{g_0}\left(|\nabla^{g_0}h_2(\tau)|_{g_0}+|\nabla^{g_0}\varphi^*h_1(\tau)|_{g_0}\right)\\
&\quad+C(n,g_0)|h_0(\tau)|_{g_0}\left(|\nabla^{g_0,\,2}h_2(\tau)|_{g_0}+|\nabla^{g_0,\,2}\varphi^*h_1(\tau)|_{g_0}+|\nabla^{g_0}h_2(\tau)|^2_{g_0}+|\nabla^{g_0}\varphi^*h_1(\tau)|^2_{g_0}\right)\\
&\quad+C(n,g_0)|h_0(\tau)|_{g_0}|\Rm(g_0)|_{g_0}\left(|h_2(\tau)|_{g_0}(1+|h_2(\tau)|_{g_0})+|\varphi^*h_1(\tau)|_{g_0}(1+|\varphi^*h_1(\tau)|_{g_0})\right)\\
&\leq C(n,g_0,\gamma)\varepsilon e^{\lambda_+\tau}\sum_{k=0}^2|\nabla^{g_0,\,k}h_0(\tau)|_{g_0}e^{-\gamma f_0},
\end{split}
\end{equation*}
where we have used Proposition \ref{prop-alm-sharp-rate} (applied to $h_2(\tau)$) and Lemma \ref{lemma-rough-equiv} in the first and second inequality.

The last term of \eqref{defn-Q-tilde} $\left(Q_{g_0}(\varphi^*h_1(\tau))-Q_{\varphi^*g_0}(\varphi^*h_1(\tau))\right)$ can be handled similarly for $\tau\leq 0$ by applying the identity $\nabla^{\varphi^*g_0}T=\nabla^{g_0}T+(\varphi^*g_0)^{-1}\ast \nabla^{g_0}h_0\ast T$ for an arbitrary tensor $T$ repeatedly.

As for the second term of \eqref{defn-Q-tilde}, one gets for $\gamma\in[0,1)$:
\begin{equation*}
\begin{split}
|\left(L_{\varphi^*f_0}-L_{f_0}\right)\varphi^*h_1(\tau)|_{g_0}&\leq |(\Delta_{\varphi^*g_0}-\Delta_{g_0})\varphi^*h_1(\tau)|_{g_0}+|(\nabla^{\varphi^*g_0}_{\varphi^*X_0}-\nabla^{g_0}_{X_0})\varphi^*h_1(\tau)|_{g_0}\\
&\quad +|\Rm(\varphi^*g_0)\ast_{\varphi^*g_0}\varphi^*h_1(\tau)-\Rm(g_0)\ast_{g_0}\varphi^*h_1(\tau)|_{g_0}\\
&\leq C(n,g_0)\sum_{k=0}^2|\nabla^{g_0,\,k}h_0(\tau)|_{g_0}|\nabla^{g_0,\,2-k}\varphi^*h_1(\tau)|_{g_0}\\
&\quad+|\nabla^{\varphi^*g_0}\varphi^*h_1(\tau)|_{g_0}|X|_{g_0}+C(n,g_0)|\nabla^{g_0}h_0(\tau)|_{g_0}|\varphi^*h_1(\tau)|_{g_0}|X_0|_{g_0}\\
&\leq C(n,g_0,\gamma)\varepsilon e^{\lambda_+\tau}\sum_{k=0}^2|\nabla^{g_0,\,k}h_0(\tau)|_{g_0}e^{-\gamma f_0}\\
&\quad+C(n,g_0,\gamma)\varepsilon e^{\lambda_+\tau} |X|_{g_0}e^{-\gamma f_0},
\end{split}
\end{equation*}
where we have used Lemma \ref{lemma-rough-equiv} in the second and last inequalities together with the polynomial growth of $|X_0|_{g_0}$ with respect to $f_0$.

Regarding the estimate on $|\nabla^{g_0}\tilde{Q}_{g_0,f_0}(\tau)|_{g_0}$, one first differentiates \eqref{defn-Q-tilde} and proceed similarly to the above reasonings that led to the pointwise estimate on $|\tilde{Q}_{g_0,f_0}(\tau)|_{g_0}$.

The second elliptic equation can be proved exactly as in \eqref{hello-eqn-non-linear} by linearizing the soliton equation satisfied by $(M,\varphi^*g_0,\varphi^*X_0)$ with respect to that satisfied by $(M,g_0,X_0:=\nabla^{g_0}f_0)$.

The third equation is exactly Lemma \ref{lem-bochner-lie-der} applied to $Y:=X$.
\end{proof}

\begin{lemma}\label{lemma-ene-cov-der-h_0}
There exists $C=C(n,g_0)>0$ such that for $\tau\leq 0$,
\begin{equation*}
\begin{split}
\sum_{k=0}^2\|\nabla^{g_0,\,k}h_0(\tau)\|_{L^2_{f_0}}&\leq C\left(\|h_0(\tau)\|_{L^2_{f_0}}+\|L_{f_0}h_0(\tau)\|_{L^2_{f_0}}\right),\\
\sum_{k=0}^3\|\nabla^{g_0,\,k}h_0(\tau)\|_{L^2_{f_0}}&\leq C\left(\|h_0(\tau)\|_{L^2_{f_0}}+\|L_{f_0}h_0(\tau)\|_{L^2_{f_0}}+\|\nabla^{g_0}L_{f_0}h_0(\tau)\|_{L^2_{f_0}}\right).
\end{split}
\end{equation*}
\end{lemma}
\begin{proof}[Proof of Lemma \ref{lemma-ene-cov-der-h_0}]
Observe that it suffices to prove Lemma \ref{lemma-ene-cov-der-h_0} for $k=1,\,2,\,3$.

For $k=1$, an integration by parts legitimated by Lemma \ref{lemma-rough-equiv} implies for $\tau\leq 0$:
\begin{equation*}
\begin{split}
\|\nabla^{g_0}h_0(\tau)\|_{L^2_{f_0}}^2&=\langle-\Delta_{f_0}h_0(\tau),h_0(\tau)\rangle_{L^2_{f_0}}\\
&=\langle-L_{f_0}h_0(\tau),h_0(\tau)\rangle_{L^2_{f_0}}+\langle 2\Rm(g_0)\ast h_0(\tau),h_0(\tau)\rangle_{L^2_{f_0}}\\
&\leq C\left(\|L_{f_0}h_0(\tau)\|_{L^2_{f_0}}^2+\|h_0(\tau)\|_{L^2_{f_0}}^2\right),
\end{split}
\end{equation*}
where $\Delta_{f_0}h_0(\tau):=\Delta_{g_0}h_0(\tau)+\nabla^{g_0}_{X_0}h_0(\tau)$. Here we have used Cauchy-Schwarz inequality together with the boundedness of the curvature of $g_0$ on $M$ in the last line.

As for $k\geq 2$, recall the following commutation formula that can be found for instance in \cite[Proposition $2.3$]{Der-Asy-Com-Egs}:
\begin{equation*}
[\Delta_{f_0},\nabla^{g_0,\,k}]h_0(\tau)=-\frac{k}{2}\nabla^{g_0,\,k}h_0(\tau)+\sum_{i=0}^k\nabla^{g_0,\,i}\Rm(g_0)\ast \nabla^{g_0,\,k-i}h_0(\tau).
\end{equation*}
This leads to the following Bochner identity:
\begin{equation*}
\begin{split}
\Delta_{f_0}|\nabla^{g_0}h_0(\tau)|^2_{g_0}&=2|\nabla^{g_0,\,2}h_0(\tau)|_{g_0}+2\langle[\Delta_{f_0},\nabla^{g_0}]h_0(\tau),\nabla^{g_0}h_0(\tau)\rangle_{g_0}+2\langle \nabla^{g_0}\Delta_{f_0}h_0(\tau),\nabla^{g_0}h_0(\tau)\rangle_{g_0}\\
&=2|\nabla^{g_0,\,2}h_0(\tau)|_{g_0}+2\left\langle-\frac{1}{2}\nabla^{g_0}h_0(\tau)+\sum_{i=0}^1\nabla^{g_0,\,i}\Rm(g_0)\ast \nabla^{g_0,\,1-i}h_0(\tau),\nabla^{g_0}h_0(\tau)\right\rangle_{g_0}\\
&\quad+2\langle \nabla^{g_0}L_{f_0}h_0(\tau),\nabla^{g_0}h_0(\tau)\rangle_{g_0}-2\langle \nabla^{g_0}(\Rm(g_0)\ast h_0(\tau)),\nabla^{g_0}h_0(\tau)\rangle_{g_0}.
\end{split}
\end{equation*}
As the (covariant derivatives of the) curvature of $g_0$ are bounded on $M$, one gets pointwise on $M$:
\begin{equation*}
\begin{split}
\Delta_{f_0}|\nabla^{g_0}h_0(\tau)|^2_{g_0}&\geq 2|\nabla^{g_0,\,2}h_0(\tau)|_{g_0}+2\langle \nabla^{g_0}L_{f_0}h_0(\tau),\nabla^{g_0}h_0(\tau)\rangle_{g_0}\\
&\quad-C|\nabla^{g_0}h_0(\tau)|^2_{g_0}-C|h_0(\tau)|_{g_0}|\nabla^{g_0}h_0(\tau)|_{g_0}.
\end{split}
\end{equation*}
Integrating the previous inequality and recalling the definition of the weighted Licherowicz operator lead to:
\begin{equation*}
\begin{split}
0&=\int_M\Delta_{f_0}|\nabla^{g_0}h_0(\tau)|^2_{g_0}\,d\mu_{f_0}\\
&\geq 2\|\nabla^{g_0,\,2}h_0(\tau)\|_{L^2_{f_0}}-2\langle L_{f_0}h_0(\tau),\Delta_{f_0}h_0(\tau)\rangle_{L^2_{f_0}}-C\left(\|h_0(\tau)\|^2_{L^2_{f_0}}+\|\nabla^{g_0}h_0(\tau)\|_{L^2_{f_0}}^2\right)\\
&\geq 2\|\nabla^{g_0,\,2}h_0(\tau)\|_{L^2_{f_0}}-C\| L_{f_0}h_0(\tau)\|^2_{L^2_{f_0}}-C\left(\|h_0(\tau)\|^2_{L^2_{f_0}}+\|\nabla^{g_0}h_0(\tau)\|_{L^2_{f_0}}^2\right),
\end{split}
\end{equation*}
where we have used Cauchy-Schwarz inequality in the second and third inequality and where $C$ is a positive constant that may vary from line to line. Here, the integration by parts are legitimated by Lemma \ref{lemma-rough-equiv}. The result follows once we invoke the previously established bound on $\|\nabla^{g_0}h_0(\tau)\|_{L^2_{f_0}}$.

The bound on $\|\nabla^{g_0,\,3}h_0(\tau)\|_{L^2_{f_0}}$ can be proved along the same lines.
\end{proof}

\begin{lemma}\label{lemma-ene-div-X-h_0-sec-7}
For each $\gamma\in (0,1)$, there exists $C=C(n,g_0,\gamma)>0$ such that,
\begin{equation*}
\left|\div_{f_0}\div_{f_0}\left(\Li_{X}g_0\right)\right|+|\div_{f_0}X|_{g_0}+|\nabla^{g_0}\div_{f_0}X|_{g_0}\leq C\,\varepsilon e^{-\gamma f_0}.
\end{equation*}
%Moreover, for each $\gamma\in (0,1)$, there exists $C_1>0$ such that $|\nabla^{g_0}\div_{f_0}X|_{g_0}\leq C_1 e^{-\gamma f_0}.$
\end{lemma}

\begin{proof}
The required bound on $\left|\div_{f_0}\div_{f_0}\left(\Li_{X}g_0\right)\right|$ follows from (repeated applications of) Lemma \ref{lemma-rough-equiv} together with the fact that $\nabla^{g_0}f_0$ grows at most polynomially in $f_0$ (since it grows like $f_0^{1/2}$).

We first claim that $\div_{f_0}X=o(f_0^{\frac{1}{2}})$, a fact that can be proved as in the proof of Lemma \ref{lemma-pt-est-wei-div}. Then, Lemma \ref{lem-equiv-gauge-eqn} applied to $Z:=(\div_{f_0}\Li_Xg_0)^{\sharp}$ ensures that $u:=\div_{f_0}X$ satisfies thanks to the previously established bound on $\left|\div_{f_0}\div_{f_0}\left(\Li_{X}g_0\right)\right|$:
\begin{equation*}
\Delta_{f_0}u-\frac{u}{2}\geq -C\varepsilon \,e^{-\gamma f},
\end{equation*}
where $\gamma\in(0,1)$ and $C=C(n,g_0,\gamma)>0$.
Together with the fact that $\div_{f_0}X=o(f_0^{\frac{1}{2}})$ as stated previously, we are then in a position to apply Lemma \ref{lemma-max-ppe-sub-lin} to ensure that $\div_{f_0}X\leq 3C\varepsilon\,e^{-\gamma f_0}$ which is of the desired form by substituting $3C$ by $C$.
A similar reasoning applied to $-\div_{f_0}X$ leads to the desired decay on $|\div_{f_0}X|$.

In order to derive the decay on the gradient of $\div_{f_0}X$, we invoke the following basic interpolation inequality whose proof can be found for instance in \cite{Aubin} for $m\geq 2$ and $x\in M$:
\begin{equation}
\begin{split}
&\sup_{B_{g_0}\left(x,\frac{\iota}{2} f_0(x)^{\frac{1}{2}}\right)}f_0(x)^{\frac{1}{2}}|\nabla^{g_0}\div_{f_0}X|_{g_0}\\
&\leq C\sup_{B_{g_0}(x,\iota f_0(x)^{\frac{1}{2}} )}|\div_{f_0}X|_{g_0}^{1-\frac{1}{m}}\cdot \sup_{B_{g_0}\left(x,\iota f_0(x)^{\frac{1}{2}}\right)}\bigg(\sum_{k=0}^mf_0(x)^{\frac{k}{2}}|\nabla^{g_0,\,k}\div_{f_0}X|_{g_0}\bigg)^{\frac{1}{m}},
\end{split}
\end{equation}
where $\inj_{g_0}(x)\geq \iota f_0(x)^{\frac{1}{2}}$ from some positive constant $\iota$ uniform in $x\in M$, a property that holds thanks to the asymptotics of $g_0$ provided by Theorem \ref{ac-str-exp}.

Now, by repeated application of Lemma \ref{lemma-rough-equiv}, covariant derivatives of $\div_{f_0}X$ are growing at most polynomially in $f_0$ in a quantitative way with respect to $\varepsilon$. This fact combined with the previous exponential decay on $\div_{f_0}X$ leads to the desired exponential decay on the gradient of $\div_{f_0}X$.
\end{proof}

\begin{lemma}\label{lemma-ene-X-h_0}
The vector field $X\in L^2_{f_0}$ and there exists $C=C(n,g_0)>0$ such that,
\begin{equation*}
\|X\|_{L^2_{f_0}}+\|\nabla^{g_0}X\|_{L^2_{f_0}}\leq C\left(\|h_0\|_{L^2_{f_0}}+\|L_{f_0}h_0\|_{L^2_{f_0}}\right).
\end{equation*}
%In particular, if $\varepsilon$ is small enough, 
%\begin{equation*}
%\|X\|_{L^2_{f_0}}+\|L_{f_0}h_0\|_{L^2_{f_0}}\leq C_0\|h_0\|_{L^2_{f_0}}.
%\end{equation*}
\end{lemma}
\begin{proof}[Proof of Lemma \ref{lemma-ene-X-h_0}]
Let us first show that $X\in L^2_{f_0}$. Applying Lemma \ref{lem-bochner-lie-der} to $Y:=X$ gives:
\begin{equation*}
\Delta_{f_0}|X|_{g_0}-\frac{1}{2}|X|_{g_0}\geq -|\nabla^{g_0}\div_{f_0}X|_{g_0}-|\div_{f_0}\Li_Xg_0|_{g_0}.
\end{equation*}
Now, invoking Lemma \ref{lemma-ene-div-X-h_0-sec-7} let us prove that for all $\gamma\in (0,1)$, there exists $C>0$ such that on $M$:
\begin{equation}\label{exp-dec-rhs-X}
|\nabla^{g_0}\div_{f_0}X|_{g_0}+|\div_{f_0}\Li_Xg_0|_{g_0}\leq Ce^{-\gamma f_0}.
\end{equation}

As $|X|_{g_0}$ is bounded on $M$ thanks to Lemma \ref{lemma-rough-equiv}, one can apply Lemma \ref{lemma-max-ppe-sub-lin} to $u:=|X|_{g_0}$, which once combined with the decay \eqref{exp-dec-rhs-X} for a given $\gamma\in(0,1)$ leads to a bound of the form $|X|_{g_0}\leq 3C e^{-\gamma f_0}$ on $M$. This proves in particular that $X\in L^2_{f_0}$ if we choose $\gamma\in(1/2,1)$.

Regarding the quantitative $L^2_{f_0}$-bounds on $X$, one can apply Lemma \ref{korn-type-id} to $Y:=X$ which leave us  with estimating $|\Li_Xg_0|_{g_0}$ pointwise from above. Recall by Lemma \ref{lem-hello-eqn-non-linear-bis} that pointwise on $M$:
\begin{equation*}
\begin{split}
|\Li_Xg_0|_{g_0}&\leq |L_{f_0}h_0|_{g_0}+C\left(|\nabla^{g_0,\,2}h_0|_{g_0}+|\nabla^{g_0}h_0|_{g_0}+|h_0|_{g_0}\right)+|\Li_Xh_0|_{g_0}\\
&\leq|L_{f_0}h_0|_{g_0}+C\left(|\nabla^{g_0,\,2}h_0|_{g_0}+(1+|X|_{g_0})|\nabla^{g_0}h_0|_{g_0}+(1+|\nabla^{g_0}X|_{g_0})|h_0|_{g_0}\right),
\end{split}
\end{equation*}
where we have used the identity $\Li_Xh_0=\nabla^{g_0}_Xh_0+h_0\ast \nabla^{g_0}X$ in the last line together with Lemma \ref{app-C-monster} and Lemma \ref{lemma-rough-equiv} (in order to bound the non-linearities).
Squaring and integrating the previous pointwise estimate gives:
\begin{equation}
\begin{split}\label{rough-L2-lie-der}
\|\Li_Xg_0\|_{L^2_{f_0}}&\leq\|L_{f_0}h_0\|_{L^2_{f_0}}+C\left(\sum_{k=0}^2\|\nabla^{g_0,\,k}h_0\|_{L^2_{f_0}}+\sum_{k=0}^1\|\nabla^{g_0,\,1-k}h_0\|_{C^0}\|\nabla^{g_0,\,k}X\|_{L^2_{f_0}}\right)\\
&\leq C\left(\|h_0\|_{L^2_{f_0}}+\|L_{f_0}h_0\|_{L^2_{f_0}}\right)+C\varepsilon\left(\|X\|_{L^2_{f_0}}+\|\nabla^{g_0}X\|_{L^2_{f_0}}\right),
\end{split}
\end{equation}
where $C=C(n,g_0)>0$ and where we have invoked Lemma \ref{lemma-rough-equiv} and Lemma \ref{lemma-ene-cov-der-h_0} in the last inequality. By choosing $0<\varepsilon\leq \varepsilon(n,g_0)$ small enough, an absorption principle obtained by feeding back the above integral estimate \eqref{rough-L2-lie-der} in [\eqref{korn-type-id-est}, Lemma \ref{korn-type-id}] lead to the required estimate on $\|X\|_{L^2_{f_0}}$.
\end{proof}

The following two lemmata are abstract statements regarding evolution equations defined on Hilbert spaces adapted to our setting: 
\begin{lemma}\label{lemma-hilb-est}
There exists $C=C(n,g_0)>0$ such that for $\tau\leq 0$,
\begin{equation}
\begin{split}\label{est-lemma-hilb-est}
&\|\pi_{-}h_0(\tau)\|_{L^2_{f_0}}^2+\|\nabla^{g_0}\pi_{-}h_0(\tau)\|^2_{L^2_{f_0}}\leq C\int_{-\infty}^{\tau}\|\tilde{Q}_{g_0,f_0}(s)\|^2_{L^2_{f_0}}\,ds,\\
&\int_{-\infty}^{\tau}\|\pi_{-}h_0(s)\|^2_{L^2_{f_0}}+\|\nabla^{g_0}\pi_{-}h_0(s)\|^2_{L^2_{f_0}}+\|L_{f_0}\pi_{-}h_0(s)\|^2_{L^2_{f_0}}\,ds\leq C\int_{-\infty}^{\tau}\|\tilde{Q}_{g_0,f_0}(s)\|^2_{L^2_{f_0}}\,ds,\\
&\|L_{f_0}\pi_{-}h_0(\tau)\|_{L^2_{f_0}}^2+\int_{-\infty}^{\tau}\|\nabla^{g_0}L_{f_0}\pi_{-}h_0(s)\|^2_{L^2_{f_0}}\,ds\leq C\int_{-\infty}^{\tau}\|\tilde{Q}_{g_0,f_0}(s)\|^2_{L^2_{f_0}}+\|\nabla^{g_0}\tilde{Q}_{g_0,f_0}(s)\|^2_{L^2_{f_0}}\,ds.
\end{split}
\end{equation}
\end{lemma}

\begin{proof}
Let us handle the component $\pi_{-}h_0(\tau)=:H_0^{-}(\tau)$ as follows: it satisfies for $\tau\leq 0$,
\begin{equation}
\begin{split}\label{bdy-cond-H-pi-minus}
(\partial_{\tau}-L_{f_0})H^{-}_0(\tau)&=\pi_{-}(\tilde{Q}_{g_0,f_0}(\tau)), \quad \tau\leq 0,\\
\|H_0^{-}(\tau)\|_{L^2_{f_0}}+\|\nabla^{g_0}H_0^{-}(\tau)\|_{L^2_{f_0}}&\leq C(n,g_0)\,\varepsilon e^{\lambda_+\tau},\quad \tau\leq 0,
\end{split}
\end{equation}
 thanks to Lemma \ref{lemma-rough-equiv}. 
 A direct computation shows that for $\tau\leq 0$:
 \begin{equation}
 \begin{split}\label{hilb-0-pi-minus}
\frac{1}{2}\partial_{\tau}\|H_0^{-}(\tau)\|_{L^2_{f_0}}^2&=\langle L_{f_0}H_0^{-}(\tau)+\tilde{Q}_{g_0,f_0}(\tau),H_0^{-}(\tau)\rangle_{L^2_{f_0}}\\
&\leq -\lambda_{-}\|H_0^{-}(\tau)\|_{L^2_{f_0}}^2+\|\tilde{Q}_{g_0,f_0}(\tau)\|_{L^2_{f_0}}\|H_0^{-}(\tau)\|_{L^2_{f_0}}\\
&\leq -\left(\lambda_{-}-\delta\right)\|H_0^{-}(\tau)\|_{L^2_{f_0}}^2+C_{\delta}\|\tilde{Q}_{g_0,f_0}(\tau)\|_{L^2_{f_0}}^2,
\end{split}
\end{equation}
where we have used Young's inequality and Cauchy-Schwarz inequality in the last line. Here $-\lambda_{-}>0$ is the biggest negative eigenvalue of $L_{f_0}$ and $\delta\in(0,\lambda_{-})$.

By integrating in time and invoking \eqref{bdy-cond-H-pi-minus}:
 \begin{equation}
 \begin{split}\label{hilb-1-pi-minus}
\|H_0^{-}(\tau)\|_{L^2_{f_0}}^2+\left(\lambda_{-}-\delta\right)\int_{-\infty}^{\tau}\|H_0^{-}(s)\|_{L^2_{f_0}}^2\,ds\leq C_{\delta}\int_{-\infty}^{\tau}\|\tilde{Q}_{g_0,f_0}(s)\|_{L^2_{f_0}}^2\,ds.
\end{split}
\end{equation}
 Going back to the first line of \eqref{hilb-0-pi-minus} and integrating by parts give:
   \begin{equation}
 \begin{split}\label{hilb-1-bis-pi-minus}
\|H_0^{-}(\tau)\|_{L^2_{f_0}}^2+\int_{-\infty}^{\tau}\|\nabla^{g_0}H_0^{-}(s)\|_{L^2_{f_0}}^2\,ds&\leq C\int_{-\infty}^{\tau}\|H_0^{-}(s)\|_{L^2_{f_0}}^2+\|\tilde{Q}_{g_0,f_0}(s)\|_{L^2_{f_0}}^2\,ds\\
&\leq C'\int_{-\infty}^{\tau}\|\tilde{Q}_{g_0,f_0}(s)\|_{L^2_{f_0}}^2\,ds,
\end{split}
\end{equation}
where we have used \eqref{hilb-1-pi-minus} in the second inequality and where $C'=C'(n,g_0)>0$.

 Similarly, a direct computation and an integration by parts  show then:
 \begin{equation}
 \begin{split}\label{hilb-2-pi-minus}
\frac{1}{2}\partial_{\tau}\|\nabla^{g_0}H_0^{-}(\tau)\|_{L^2_{f_0}}^2&=\langle \nabla^{g_0}\left(L_{f_0}H_0^{-}(\tau)+\tilde{Q}_{g_0,f_0}(\tau)\right),\nabla^{g_0}H_0^{-}(\tau)\rangle_{L^2_{f_0}}\\
&=-\langle L_{f_0}H_0^{-}(\tau)+\tilde{Q}_{g_0,f_0}(\tau),\Delta_{f_0}H_0^{-}(\tau)\rangle_{L^2_{f_0}}\\
&=-\|L_{f_0}H_0^{-}(\tau)\|_{L^2_{f_0}}^2+\langle L_{f_0}H_0^{-}(\tau),2\Rm(g_0)\ast H_0^{-}(\tau)-\tilde{Q}_{g_0,f_0}(\tau)\rangle_{L^2_{f_0}}\\
&\quad +\langle \tilde{Q}_{g_0,f_0}(\tau),2\Rm(g_0)\ast H_0^{-}(\tau)\rangle_{L^2_{f_0}}\\
&\leq -\frac{1}{2}\|L_{f_0}H_0^{-}(\tau)\|_{L^2_{f_0}}^2+C(n,g_0)\left(\|H_0^{-}(\tau)\|_{L^2_{f_0}}^2+\|\tilde{Q}_{g_0,f_0}(\tau)\|^2_{L^2_{f_0}}\right),
\end{split}
\end{equation}
where we have used Young's inequality and Cauchy-Schwarz inequality in the last line.

Integrating in time and invoking \eqref{bdy-cond-H-pi-minus} leads to:
\begin{equation}
\begin{split}
\|\nabla^{g_0}H_0^{-}(\tau)\|^2_{L^2_{f_0}}+\int_{-\infty}^{\tau}\|L_{f_0}H_0^{-}(s)\|^2_{L^2_{f_0}}\,ds&\leq C\int_{-\infty}^{\tau}\|H_0^{-}(s)\|_{L^2_{f_0}}^2+\|\tilde{Q}_{g_0,f_0}(s)\|^2_{L^2_{f_0}}\,ds\\
&\leq C'\int_{-\infty}^{\tau}\|\tilde{Q}_{g_0,f_0}(s)\|^2_{L^2_{f_0}}\,ds
\end{split}
\end{equation}
where we have used \eqref{hilb-1-pi-minus} in the last line again and where $C'=C'(n,g_0)>0$.

A similar reasoning leads to the third required estimate by differentiating $\|L_{f_0}H_0^{-}(\tau)\|^2_{L^2_{f_0}}$.
\end{proof}

\begin{lemma}\label{lemma-hilb-est-II}
For $\delta\in (0,\lambda_+)$, there exists $C=C(n,g_0,\delta)>0$ such that for $\tau\leq 0$,
\begin{equation}
\begin{split}\label{est-lemma-hilb-est-II}
e^{-2(\lambda_+-\delta)\tau}\|\pi_{+}h_0(\tau)-P_{\tau}(\pi_+h_0)\|^2_{L^2_{f_0}}&\leq C\int_{\tau}^{0}e^{-2(\lambda_+-\delta)s}\|\tilde{Q}_{g_0,f_0}(s)\|^2_{L^2_{f_0}}\,ds,\\
\int_{-\infty}^{0}\|\pi_{+}h_0(s)-P_{s}(\pi_+h_0)\|^2_{L^2_{f_0}}\,ds&\leq C\int_{-\infty}^0\|\tilde{Q}_{g_0,f_0}(s)\|^2_{L^2_{f_0}}\,ds,
\end{split}
\end{equation}
where $P_\tau$ is the heat semigroup generated $L_{f_0}$ see Section \ref{sec-spec-l-f}. Furthermore, there exists $C=C(n,g_0)>0$ such that for $\tau\leq 0$,
\begin{equation}
\begin{split}\label{est-lemma-hilb-est-III}
\|\pi_{0}h_0(\tau)\|_{L^2_{f_0}}&\leq C\int_{-\infty}^{\tau}\|\tilde{Q}_{g_0,f_0}(s)\|_{L^2_{f_0}}\,ds.
\end{split}
\end{equation}
\end{lemma}

\begin{proof}
The proof of this lemma is similar to that of Lemma \ref{lemma-hilb-est} and will therefore be omitted. 
\end{proof}

\begin{lemma}\label{lemma-control-l2-pi+-h_0}
There exists $C=C(n,g_0)>0$ and $\varepsilon(n,g_0)>0$ such that for $0<\varepsilon\leq \varepsilon(n,g_0)$,
\begin{equation*}
\begin{split}
&\|h_0-\pi_+h_0\|^2_{L^2_{f_0}}+\|\nabla^{g_0}\pi_{-}h_0\|^2_{L^2_{f_0}}+\|L_{f_0}\pi_{-}h_0\|^2_{L^2_{f_0}}\leq C\varepsilon^2\left(\|\pi_+h_0\|^2_{L^2_{f_0}}+\|X\|^2_{L^2_{f_0}}+\|\nabla^{g_0}X\|^2_{L^2_{f_0}}\right).
\end{split}
\end{equation*}
\end{lemma}

\begin{proof}
According to Lemma \ref{lem-hello-eqn-non-linear-bis} applied with $\gamma=0$,
\begin{equation}
\begin{split}\label{feed-me-back}
\|\tilde{Q}_{g_0,\,f_0}(\tau)\|_{L^2_{f_0}}&\leq C\varepsilon e^{\lambda_+\tau}\sum_{k=0}^2\|\nabla^{g_0,\,k}h_0(\tau)\|_{L^2_{f_0}}+C\varepsilon e^{\lambda_+\tau} \|X\|_{L^2_{f_0}},\quad\tau\leq 0,
\end{split}
\end{equation}
which gives after squaring and integration in time for $\tau\leq 0$:
\begin{equation*}
\begin{split}
\int_{-\infty}^{0}\|\tilde{Q}_{g_0,\,f_0}(s)\|^2_{L^2_{f_0}}\,ds&\leq C\varepsilon^2 \int_{-\infty}^{0}\sum_{k=0}^2\|\nabla^{g_0,\,k}h_0(s)\|^2_{L^2_{f_0}}\,ds+C\varepsilon^2  \|X\|^2_{L^2_{f_0}}\\
&\leq C\varepsilon^2 \int_{-\infty}^{0}\|h_0(s)\|^2_{L^2_{f_0}}+\|L_{f_0}h_0(s)\|^2_{L^2_{f_0}}\,ds +C\varepsilon^2 \|X\|^2_{L^2_{f_0}},
\end{split}
\end{equation*}
where we have used Lemma \ref{lemma-ene-cov-der-h_0} in the last line. By choosing $\varepsilon\leq \varepsilon (n,g_0)$ small enough, the previous estimate combined with Lemma \ref{lemma-hilb-est} leads to the following intermediate estimate:
\begin{equation}
\begin{split}\label{int-est-Q}
\int_{-\infty}^{0}\|\tilde{Q}_{g_0,\,f_0}(s)\|^2_{L^2_{f_0}}\,ds&\leq C\varepsilon^2 \int_{-\infty}^{0}\|h_0(s)\|^2_{L^2_{f_0}}\,ds +C\varepsilon^2  \|X\|^2_{L^2_{f_0}}\\
&\leq C\varepsilon^2 \int_{-\infty}^{0}\|h_0(s)-P_s(\pi_+h_0)\|^2_{L^2_{f_0}}\,ds \\
&\quad+C\varepsilon^2 \int_{-\infty}^{0}\|P_s(\pi_+h_0)\|^2_{L^2_{f_0}}\,ds+C\varepsilon^2  \|X\|^2_{L^2_{f_0}}.
\end{split}
\end{equation}
Now, by Lemmata \ref{lemma-hilb-est} and \ref{lemma-hilb-est-II}:
\begin{equation*}
\begin{split}
\int_{-\infty}^{0}\|h_0(s)-P_{s}(\pi_+h_0)\|_{L^2_{f_0}}^2\,ds&\leq C \int_{-\infty}^{0}\|\tilde{Q}_{g_0,\,f_0}(s)\|^2_{L^2_{f_0}}\,ds+\int_{-\infty}^{0}\|\pi_0h_0(s)\|^2_{L^2_{f_0}}\,ds,
\end{split}
\end{equation*}
which gives thanks to \eqref{int-est-Q} and  after choosing $0<\varepsilon\leq \varepsilon(n,g_0)$ small enough, 
\begin{equation}
\begin{split}\label{first-est-int-l2-h}
\int_{-\infty}^{0}\|\tilde{Q}_{g_0,\,f_0}(s)\|^2_{L^2_{f_0}}\,ds&\leq C\varepsilon^2 \int_{-\infty}^{0}\|\pi_0h_0(s)\|^2_{L^2_{f_0}}+\|P_s(\pi_+h_0)\|^2_{L^2_{f_0}}\,ds +C\varepsilon^2  \|X\|^2_{L^2_{f_0}},\\
\int_{-\infty}^{0}\|h_0(s)-P_{s}(\pi_+h_0)\|_{L^2_{f_0}}^2\,ds&\leq C\int_{-\infty}^{0}\|\pi_0h_0(s)\|^2_{L^2_{f_0}}\,ds\\
&\quad+C\varepsilon^2\int_{-\infty}^{0}\|P_s(\pi_+h_0)\|^2_{L^2_{f_0}}\,ds +C\varepsilon^2  \|X\|^2_{L^2_{f_0}}.
\end{split}
\end{equation}
Now, by \eqref{est-lemma-hilb-est-III}, Lemma \ref{lemma-hilb-est-II} and \eqref{feed-me-back},
\begin{equation*}
\begin{split}
\int_{-\infty}^{0}\|\pi_0h_0(s)\|^2_{L^2_{f_0}}\,ds&\leq C\int_{-\infty}^{0}\left(\int_{-\infty}^{\tau}\|\tilde{Q}_{g_0,\,f_0}(s)\|_{L^2_{f_0}}\,ds\right)^2d\tau\\
&\leq C'\varepsilon^2\left(\int_{-\infty}^0\left(\int_{-\infty}^{\tau}e^{\lambda_+s}\sum_{k=0}^2\|\nabla^{g_0,\,k}h_0(s)\|_{L^2_{f_0}}\,ds\right)^2d\tau+\|X\|^2_{L^2_{f_0}}\right)\\
&\leq C'\varepsilon^2\left(\int_{-\infty}^0e^{2\lambda_+\tau}\left(\sum_{k=0}^2\int_{-\infty}^{\tau}\|\nabla^{g_0,\,k}h_0(s)\|^2_{L^2_{f_0}}\,ds\right)d\tau+\|X\|^2_{L^2_{f_0}}\right)\\
&\leq C'\varepsilon^2\left(\sum_{k=0}^2\int_{-\infty}^{0}\|\nabla^{g_0,\,k}h_0(s)\|^2_{L^2_{f_0}}\,ds+\|X\|^2_{L^2_{f_0}}\right),
\end{split}
\end{equation*}
where we have used Cauchy-Schwarz inequality in the penultimate inequality. Here $C'=C'(n,g_0)>0$. Invoking Lemma \ref{lemma-ene-cov-der-h_0} gives then:
\begin{equation}
\begin{split}\label{int-pi0}
\int_{-\infty}^{0}\|\pi_0h_0(s)\|^2_{L^2_{f_0}}\,ds
&\leq C'\varepsilon^2\left(\int_{-\infty}^{0}\|h_0(s)\|^2_{L^2_{f_0}}+\|L_{f_0}h_0(s)\|^2_{L^2_{f_0}}\,ds+\|X\|^2_{L^2_{f_0}}\right)\\
&\leq C'\varepsilon^2 \int_{-\infty}^{0}\|\pi_0h_0(s)\|^2_{L^2_{f_0}}\,ds +C'\varepsilon^2 \int_{-\infty}^{0}\|P_s(\pi_+h_0)\|^2_{L^2_{f_0}}\,ds
+C'\varepsilon^2  \|X\|^2_{L^2_{f_0}},
\end{split}
\end{equation}
where we have used Lemma \ref{lemma-hilb-est} together with \eqref{first-est-int-l2-h} in the last line.

Feeding back \eqref{first-est-int-l2-h} with \eqref{int-pi0} and choosing $\varepsilon$ small enough gives the required two first estimates on $\|\pi_{-}h_0\|_{L^2_{f_0}}$ and $\|\nabla^{g_0}\pi_{-}h_0\|_{L^2_{f_0}}$ by evaluating Lemma \ref{lemma-hilb-est} at $\tau=0$, since $$\int_{-\infty}^{\tau}\|P_{s}(\pi_+h_0)\|^2_{L^2_{f_0}}\,ds\leq Ce^{2\lambda_+\tau}\|\pi_+h_0\|^2_{L^2_{f_0}}.$$

The required estimate on $\|\pi_0h_0\|_{L^2_{f_0}}$ can be obtained along the same lines by invoking \eqref{est-lemma-hilb-est-III}, Lemma \ref{lemma-hilb-est-II} at $\tau=0$ together with \eqref{feed-me-back}.

The last estimate on $\|L_{f_0}\pi_{-}h_0\|_{L^2_{f_0}}$ follows the same lines of reasoning by relying on the third integral estimate of Lemma \ref{lemma-hilb-est}.
\end{proof}

\begin{lemma}\label{lemma-pi-plus-X}
There exists $C=C(n,g_0)>0$ and $\varepsilon(n,g_0)>0$ such that for $0<\varepsilon\leq \varepsilon(n,g_0)$,
\begin{equation*}
\begin{split}
\|\pi_+h_0\|_{L^2_{f_0}}^2&\leq C\varepsilon\left(\|X\|_{L^2_{f_0}}^2+\|\nabla^{g_0}X\|_{L^2_{f_0}}^2+\|h_0-\pi_+h_0\|_{L^2_{f_0}}^2+\|\nabla^{g_0}\pi_{-}h_0\|^2_{L^2_{f_0}}\right).
\end{split}
\end{equation*}
\end{lemma}

\begin{proof}
Thanks to Lemma \ref{lem-hello-eqn-non-linear-bis},
\begin{equation}
\begin{split}\label{I-II-III-sec-7}
\lambda_+\|\pi_+h_0\|_{L^2_f}^2&\leq \langle L_{f_0}\pi_+h_0,\pi_+h_0\rangle_{L^2_{f_0}}=\langle L_{f_0}h_0,\pi_+h_0\rangle_{L^2_{f_0}}\\
&=\langle \Li_{B_{g_0}(h_0)^{\sharp}-X}g_0-\Li_{X}h_0-Q_{g_0}(h_0),\pi_+h_0\rangle_{L^2_{f_0}}=:I-II-III.
\end{split}
\end{equation}

Observe that $I$, $II$ and $III$ are well-defined in the Lebesgue sense thanks to Lemma \ref{lemma-rough-equiv} and Lemma \ref{lemma-ene-X-h_0}.

As for the proof of \eqref{est-I}, $I=0$ since $\div_{f_0}\pi_+h_0=0$ according to Lemma \ref{lemma-div-free}.

Regarding the integral $II$, observe that pointwise on $M$:
\begin{equation}\label{pt-id-lie-der-bis}
\begin{split}
|\langle\Li_{X}h_0,\pi_+h_0\rangle_{L^2_{f_0}}|&\leq C\int_M\left(|\nabla^{g_0}h_0|_{g_0}|X|_{g_0}+|\nabla^{g_0}X|_{g_0}|h_0|_{g_0}\right)|\pi_+h_0|_{g_0}\,d\mu_{f_0}\\
&\leq C\varepsilon\left(\|X\|^2_{L^2_{f_0}}+\|\nabla^{g_0}X\|_{L^2_{f_0}}^2+\|\pi_+h_0\|^2_{L^2_{f_0}}\right),
\end{split}
\end{equation}
where we have used Lemma \ref{lemma-rough-equiv} in the last line together with Young's inequality.

Finally, the integral $III$ can be handled thanks to the definition of $Q_{g_0}(h)$ given in Lemma \ref{app-C-monster} as follows:
\begin{equation}
\begin{split}\label{III-sec-7}
|III|&\leq C\int_M\left(|\nabla^{g_0,\,2}h_0|_{g_0}|h_0|_{g_0}+|\nabla^{g_0}h_0|^2_{g_0}\right)|\pi_+h_0|_{g_0}\,d\mu_{f_0}\\
&\leq C\varepsilon \left(\|h_0\|_{L^2_{f_0}}^2+\|\nabla^{g_0}h_0\|^2_{L^2_{f_0}}+\|\pi_+h_0\|_{L^2_{f_0}}^2\right)\\
&\leq C\varepsilon \left(\|h_0-\pi_+h_0\|_{L^2_{f_0}}^2+\|\nabla^{g_0}\pi_{-}h_0\|^2_{L^2_{f_0}}+\|\pi_+h_0\|_{L^2_{f_0}}^2\right),
\end{split}
\end{equation}
where we have used Lemma \ref{lemma-rough-equiv} in the second line together with Young's inequality together with Lemma \ref{lemma-eq-norms} in the last line.
The combination of $I=0$, \eqref{pt-id-lie-der-bis} and \eqref{III-sec-7} leads to the required estimate once $\varepsilon\leq \varepsilon(n,g_0)$ is small enough so that the norm $\|\pi_+h_0\|_{L^2_{f_0}}^2$ can be absorbed by the lefthand side of \eqref{I-II-III-sec-7}.
\end{proof}
We are in a good position to give a proof of Theorem \ref{theo-main}, item \eqref{item-mod-spa}:
\begin{proof}[Proof of Theorem \ref{theo-main}, item \eqref{item-mod-spa}]

By Lemma \ref{lemma-ene-X-h_0}, there exists $C(n,g_0)>0$ such that:
\begin{equation*}
\begin{split}
\|X\|^2_{L^2_{f_0}}+\|\nabla^{g_0}X\|^2_{L^2_{f_0}}&\leq C\left(\|h_0\|^2_{L^2_{f_0}}+\|L_{f_0}h_0\|^2_{L^2_{f_0}}\right)\\
&\leq C\left(\|h_0-\pi_+h_0\|^2_{L^2_{f_0}}+\|\pi_+h_0\|^2_{L^2_{f_0}}+\|L_{f_0}\pi_{-}h_0\|^2_{L^2_{f_0}}\right)\\
&\leq C\varepsilon^2\left(\|X\|^2_{L^2_{f_0}}+\|\nabla^{g_0}X\|^2_{L^2_{f_0}}\right)+C\|\pi_+h_0\|^2_{L^2_{f_0}},
\end{split}
\end{equation*}
where we have used Lemma \ref{lemma-control-l2-pi+-h_0} in the last line.
Therefore, if $0<\varepsilon\leq \varepsilon(n,g_0)$ is small enough, the previous estimate ensures that for some positive constant $C=C(n,g_0)$:
\begin{equation}
\begin{split}\label{1st-est-sec-7}
\|X\|^2_{L^2_{f_0}}+\|\nabla^{g_0}X\|^2_{L^2_{f_0}}&\leq C\|\pi_+h_0\|^2_{L^2_{f_0}}.
\end{split}
\end{equation}
Feeding back Lemma \ref{lemma-control-l2-pi+-h_0} with \eqref{1st-est-sec-7} leads to:
\begin{equation}
\begin{split}\label{2nd-est-sec-7}
&\|h_0-\pi_+h_0\|^2_{L^2_{f_0}}+\|\nabla^{g_0}\pi_{-}h_0\|^2_{L^2_{f_0}}\leq C\varepsilon^2\|\pi_+h_0\|^2_{L^2_{f_0}}.
\end{split}
\end{equation}

Finally, Lemma \ref{lemma-pi-plus-X} ensures:
\begin{equation}
\begin{split}\label{3rd-est-sec-7}
\|\pi_+h_0\|_{L^2_{f_0}}^2&\leq C\varepsilon\left(\|X\|_{L^2_{f_0}}^2+\|\nabla^{g_0}X\|_{L^2_{f_0}}^2+\|h_0-\pi_+h_0\|_{L^2_{f_0}}^2+\|\nabla^{g_0}\pi_{-}h_0\|^2_{L^2_{f_0}}\right)\\
&\leq C'\varepsilon \|\pi_+h_0\|^2_{L^2_{f_0}},
\end{split}
\end{equation}
where we have used \eqref{1st-est-sec-7} and \eqref{2nd-est-sec-7} in the last line. Here $C'=C'(n,g_0)>0$. Estimate \eqref{3rd-est-sec-7} implies that $\pi_+h_0=0$ provided $0<\varepsilon\leq \varepsilon(n,g_0)$ is chosen small enough. Therefore, estimates \eqref{1st-est-sec-7} and \eqref{2nd-est-sec-7} allow us to infer that $X=0$ and $h_0=0$, i.e. $\varphi$ is an isometry of $g_0$ and therefore $h_2=\varphi^*h_1$. 

As a consequence, $\varphi^*f_0-f_0$ is constant on $M$. By the soliton identities of Lemma \ref{id-sol-egs} and the normalization \eqref{norm-f-pot-fct} satisfied both by $f_0$ and $\varphi^*f_0$, we conclude that $\varphi^*f_0=f_0$ as desired.
\end{proof}

\section{Existence of unstable asymptotically conical expanding gradient Ricci solitons}\label{sec:unstable-examples}

%\todo[inline]{Decide where to we want to put this subsection, and adjust above discussion}

Let us now justify that Theorem \ref{theo-main} is not vacuous: unstable asymptotically conical expanding gradient Ricci solitons exist. We indicate here that such examples arise in a suitable limit, within the $\mathrm{SO}(3)\times\mathrm{SO}(2)$-invariant family of cohomogeneity-one expanders on $\mathbb S^2\times\mathbb{R}^2$ constructed in \cite{Buz-Dan-Gal-Wan,Nie-Win} and revisited by Rajan \cite{Rajan}. Following \cite{Rajan}, these solitons form a two-parameter family
\begin{equation*}
g_{q,\varepsilon}=ds^2+a(s)^2\,d\theta^2+b(s)^2\,g_{\mathbb S^2},\qquad q:=b(0)>0,\quad \varepsilon:=f''(0)>0,
\end{equation*}
with  $a(0)=0$, $a'(0)=1$, $b'(0)=0$, $f(0)=f'(0)=0$, and asymptotic to a cone over $\mathbb S^2\times \mathbb S^1$. In this parametrization, according to \cite[Lemma 4.3]{Rajan}, $\varepsilon=0$ corresponds to $f\equiv 0$ and by uniqueness, \cite[Appendix B]{Rajan}, the limiting family with $\varepsilon=0$ is the well-known Euclidean Schwarzschild-anti-de Sitter Einstein metric $\bar g_q$ on $\mathbb S^2\times\mathbb{R}^2$. By analytic dependence on initial data near the bolt proved in \cite[Appendix B]{Rajan} we have the following limit in the $C^{\infty}_{loc}$-topology,
\begin{equation}\label{limit-eps-zero}
\lim_{\varepsilon\,\to\, 0}(g_{q,\varepsilon},f_{q,\varepsilon})= (\bar g_q,0).
\end{equation}
In the family of Euclidean Schwarzschild-anti-de Sitter metrics, as $q\to 0$, $q^{-2}\bar g_q$ converges smoothly on compact subsets to the Euclidean Schwarzschild metric $g_{\operatorname{Sch}}$. 
%In conclusion, for $\varepsilon$ and $q$ chosen small enough, 

On $(\mathbb{R}^2\times \mathbb S^2, g_{\operatorname{Sch}})$, by \cite{BiquardOzuch,Ozuch-instab-S1,Kim-Ozuch}, there exists an explicit smooth transverse-traceless symmetric $2$-tensor $h$ with $|\nabla^{g_{\operatorname{Sch}},\,k}h|_{g_{\operatorname{Sch}}}=O(r^{-2-k})$ for all $k\geq 0$, satisfying
\begin{equation}\label{ineq-BO}
\mathcal{Q}_{\operatorname{Sch}}(h):=\int_{\mathbb{R}^2\times \mathbb S^2}\left(|\nabla^{g_{\operatorname{Sch}}}h|_{g_{\operatorname{Sch}}}^2-2\langle\Rm(g_{\operatorname{Sch}})\ast h,h\rangle_{g_{\operatorname{Sch}}}\right)\,d\mu_{g_{\operatorname{Sch}}}<0.
\end{equation}
Because $h$ decays quadratically and $g_{\operatorname{Sch}}$ has quadratic curvature decay, a smooth radial cut-off $\chi_A$ supported in $\{r\leq 2A\}$ with $\chi_A\equiv 1$ on $\{r\leq A\}$ and $|\nabla \chi_A|\leq CA^{-1}$ satisfies
$\mathcal{Q}_{\operatorname{Sch}}(\chi_A h)=\mathcal{Q}_{\operatorname{Sch}}(h)+O(A^{-3})$, hence, for $A$ large enough, we have $\mathcal{Q}(\chi_A h)<0$. Denote somewhat abusively $k:=\chi_A h$ the resulting tensor on $(\mathbb S^2\times \mathbb R^2,g_{q,\varepsilon})$. 

By \eqref{limit-eps-zero}, for $\varepsilon$ and $q$ small enough, this yields, 
\begin{equation*}
\mathcal{Q}_{q,\varepsilon}(k):=\int_{\mathbb{R}^2\times \mathbb S^2}\left(|\nabla^{g_{q,\varepsilon}}k|_{g_{q,\varepsilon}}^2-2\langle\Rm(g_{q,\varepsilon})\ast k,k\rangle_{g_{q,\varepsilon}}\right)\,e^{f_{q,\varepsilon}}d\mu_{g_{q,\varepsilon}}<0.
\end{equation*}
By the min-max principle, this means that there exists an eigentensor $h_0\in E_+\backslash\{0\}$. Moreover, by Lemma \ref{lem-equiv-gauge-eqn}, $h_0$ must be in $\ker\operatorname{div}_f$. 

In particular, Theorem \ref{theo-main} applies to the above family we just described and produces non self-similar Ricci flows coming out of the corresponding $\mathrm{SO}(3)\times\mathrm{SO}(2)$-invariant cones over $\mathbb S^2\times \mathbb S^1$. 

\newpage
\appendix

\section{Killing vector fields on a metric cone}\label{sec-kill-metric-cone}
%\todo[inline]{Killing vectors are orthogonal. I followed the notations of K-S (I couldn't find another reference), but we can obviously change them.}

\begin{prop}[Killing fields on a nonflat metric cone are horizontal]\label{claim-cone-killing-orth}
Let
$$
(C_0,g_c)=\big((0,\infty)\times\Sigma,\;dr^2+r^2g_{\Sigma}\big)
$$
be a non flat metric cone over a closed Riemannian manifold $(\Sigma^{n},g_{\Sigma})$.
If $\tilde{K}$ is a Killing vector field on $(C_0,g_c)$ and
$\omega:=g_c(\tilde{K},\cdot)$ its dual $1$-form, then
$$
g_{c}(\tilde{K},r\partial_{r})=\omega(r\partial_{r})= 0\qquad\text{ on }C_{0}.
$$
In other words, $\tilde{K}$ is tangent to every link $\{r=\mathrm{const}\}$
and is precisely the extension of a Killing field on
$(\Sigma,g_{\Sigma})$ whose length grows linearly.
\end{prop}
\begin{proof}

We will use the following convention throughout.  A vector field
$X\in \Gamma(T\Sigma)$ will also denote its $r$-independent lift to the cylinder
$(0,\infty)\times \Sigma$, namely the vector field satisfying
\[
[\partial_r,X]=0.
\]
In particular, $g_c(X,Y) = r^2 g_{\Sigma}(X,Y)$ on $\{r\}\times \Sigma$, so the length of these vectors grows linearly on the cone.
\\

Decompose
$$\tilde K=f\,\partial_r+Y,$$
where $f\in C^\infty(C_0)$ and $Y$ is tangent to the slices
$\Sigma_r=\{r=\mathrm{const}\}$.

For any vector fields $U,V$, the Lie derivative of the metric is
$$(\mathcal L_{\tilde K}g_c)(U,V)=g_c(\nabla_U\tilde K,V)+g_c(U,\nabla_V\tilde K).$$
Since $\tilde K$ is Killing, this vanishes for all $U,V$.

Now let $X,Z$ be vector fields on $\Sigma$, extended to $(0,\infty)\times\Sigma$ constant in $r$.
%\todo[inline]{F: We should say something about how to extend these vectorfields.}
 For the cone metric
$g_c=dr^2+r^2g_\Sigma,$ recall that
the Levi-Civita connection satisfies
$$\nabla_{\partial_r}\partial_r=0,\qquad
\nabla_{\partial_r}X=\nabla_X\partial_r=\frac1r X,\qquad \nabla_X Z=\nabla^\Sigma_X Z-r\,g_\Sigma(X,Z)\partial_r.$$

Taking $U=V=\partial_r$, we get
$$0=(\mathcal L_{\tilde K}g_c)(\partial_r,\partial_r)=2\,g_c(\nabla_{\partial_r}\tilde K,\partial_r)=2\,\partial_r f.$$

Taking $U=\partial_r$ and $V=X$, we get
$$0=(\mathcal L_{\tilde K}g_c)(\partial_r,X)=g_c(\nabla_{\partial_r}\tilde K,X)+g_c(\partial_r,\nabla_X\tilde K).$$
Using the connection formulas and $\partial_r f=0$, and writing
$Y(r,\cdot)$ as an $r$-dependent vector field on $\Sigma$, we have
$$\nabla_{\partial_r}Y=\partial_r Y+\frac1r Y,$$
where $\partial_rY$ denotes the coefficientwise derivative with respect to $r$.
Also,
$$g_c(\partial_r,\nabla_XY)=-g_c(\nabla_X\partial_r,Y)=-\frac1r\,g_c(X,Y).$$
Hence
\begin{align*}
0
&=(\mathcal L_{\tilde K}g_c)(\partial_r,X) \\
&=g_c(\nabla_{\partial_r}\tilde K,X)+g_c(\partial_r,\nabla_X\tilde K) \\
&=g_c\!\left(\partial_rY+\frac1rY,X\right)+df(X)-\frac1r\,g_c(X,Y) \\
&=g_c(\partial_rY,X) +df(X) \\
&=r^2 g_\Sigma(\partial_rY,X)+df(X).
\end{align*}
Since this holds for every $X$ tangent to $\Sigma$, we obtain
$$\partial_rY=-\frac1{r^2}\nabla^\Sigma f.$$
Because $\partial_r f=0$, integrating in $r$ gives
\begin{equation}\label{eq:Y(r)}
Y(r,\cdot)=Y_0+\frac1r\nabla^\Sigma f
\end{equation}
for some vector field $Y_0$ on $\Sigma$.

Finally, taking $U=X$ and $V=Z$, we obtain
$$0=(\mathcal L_{\tilde K}g_c)(X,Z)=g_c(\nabla_X\tilde K,Z)+g_c(X,\nabla_Z\tilde K).$$
Substituting \eqref{eq:Y(r)} and using the cone connection formulas yields
$$0=r^2(\mathcal L_{Y_0}g_\Sigma)(X,Z)+2r\bigl(\nabla_\Sigma^2 f+f\,g_\Sigma\bigr)(X,Z).$$
Since this identity holds for all $r>0$, we must have
$$\mathcal L_{Y_0}g_\Sigma=0,\qquad\nabla_\Sigma^2 f=-f\,g_\Sigma.$$

If $f\not\equiv 0$, then by Obata's theorem, \cite[Theorem A]{Obata}, $(\Sigma,g_\Sigma)$ is the unit and round sphere,
and therefore the cone $(C_0,g_c)$ is flat, contradicting the hypothesis.
Hence $f\equiv 0$.

Therefore
$$g_c(\tilde K,r\partial_r)=r\,g_c(\tilde K,\partial_r)=r f=0$$
on $C_0$, so $\tilde K$ is tangent to every link $\{r=\mathrm{const}\}$.
Moreover, \eqref{eq:Y(r)} reduces to $Y(r,\cdot)=Y_0$, with $Y_0$ a Killing
field on $(\Sigma,g_\Sigma)$.
\end{proof}

\section{A Gagliardo-Nirenberg inequality}\label{sec-gag-nir-ine}

 The following lemma is an ad-hoc Gagliardo-Nirenberg inequality. As the proof is short, we include it for the sake of completeness.  
 
 \begin{lemma}\label{lemma-inter-inequ-Appleton}
 Let $(M^n,g)$ be a complete Riemannian manifold such that here exists $c_0>0$ such that for all $x\in M$ and $r>0$, 
 \begin{equation}\label{eq:noncollapsing}
 \vol_{g}B_{g}(x_0,r)\geq c_0r^n.
 \end{equation}
Then there exists a constant $C=C(n,c_0)$ such that if $h$ is a smooth symmetric $2$-tensor with $h\in L^2(g)$ and such that $|\nabla^{g}h|_{g}\in L^{\infty}(g)$, the following interpolation inequality holds:
 \begin{eqnarray}
\|h\|_{L^{\infty}(g)}\leq C\|h\|_{L^2(g)}^{\frac{2}{n+2}}\|\nabla^{g}h\|^{\frac{n}{n+2}}_{L^{\infty}(g)}.\label{inter-inequ-Appleton}
\end{eqnarray}
In particular, if $(M^n,g)$ is asymptotically conical, then there exists $c_0>0$ such that \eqref{eq:noncollapsing} holds true, hence so does \eqref{inter-inequ-Appleton}.
\end{lemma}
 \begin{proof}
 Assume $h$ is not identically $0$ and let $x_0\in M$ such that $|h|_{g}(x_0)>0$. Then consider the geodesic ball $B_{g}(x_0,r_0)$ where $r_0:=|h|_{g}(x_0)/\|\nabla^{g}h\|_{L^{\infty}(g)}$,  and observe that
 \begin{eqnarray*}
|h|_{g}(x)\geq |h|_{g}(x_0)-\|\nabla^{g}h\|_{L^{\infty}(g)}d_{g}(x_0,x), \quad\text{for all $x\in B_{g}(x_0,r_0)$},
\end{eqnarray*}
where we used Kato inequality to bound $|\nabla^{g}|h||$ almost everywhere from above by $|\nabla^{g}h|_{g}$.
 
By integrating the previous inequality with respect to the variable $x$:
\begin{equation}
\begin{split}\label{series-est-easy}
\|h\|_{L^2(g)}^2&\geq \int_{B_{g}(x_0,r_0)}|h|^2_{g}d\mu_{g}\\
&\geq\int_{B_{g}(x_0,r_0)}\left(|h|_{g}(x_0)-\|\nabla^{g}h\|_{L^{\infty}(g)}d_{g}(x_0,x)\right)^2d\mu_{g}\\
&=\|\nabla^{g}h\|_{L^{\infty}(g)}^2\int_0^{r_0}\left(r_0-r\right)^2A(x_0,r)\,dr\\
&=\|\nabla^{g}h\|_{L^{\infty}(g)}^2\int_0^{r_0}\left(r_0-r\right)\vol_gB_g(x_0,r)\,dr,
\end{split}
\end{equation}
where we used the co-area formula in the penultimate inequality together with an integration by parts in the last line. Here $A(x_0,r)$ denotes the $(n-1)$-Hausdorff measure of $S(x_0,r)$.

The estimate \eqref{series-est-easy} together with the assumption on the volume of geodesic balls allows us to estimate as follows:
\begin{equation}
\begin{split}
\|h\|_{L^2(g)}^2&\geq c_0r_0^{n+2}\|\nabla^{g}h\|_{L^{\infty}(g)}^2\underbrace{\int_0^{1}\left(1-s\right)s^n\,ds}_{=:c_n}\\
&=c_0c_n|h(x_0)|_g^{n+2}\|\nabla^{g}h\|_{L^{\infty}(g)}^{-n}\int_0^{1}\left(1-s\right)s^n\,ds.
\end{split}
\end{equation}
%\todo[inline]{TO: should it be $=c_0c_nr_0^{n+2}|h(x_0)|_g^{n+2}\|\nabla^{g}h\|_{L^{\infty}(g)}^{-n}$?}
Therefore, one gets the expected inequality since the point $x_0$ is arbitrary.
 \\
 
Finally, if $(M^n,g)$ is asymptotically conical, we claim that  there exists $c_0>0$ such that for all $x_0\in M$ and $r>0$, $\vol_gB_g(x_0,r)\geq c_0\,r^n$.

Indeed, as $(M^n,g)$ is asymptotically conical, there exists $\iota_0>0$ such that $\inj_g(x_0)\geq \iota_0\, d_g(p,x_0)$ where $p\in M$ is fixed. In particular, there exists $c>0$ independent of $x_0$ such that $\vol_gB_g(x_0,r)\geq c\,r^n$ as long as $r< \iota_0\,d_g(p,x_0)$. Moreover, as $(M^n,g)$ is asymptotically conical again, there exists $c'>0$ such that $\vol_gB_g(p,r)\geq c'\,r^n$ for all $r>0$. The previous estimate ensures that if $r\geq 2d_g(p,x_0)$, $\vol_gB_g(x_0,r)\geq \vol_gB_g(p,r/2)\geq 2^{-n}c'\,r^n.$

Finally, if $r\in [\min\{\iota_0,2\}d_g(p,x_0),2d_g(p,x_0)]$ then, $$\vol_gB_g(x_0,r)\geq \vol_gB_g(x_0,\iota_0d_g(p,x_0))\geq c\,(\iota_0d_g(p,x_0))^n\geq \underbrace{c(\iota_0/2)^n}_{=:c_0}r^n,$$
which proves the remaining cases to consider in order to get the desired claim.
 \end{proof}

Since an (asymptotically conical) expanding gradient Ricci soliton $(M^n,g(t))_{t>0}$ is a self-similar solution to the Ricci flow, there exists $c_0>0$ such that for all positive time $t$, $x\in M$ and $r>0$, $\vol_{g(t)}B_{g(t)}(x_0,r)\geq c_0r^n$. 
 
 Consequently, by Lemma \ref{lemma-inter-inequ-Appleton}, one has the following uniform-in-time Gagliardo-Nirenberg inequality:
  \begin{coro}\label{coro-inter-inequ-Appleton}
 Let $(M^n,g(t))_{t>0}$ be an asymptotically conical expanding gradient Ricci soliton. Then there exists $C>0$ such that if $h$ is a smooth symmetric $2$-tensor with $h\in L^2(g(t))$ and such that $|\nabla^{g(t)}h|_{g(t)}\in L^{\infty}(g(t))$ for every $t>0$, the following interpolation inequality holds:
 \begin{eqnarray}
\|h\|_{L^{\infty}(g(t))}\leq C\|h\|_{L^2(g(t))}^{\frac{2}{n+2}}\|\nabla^{g(t)}h\|^{\frac{n}{n+2}}_{L^{\infty}(g(t))},\label{bis-inter-inequ-Appleton}\quad t>0.
\end{eqnarray}
In particular, if $h\in L^2_f\cap L^{\infty}(g)$ and $\nabla^{g}h\in L^{\infty}(g)$, then, for all $\varepsilon>0$,
\begin{equation*}
\sup_{M}e^{\frac{f}{(1+\varepsilon)(n+2)}}|h|_{g}<\infty.
\end{equation*}
\end{coro}

\begin{proof}
The proof of \eqref{bis-inter-inequ-Appleton} is straightforward on the basis of Lemma \ref{lemma-inter-inequ-Appleton}.

The second estimate follows by applying \eqref{bis-inter-inequ-Appleton} at $t=1$ to a $2$-tensor of the kind $\chi_{R,\varepsilon}\cdot h$ where $\chi_{R,\varepsilon}(x):=\chi_{\varepsilon}(f(x)/R)$, $R\geq 1$, and where $\chi_{\varepsilon}:\RR\rightarrow [0,1]$ is a smooth function such that $\chi_{\varepsilon}\equiv0$ on $(-\infty,1]$, $\chi_{\varepsilon}\equiv1$ on $[1+\varepsilon,+\infty)$ and $|\chi'_{\varepsilon}|\leq C/\varepsilon$ on $\RR$.
Indeed, one gets on the one hand,
\begin{equation*}
\begin{split}
\sup_{f\,\geq\, (1+\varepsilon)R}|h|_g&\leq \sup_{M}|\chi_{R,\varepsilon}\cdot h|_g\\
&\leq C\|\chi_{R,\varepsilon}\cdot h\|_{L^2(g)}^{\frac{2}{n+2}}\|\nabla^{g}(\chi_{R,\varepsilon}\cdot h)\|^{\frac{n}{n+2}}_{L^{\infty}(g)}\\
&\leq C_{\varepsilon}\|\chi_{R,\varepsilon}\cdot h\|_{L^2(g)}^{\frac{2}{n+2}},
\end{split}
\end{equation*}
since $h\in L^{\infty}(g)$ and $\nabla^gh\in L^{\infty}(g)$ and where $C_{\varepsilon}$ is a positive constant that depends on $\varepsilon$ and $\|h\|_{L^\infty(g)}$ and $\|\nabla h\|_{L^\infty(g)}$.

On the other hand, as $h\in L^2_f$,
\begin{equation*}
\begin{split}
\|\chi_{R,\varepsilon}\cdot h\|_{L^2(g)}^{\frac{2}{n+2}}&\leq e^{-\frac{R}{n+2}}\|\chi_{R,\varepsilon}\cdot h\|_{L^2_f}^{\frac{2}{n+2}}\leq e^{-\frac{R}{n+2}}\|h\|_{L^2_f}^{\frac{2}{n+2}}.
\end{split}
\end{equation*}
Combining the two previous estimates leads to: 
\begin{equation*}
\begin{split}
\sup_{f\,=\, (1+\varepsilon)R}e^{\frac{f}{(1+\varepsilon)(n+2)}}|h|_g&=e^{\frac{R}{n+2}}\sup_{f\,=\, (1+\varepsilon)R}|h|_g\\
&\leq e^{\frac{R}{n+2}}\sup_{f\,\geq\, (1+\varepsilon)R}|h|_g\leq C_{\varepsilon}\|h\|_{L^2_f}^{\frac{2}{n+2}}.
\end{split}
\end{equation*}
This estimate being true for an arbitrary $R>0$, the result follows.
\end{proof}

\section{Soliton identities}
We recall the following identities and properties of an expanding gradient Ricci soliton:
\begin{lemma}\label{id-sol-egs}
Let $(M^n,g,\nabla^g f)$ be a non trivial expanding gradient Ricci soliton. Then:
\begin{eqnarray}
&&\Delta_g f = \R_g+\frac{n}{2}>0, \label{equ:1} \\
&&\nabla^g \R_g+ 2\Ric(g)(\nabla^g f,\cdot)^{\sharp}=0, \label{equ:2} \\
&&\arrowvert \nabla^g f \arrowvert_g^2+\R_g=f+Cst, \label{equ:3}\\
&&\div_g\Rm(g)(Y,Z,U)=\Rm(g)(Y,Z, \nabla^g f,U),\label{equ:4}
\end{eqnarray}
for any vector fields $Y$, $Z$, $U$. 

In particular, $\nabla^gf$ is an eigenvector field of $\Delta_f$ associated to the eigenvalue $\frac{1}{2}$, i.e. on $M$, 
\begin{equation}
\Delta_f\nabla^gf-\frac{1}{2}\nabla^gf=0.\label{eigen-vec-nabla^gf}
\end{equation}

Moreover, the elliptic equations satisfied by the scalar curvature and the Ricci tensor are:
\begin{eqnarray}
%\Delta_f \Rm(g)&=&-\Rm(g)+\Rm(g)\ast\Rm(g),\label{ell-equ-curv-op}\\
\Delta_f\R_g&=&-\R_g-2|\Ric(g)|^2_g,\label{ell-equ-scal-curv}\\
\Delta_f\Ric(g)&=&-\Ric(g)-2\Rm(g)\ast\Ric(g).\label{ell-equ-ric-curv}
\end{eqnarray}
where $(\Rm(g)\ast \Ric(g))_{ij}:=\Rm(g)_{iklj}\Ric(g)^{kl}$.
%where, if $A$ and $B$ are two tensors, $A\ast B$ denotes any linear combination of contractions of the tensorial product of $A$ and $B$.
\end{lemma}
The proof of the above soliton identities can be found for instance in \cite[Chapter $4$]{Cho-Lu-Ni-Boo}. The proof of the strict positivity of $\Delta_gf$ is a consequence of the strong minimum principle applied to \eqref{ell-equ-scal-curv} as explained for instance \cite{Pig-Rim-Set}.
Similarly, one can prove:
\begin{lemma}\label{id-sol-es}
Let $(M^n,g,X)$ be an expanding Ricci soliton. Then:
\begin{eqnarray}
&&\div_gX = \R_g+\frac{n}{2}>0, \label{equ:1bis} \\
&&\Delta_gX+\Ric(g)(X,\cdot)^{\sharp}=0 \label{equ:2bis}.
\end{eqnarray}
\end{lemma}

\section{Extending Killing fields canonically}\label{sec-kill-ext}
In this section we explain how to extend Killing vector fields on the tangent cone at infinity of an expanding gradient Ricci soliton in a canonical way:
\begin{prop}\label{can-kill-vec-field}
Let $(M^n,g,\nabla^gf)$ be an asymptotically conical expanding gradient Ricci soliton with asymptotic cone $(C_0,g_c)$. Let $\tilde{K}$ be a Killing vector field orthogonal to $r\partial_r$. Then, there exists a smooth vector field $K$ and $\varepsilon>0$ on $M$ such that:
\begin{equation*}
\begin{split}
\Delta_{f}K-\frac{K}{2}&=0,\quad\textit{  on $M$,}\\
\sup_Mf^{\frac{k-1}{2}}|\nabla^{g,\,k}K|_{g}&<\infty,\quad k\geq 0,\\
L_{f}\left(\Li_Kg\right)&=0,\quad\textit{  on $M$,}\quad \left|\nabla^{g,\,k}\left( \Li_Kg\right)\right|_{g}=O(f^{\frac{k-n}{2}}e^{-f}),\quad k\geq 0,\\
\div_{f}K&=0,\quad\textit{  on $M$,}\\
\iota^*K&=\tilde{K}+O(r^{1-\varepsilon}),\quad\text{as $r\rightarrow+\infty$}.
\end{split}
\end{equation*}
\end{prop}

The proof of Proposition \ref{can-kill-vec-field} is essentially taken from \cite{Cho-Rot-Sym} which in turn is based on \cite{Bre-Rot} establishing Perelman's conjecture on the steady Bryant soliton. 

\begin{proof}
For the reasons above, we will be sketchy.

Let us consider a Killing vector field $\tilde{K}$ on the link of the cone. 
%The pull-back of a Killing field on the link of the cone. Let us denote this approximate Killing field by $U$.

First, one ensures the existence of a global vector field satisfying the same asymptotic properties as those of the vector field $\tilde{K}$. Its proof can be adapted from the proof of \cite[Proposition $3.2$]{Bre-Rot} adapted to our setting: there exists a smooth vector field $\overline{K}$ on $M$, 

\begin{equation}
\begin{split}
&|\nabla^{g,\,k}\overline{K}|_g=O(f^{-\frac{k-1}{2}}),\quad k\geq 0,\label{cov-der-U_i}\\
&\arrowvert\Li_{\overline{K}}(g)\arrowvert_g=o(1),\quad g(\overline{K},\nabla^g f)=\textit{o}\left(f\right),\\
&\left\arrowvert\nabla^{g,\,k}\left(\Delta_f \overline{K}-\frac{\overline{K}}{2}\right)\right\arrowvert_g=O(f^{-\frac{k}{2}-\varepsilon}),\quad \quad k\geq 0,
\end{split}
\end{equation}
as $f$ tends to $+\infty$ for some $\varepsilon>0$.

By lack of ellipticity of the equation satisfied by a Killing field, it turns out that it is easier to perturb such an almost Killing field with the help of the operator $\Delta_g+\Ric(g)$: recall that a Killing vector field lies in the kernel of this operator for an arbitrary Riemannian manifold. Since we are dealing with expanding gradient Ricci solitons, it is more adequate to estimate the kernel of the operator $\Delta_f -\Id/2$: one can show that there exists a smooth vector field $K$ on $M$ such that:
\begin{equation}
\begin{split}\label{prop-V}
&\Delta_fK-\frac{K}{2}=0,\quad |\nabla^{g,\,k}K|_g=O(f^{-\frac{k-1}{2}}),\quad k\geq 0,\\
&\arrowvert\Li_{K}(g)\arrowvert_g=o(1),\quad g(K,\nabla^gf)=\textit{o}\left(f\right),
\end{split}
\end{equation}
as $f$ tends to $+\infty$ and where $C$ is a positive constant. Finally, it can be shown that the vector fields $\overline{K}$ and $K$ are asymptotic to each other in the following sense: there exists $\varepsilon>0$ such that,
\begin{equation*}
|\nabla^{g,\,k}(K-\overline{K})|_g=O(f^{-\frac{k-1}{2}-\varepsilon}), \quad k\geq 0.\label{conv-rate-alm-kill-fie}
\end{equation*}
Thanks to Lemma \ref{lem-equiv-gauge-eqn} applied to $Y:=K$, one concludes that $h:=\Li_Kg$ lies in the kernel of $L_f$ and that thanks to \eqref{prop-V}, $|h|_g=o(1)$ as $f$ tends to $+\infty$. Invoking \cite{Der-Asy-Com-Egs}, one gets the improved decay $|h|_g=O(f^{-\frac{n}{2}}e^{-f})$ which implies in particular that $h\in L^2_f$. The corresponding decay on the covariant derivatives follow from Lemma \ref{lem-dec-eigen-lic}.
Recall the following Bochner-type formula going back to Lichnerowicz adapted here to the setting of expanding gradient Ricci solitons: for any symmetric $2$-tensor $h$ on $M$,
\begin{equation*}
\div_{f}(L_{f}h)=\Delta_{f}(\div_{f}h)-\frac{1}{2}\div_{f}h.
\end{equation*}
In particular, as $L_{f}\left(\Li_Kg\right)=0$, then, $\Delta_{f}\alpha=\frac{1}{2}\alpha$ where $\alpha:=\div_{f}\Li_Kg$ is a $1$-form. Thanks to Lemma \ref{lem-dec-eigen-lic}, $|\nabla^{g}\Li_Kg|_g=O(f^{\frac{1-n}{2}}e^{-f})$ so that $|\alpha|_{g}=o(1)$ as $f\rightarrow +\infty$. The use of the maximum principle allows us to conclude that $\alpha\equiv 0$. Now, Lemma \ref{lem-bochner-lie-der} gives us that $\div_fK$ is constant on $M$. An application of Lemma \ref{lem-equiv-gauge-eqn} applied to $Z:=\div_f\Li_Kg=0$ gives that $0=2\Delta_f\div_fK=\div_fK$ as required.
\end{proof}

\section{Linear growth of the soliton vector field}\label{sec-lin-growth}

Quadratic decay of the Ricci curvature implies the required linear growth of the soliton vector field.

\begin{prop}\label{prop-lin-growth}
Let $(M^n,g)$ be a complete Riemannian manifold, $p\in M$, $K\geq0$, and $X$ a smooth vector field on $M$ such that, on $M$,
\begin{equation}\label{hyp-lin-growth}
\Ric(g)+\frac{g}{2}=\frac{1}{2}\Li_Xg
\qquad\text{and}\qquad
|\Ric(g)|_g\leq\frac{K}{(1+d_g(p,\cdot))^{2}}.
\end{equation}
Then there exists $C>0$, depending only on $K$ and on $|X|_g$ restricted to $\overline{B_g(p,R)}$, where $R:=1+4(|X|_g(p)+K)$, such that $|X|_g\leq C(1+d_g(p,\cdot))$ on $M$.
\end{prop}

\begin{proof}
Set $r:=d_g(p,\cdot)$ and $A:=|X|_g(p)+K$, so that $R=1+4A$, and let $\Phi_t$ denote the local flow of $X$.

Let $\gamma$ be a unit speed minimizing geodesic from $p$ to $y\neq p$. As in the proof of \cite[Lemma 2.2]{Che-Der}, where $X=\nabla^gf$, \eqref{hyp-lin-growth} gives
\[
\frac{d}{ds}g(X,\dot\gamma)
=\frac{1}{2}\Li_Xg(\dot\gamma,\dot\gamma)
=\frac{1}{2}+\Ric(g)(\dot\gamma,\dot\gamma).
\]
Since $r(\gamma(s))=s$, integration gives
\begin{equation}\label{radial-lin-growth}
\begin{aligned}
g\big(X(y),\dot\gamma(r(y))\big)
&\geq\frac{r(y)}{2}-|X|_g(p)
-\int_0^{+\infty}\frac{K}{(1+s)^2}\,ds\\
&=\frac{r(y)}{2}-A.
\end{aligned}
\end{equation}
In particular, $|X|_g\geq r/2-A>0$ on $\{r\geq R\}$.

Let $x\in M$ and $\rho(t):=r(\Phi_t(x))$. The function $\rho$ is locally Lipschitz. At a differentiability time $t$ with $\rho(t)>0$, let $\gamma$ be a unit speed minimizing geodesic from $p$ to $\Phi_t(x)$. For $0<s<\rho(t)$, the functions $s+d_g(\gamma(s),\cdot)$ are smooth upper supports for $r$ at $\Phi_t(x)$, so their derivatives along the flow agree with $\rho'(t)$. Thus \eqref{radial-lin-growth} gives
\begin{equation}\label{dist-lin-growth}
\rho'(t)
=g\big(X(\Phi_t(x)),\dot\gamma(\rho(t))\big)
\geq\frac{\rho(t)}{2}-A
\end{equation}
almost everywhere where $\rho>0$.

Suppose that $r(x)>R$. Whenever $\rho\geq R$, \eqref{dist-lin-growth} gives $\rho'\geq R/2-A=1/2+A$. Consequently, the backward orbit remains in the compact ball $\overline{B_g(p,r(x))}$ until it reaches $\{r=R\}$. By the Escape Lemma \cite[Lemma 9.19]{Lee2} and this differential inequality, it reaches this level set at some finite time $-T$, with $\rho\geq R$ on $[-T,0]$.

The function $u(t):=|X|_g(\Phi_t(x))$ is positive on $[-T,0]$ by \eqref{radial-lin-growth}. As in the proof of \cite[Lemma 2.9]{Bam-Che}, the soliton equation gives
\[
(u^2)'=\Li_Xg(X,X)=u^2+2\Ric(g)(X,X),
\]
and hence
\[
(\log u)'=\frac{1}{2}+\frac{\Ric(g)(X,X)}{u^2}
\leq\frac{1}{2}+\frac{K}{(1+\rho)^2}.
\]
Since $K\leq A$ and $0<\rho-4A\leq(1+\rho)^2$ on $[-T,0]$, \eqref{dist-lin-growth} yields
\begin{equation}\label{exp-lin-growth}
\begin{aligned}
\left(\log\frac{u}{\rho-4A}\right)'
&\leq
\frac{1}{2}+\frac{K}{(1+\rho)^2}
-\frac{\rho/2-A}{\rho-4A}\\
&=\frac{K}{(1+\rho)^2}-\frac{A}{\rho-4A}
\leq0
\end{aligned}
\end{equation}
almost everywhere. Since $\rho(-T)=R$ and $R-4A=1$, it follows that
\[
|X|_g(x)
\leq
|X|_g(\Phi_{-T}(x))\frac{r(x)-4A}{R-4A}
\leq
r(x)\sup_{\overline{B_g(p,R)}}|X|_g.
\]
Together with the bound on $\overline{B_g(p,R)}$, this proves the proposition with
\[
C:=1+\sup_{\overline{B_g(p,R)}}|X|_g.
\]
\end{proof}

\begin{coro}\label{cor-lin-growth}
Let $(M^n,g_1,X_1)$ be an expanding Ricci soliton as in item~\eqref{item-non-grad} of Theorem~\ref{theo-main}, with $\varepsilon=\varepsilon(n,g_0)$ decreased if necessary. Then there exists $\chi_1>0$ such that $|X_1|_{g_0}\leq\chi_1f_0^{1/2}$ on $M$. In particular, $X_1$ is a complete vector field.
\end{coro}

\begin{proof}
Let $p$ be a critical point of $f_0$. For $\varepsilon$ sufficiently small in terms of $n$ and $g_0$, Proposition~\ref{prop-alm-sharp-rate} at $\tau=0$ and Lemma~\ref{app-C-monster} give
\[
\frac{1}{2}g_0\leq g_1\leq2g_0
\qquad\text{and}\qquad
|\Ric(g_1)|_{g_1}
\leq\frac{K}{(1+d_{g_1}(p,\cdot))^2}
\]
for some $K=K(n,g_0)$, since $g_0$ has quadratic curvature decay. In particular, $g_1$ is complete. Proposition~\ref{prop-lin-growth} and the metric comparison give the first inequality below, while \eqref{basic-pot-fct-f_0} and $f_0>0$ give the second:
\[
|X_1|_{g_0}
\leq C(1+d_{g_0}(p,\cdot))
\leq\chi_1f_0^{1/2}
\]
for some $C,\chi_1>0$.

The fact that a vector field $X_1$ with at most linear growth is complete is standard.
\end{proof}

%\newpage
%\bibliographystyle{amsalpha}

%\bibliography{ref-WPIC1}

\def\cprime{$'$} \def\cprime{$'$}
\providecommand{\bysame}{\leavevmode\hbox to3em{\hrulefill}\thinspace}
\providecommand{\MR}{\relax\ifhmode\unskip\space\fi MR }
\providecommand{\MRhref}[2]{%
  \href{http://www.ams.org/mathscinet-getitem?mr=#1}{#2}
}
\providecommand{\href}[2]{#2}

\end{document}